%% file: WKO2.tex
\documentclass[12pt,reqno]{amsart}
\def\draft{n}
\usepackage{
  fullpage,graphicx,dbnsymb,amsmath,amssymb,floatflt,picins,
  multicol,mathtools,stmaryrd,xcolor,tensor
}

\usepackage[breaklinks=true]{hyperref}\hypersetup{colorlinks,
   linkcolor={green!50!black},
   citecolor={green!50!black},
   urlcolor=blue
}

\usepackage{xr}
\usepackage[all]{xy}

\if\draft y
  \usepackage[color]{showkeys}
  \usepackage{everypage,draftwatermark}
  \SetWatermarkScale{4}
  \SetWatermarkLightness{0.9}
\fi

\input macros.tex

\input defs.tex

\begin{document}
\newdimen\captionwidth\captionwidth=\hsize
\setcounter{secnumdepth}{4}

\pagestyle{empty}
\thispagestyle{empty}

\begin{center}
  {\Large\bf Corrigendum to ``Finite Type Invariants of w-Knotted Objects II''}

  \vskip 5mm By Dror Bar-Natan and Zsuzsanna Dancso

  \vskip 3mm Feb.~29,~2024

  \vskip 3mm (a revised version of the paper is attached below this corrigendum)

\end{center}

\vskip 5mm
\def\figsuffix{pstex_t}
\input Corrigendum-Text.tex

\newpage
\setcounter{page}{1}

\pagestyle{plain}

\title[Finite Type Invariants of w-Knotted Objects]{Finite Type Invariants
  of w-Knotted Objects II (V2): Tangles, Foams and the Kashiwara-Vergne Problem}

\author{Dror~Bar-Natan}
\address{
  Department of Mathematics\\
  University of Toronto\\
  Toronto Ontario M5S 2E4\\
  Canada
}
\email{drorbn@math.toronto.edu}
\urladdr{http://www.math.toronto.edu/~drorbn}

\author{Zsuzsanna Dancso}
\address{
  School of Mathematics and Statistics\\
  The University of Sydney\\
  Carslaw Building \\
  Camperdown, NSW, 2006, Australia
}
\email{zsuzsanna.dancso@sydney.edu.au}

\date{first edition May 5, 2014, this edition Feb.~29,~2024. Published in Mathematische Annalen
{\bf 367} (2017) 1517--1586. See
\url{http://www.math.toronto.edu/drorbn/LOP.html\#WKO2};
the \arXiv{1405.1955} edition may be older}


\thanks{This work was partially supported by NSERC grant RGPIN 262178.
  This paper is part 2 of a 4-part series whose first two parts originally
  appeared as a combined preprint,~\cite{WKO}.
}

\begin{abstract}
  \input abstract.tex

\end{abstract}

\maketitle

{\small \tableofcontents}

\input intro.tex
\input alg.tex

\input tangles.tex
\input foams.tex
\input odds.tex
\input glossary.tex
\input refs.tex

\end{document}

%% file: macros.tex
\theoremstyle{plain}
\newtheorem{theorem}{Theorem}[section]

\newtheorem{proposition}[theorem]{Proposition}
\newtheorem{fact}[theorem]{Fact}
\newtheorem{lemma}[theorem]{Lemma}

\theoremstyle{definition}
\newtheorem{definition}[theorem]{Definition}

\theoremstyle{remark}
\newtheorem{comment}[theorem]{Comment}

\newtheorem{example}[theorem]{Example}

\newtheorem{remark}[theorem]{Remark}
\newtheorem{warning}[theorem]{Warning}

\newlength{\standardunitlength}
\catcode`\@=11
\long\def\@makecaption#1#2{%
    \vskip 10pt
    \setbox\@tempboxa\hbox{
      \small\sf{\bfcaptionfont #1. }\ignorespaces #2}%
    \ifdim \wd\@tempboxa >\captionwidth {%
        \rightskip=\@captionmargin\leftskip=\@captionmargin
        \unhbox\@tempboxa\par}%
      \else
        \hbox to\hsize{\hfil\box\@tempboxa\hfil}%
    \fi}
\font\bfcaptionfont=cmssbx10 scaled \magstephalf
\newdimen\@captionmargin\@captionmargin=2\parindent
\newdimen\captionwidth\captionwidth=\hsize
\catcode`\@=12

\newcommand{\gr}{\operatorname{gr}}
\newcommand{\im}{\operatorname{im}}

\newcommand{\tr}{\operatorname{tr}}

\def\qed{{\hfill\text{$\Box$}}}

\newlength{\globalparindent}
\def\arXiv#1{{\href{http://arXiv.org/abs/#1}{arXiv:#1}}}

%% file: defs.tex
\def\bbN{{\mathbb N}}
\def\bbQ{{\mathbb Q}}
\def\bbR{{\mathbb R}}
\def\bbZ{{\mathbb Z}}
\def\calA{{\mathcal A}}
\def\calB{{\mathcal B}}

\def\calC{{\mathcal C}}
\def\calD{{\mathcal D}}

\def\calI{{\mathcal I}}

\def\calO{{\mathcal O}}
\def\calP{{\mathcal P}}
\def\calR{{\mathcal R}}
\def\calS{{\mathcal S}}
\def\calT{{\mathcal T}}
\def\calU{{\mathcal U}}
\def\fraka{{\mathfrak a}}

\def\ol{\overline}

\def\fil{\operatorname{fil}\,}
\def\gr{\operatorname{gr}\,}

\def\grad{{\operatorname{grad}\,}}
\def\grads{{\operatorname{grad}}}
\newcommand{\Ass}{\operatorname{Ass}}
\def\divop{\operatorname{div}}
\def\CA{\operatorname{CA}}

\def\attr{\operatorname{\mathfrak{tr}}}
\def\lie{\operatorname{\mathfrak{lie}}}
\def\der{\operatorname{\mathfrak{der}}}
\def\sder{\operatorname{\mathfrak{sder}}}
\def\tder{\operatorname{\mathfrak{tder}}}
\def\TAut{\operatorname{TAut}}
\def\SAut{\operatorname{SAut}}

\def\vD{{\mathit v\!D}}
\def\vT{{\mathit v\!T}}

\def\uT{{\mathit u\!T}}
\def\wT{{\mathit w\!T}}
\def\wTF{{\mathit w\!T\!F}}
\def\wTFo{{\mathit w\!T\!F^o}}
\def\sKTG{{\mathit s\!K\!T\!G}}

\def\aft{{\overrightarrow{4T}}}
\def\aAS{{\overrightarrow{AS}}}
\def\aSTU{{\overrightarrow{STU}}}
\def\aIHX{{\overrightarrow{IHX}}}
\def\Rs{{R$1^{\!s}$}}

\def\up{{\uparrow}}

\def\posvertex{\raisebox{-1mm}{\input{figs/PlusVertex.pstex_t}}}

\def\negvertex{\raisebox{-1mm}{\input{figs/MinusVertex.pstex_t}}}

\def\pstex#1{\begin{array}{c}
  \input figs/#1.pstex_t
  \if\draft y
    \smash{\makebox[0in]{\color{labelkey}{\tt figs/#1}}}
  \fi
\end{array}}

\def\draftcut{\if\draft y \newpage \fi}

\def\inverted#1{{\rotatebox[origin=c]{180}{#1}}}

\def\glos#1{\setlength{\fboxsep}{0pt}\colorbox{yellow}{$#1$}}
\def\glost#1{\setlength{\fboxsep}{0pt}\colorbox{yellow}{#1}}

%% file: figs/PlusVertex.pstex_t
\begin{picture}(0,0)%
\includegraphics{figs/PlusVertex.pstex}%
\end{picture}%
%
%
\setlength{\unitlength}{3158sp}%
\begingroup\makeatletter\ifx\SetFigFont\undefined%
\gdef\SetFigFont#1#2#3#4#5{%
  \reset@font\fontsize{#1}{#2pt}%
  \fontfamily{#3}\fontseries{#4}\fontshape{#5}%
  \selectfont}%
\fi\endgroup%
\begin{picture}(211,211)(3439,-1160)
\end{picture}%

%% file: figs/MinusVertex.pstex_t
\begin{picture}(0,0)%
\includegraphics{figs/MinusVertex.pstex}%
\end{picture}%
%
%
\setlength{\unitlength}{3158sp}%
\begingroup\makeatletter\ifx\SetFigFont\undefined%
\gdef\SetFigFont#1#2#3#4#5{%
  \reset@font\fontsize{#1}{#2pt}%
  \fontfamily{#3}\fontseries{#4}\fontshape{#5}%
  \selectfont}%
\fi\endgroup%
\begin{picture}(211,211)(3065,-1347)
\end{picture}%

%% file: Corrigendum-Text.tex
In Section 4 of ~\cite{WKO2} we introduce and study {\em w-tangled
foams}.  These are defined combinatorially, as a finitely generated {\em
circuit algebra} with certain Reidemeister relations. For motivation, we
present a local topological interpretation of w-tangled foams as tangled
tubes in $\bbR^4$.  Unfortunately, some errors have occurred in this
interpretation, stemming from our lack of care around ``2D orientations'',
as below. We wish to thank Yusuke Kuno and Haruko Miyazawa for noting
these issues.

\begin{itemize}
\parpic[r]{\input{figs/R4Issue.\figsuffix}}
\item If foams are to have 2D orientations, they cannot be glued to
each other unless these orientations match.  This breaks the circuit
algebra structure. This can be corrected by switching to {\em coloured}
circuit algebras, in which strands and strand-ends are coloured by
their 2D orientations, and gluings are allowed only if colours match. However, we
choose a different resolution explained below.
\item We missed that the ``$A_e$'' operations of flipping only the 1D but not the 2D
orientations of a strand interact unpleasantly with R4 moves, as shown
on the right. This can be corrected by giving more care to the colours (2D orientations)
of the edges next to a vertex and making the precise meaning of
the R4 move depend on these colours. Since this is notationally tedious, we choose 
a different resolution below.
\end{itemize}

Although the resolution is a small change to the content of the paper, it necessitates some notation and language changes, and for the convenience of our readers we have incorporated these into the body of the paper as a revised version, available at~\cite{WKO2v2}.

\parpic[r]{\input{figs/WenShort.\figsuffix}}
The easier solution to both of the problems above is to forgo of
2D orientations and of the operations $A_e$ altogether, and to
replace them with more serious attention to the {\em wens}, which are cut open Klein-bottles 
(\cite[Section~4.5]{WKO2}, and \cite[Section~4.1]{WKO2v2}). Indeed, conjugating a 
strand by a wen and reversing it (strand reversal is the operation $S_e$, unchanged from the original 
paper) has the same effect on crossings as the operation $A_e$ 
was meant to have. 

The main result of the paper is identifying {\em homomorphic expansions} for w-tangled foams with 
solutions to the Kashiwara-Vergne (KV) equations. A crucial equation relating the values of such an expansion, the {\em unitarity equation} $V\cdot
A_1A_2(V)=1$ \cite[Equation~(12)]{WKO2} is replaced by an equivalent equation,
$V^*V=1$ \cite[Equation~(U), page 41]{WKO2v2} in which the ``global adjoint'' $V^*$ is obtained from $V$
by multiplying it by a wen on all three sides, and reversing all strand orientations. 

In \cite{WKO2} we claimed that the value of the wen under an expansion must be 1, and that expansions for w-tangled foams are in bijection with Kashiwara-Vergne solutions with {\em even Duflo function}. In fact, the value of the wen could be non-trivial \cite[Lemma 4.9]{WKO2v2}, but the statement regarding KV solutions remains true for homomorphic expansions where the value of the wen is set to be $1$. 

Theorem~4.9 of \cite{WKO2} stated that for ``orientable w-tangled foams'', i.e. where wens are not included, homomorphic expansions are in one to one correspondence with KV solutions in general. This theorem remains true: in \cite[Section 4.6]{WKO2v2} we present a more careful definition of orientable w-tangled foams as a sub-circuit algebra. Edge-wise adjoint operations replaced by the well-defined global adjoint (global orientation reversal and composition with wens at every tangle end) as above, noting that these wens all cancel. Homomorphic expansions for these foams are indeed in one to one correspondence with KV solutions.

We note that the results of papers \cite{BDS, DHR} which build on Section 4 of \cite{WKO2} are unaffected by this Corrigendum: in these papers one can simply mechanically replace $A_1A_2(V)$ with $V^*$.

\section*{Conflict of interest statement}
On behalf of all authors, the corresponding author states that there is no conflict of interest.

%% file: figs/WenShort.pstex_t
\begin{picture}(0,0)%
\includegraphics{figs/WenShort.pstex}%
\end{picture}%
%
%
\setlength{\unitlength}{3947sp}%
\begingroup\makeatletter\ifx\SetFigFont\undefined%
\gdef\SetFigFont#1#2#3#4#5{%
  \reset@font\fontsize{#1}{#2pt}%
  \fontfamily{#3}\fontseries{#4}\fontshape{#5}%
  \selectfont}%
\fi\endgroup%
\begin{picture}(2348,1299)(53,-448)
\put(226, 14){\makebox(0,0)[lb]{\smash{{\SetFigFont{12}{14.4}{\rmdefault}{\mddefault}{\itdefault}{\color[rgb]{0,0,0}wen}%
}}}}
\end{picture}%

%% file: abstract.tex
This is the second in a series of papers dedicated to studying 
w-knots, and more generally, w-knotted objects (w-braids, w-tangles,
etc.). These are classes of knotted objects that are \underline{w}ider but
\underline{w}eaker than their ``\underline{u}sual'' counterparts. To
get (say) w-knots from usual knots (or u-knots), one has to allow non-planar ``virtual''
knot diagrams, hence enlarging the the base set of knots. But then one
imposes a new relation beyond the ordinary collection of Reidemeister moves,
called the ``overcrossings commute'' relation, making w-knotted
objects a bit weaker once again.
Satoh~\cite{Satoh:RibbonTorusKnots} studied several classes of
w-knotted objects (under the name ``\underline{w}eakly-virtual'') and has
shown them to be closely related to certain classes of knotted surfaces
in $\bbR^4$. 

In this article we study finite type invariants of w-tangles and
w-trivalent graphs (also referred to as w-tangled foams).
Much as the spaces $\calA$ of chord diagrams for ordinary knotted
objects are related to metrized Lie algebras, the spaces
$\calA^w$ of ``arrow diagrams'' for w-knotted objects are related to
not-necessarily-metrized Lie algebras. Many questions concerning w-knotted
objects turn out to be equivalent to questions about Lie algebras. Most
notably we find that a homomorphic universal finite type invariant of
w-foams is essentially the same as a solution of
the Kashiwara-Vergne~\cite{KashiwaraVergne:Conjecture} conjecture and
much of the Alekseev-Torossian~\cite{AlekseevTorossian:KashiwaraVergne}
work on Drinfel'd associators and Kashiwara-Vergne can be re-interpreted
as a study of w-foams.

%% file: intro.tex
\draftcut
\section{Introduction} \label{sec:intro}
This is the second in a series of papers on w-knotted objects.  In the
first paper \cite{Bar-NatanDancso:WKO1}, we took a classical approach
to studying finite type invariants of w-braids and w-knots and proved
that the universal finite type invariant for w-knots is essentially the
Alexander polynomial. In this paper we will study finite type invariants
of w-tangles and w-tangled foams from a more algebraic point of view,
and prove that ``homomorphic'' universal finite type invariants of
w-tangled foams are in one-to-one correspondence with solutions to
the (Alekseev-Torossian version of) the Kashiwara-Vergne problem in
Lie theory.  Mathematically, this paper does not depend on the results
of \cite{Bar-NatanDancso:WKO1} in any significant way, and the reader
familiar with the theory of finite type invariants will have no difficulty
reading this paper without having read \cite{Bar-NatanDancso:WKO1}.
However, since this paper starts with an abstract re-phrasing of the
well-known finite type story in terms of general algebraic structures,
readers who need an introduction to finite type invariants may find
it more pleasant to read \cite{Bar-NatanDancso:WKO1} first (especially
Sections \ref{1-sec:intro}, \ref{1-sec:w-braids} and
\ref{1-subsec:VirtualKnots}--\ref{1-subsec:LieAlgebras}).

\subsection{Motivation and hopes}
This article and its siblings \cite{Bar-NatanDancso:WKO1} and \cite{Bar-NatanDancso:WKO3} are
efforts towards a larger goal. Namely, we believe many of the difficult
algebraic equations in mathematics, especially those that are
written in graded spaces, more especially those that are related in
one way or another to quantum groups~\cite{Drinfeld:QuantumGroups},
and to the work of Etingof and
Kazhdan~\cite{EtingofKazhdan:BialgebrasI}, can be understood, and indeed
would appear more natural, in terms of finite type invariants of various
topological objects.

This work was inspired by Alekseev and Torossian's results \cite{AlekseevTorossian:KashiwaraVergne} on
Drinfel'd associators and the Kashiwara-Vergne conjecture, both of which fall into the aforementioned
class of ``difficult equations in graded spaces''. The Kashiwara-Vergne conjecture --- 
proposed in 1978 \cite{KashiwaraVergne:Conjecture} and proven in 2006 by Alekseev and Meinrenken \cite{AlekseevMeinrenken:KV} --- 
has strong implications in Lie theory and harmonic analysis, and is a cousin of the Duflo isomorphism, which was shown to be
knot-theoretic in~\cite{Bar-NatanLeThurston:TwoApplications}.
We also know that Drinfel'd's theory
of associators~\cite{Drinfeld:QuasiHopf} can be
interpreted as a theory of well-behaved universal finite type
invariants of parenthesized tangles\footnote{``$q$-tangles''
in~\cite{LeMurakami:Universal}, ``non-associative tangles''
in~\cite{Bar-Natan:NAT}.}~\cite{LeMurakami:Universal, Bar-Natan:NAT},
or of knotted trivalent graphs~\cite{Dancso:KIforKTG}.

In Section \ref{sec:w-foams} we will re-interpret the Kashiwara-Vergne
conjecture as the problem of finding a ``homomorphic'' universal finite type invariant of
a class of w-knotted trivalent graphs (more accurately named w-tangled foams).
This result fits into a bigger picture incorporating usual, virtual and w-knotted objects and
their theories of finite type invariants, connected by the inclusion map from usual
to virtual, and the projection from virtual to w-knotted objects. In a sense that will be made precise in
Section \ref{sec:generalities}, usual and w-knotted objects with this mapping form
a unified algebraic structure, and the relationship between Drinfel'd associators
and the Kashiwara-Vergne conjecture is explained as a theory of finite type invariants for
this larger structure. This will be the topic of Section \ref{subsec:KTG}.

We are optimistic that this paper is a step towards re-interpreting the
work of Etingof and Kazhdan~\cite{EtingofKazhdan:BialgebrasI}
on quantization of Lie bi-algebras as a
construction of a well-behaved universal finite type invariant of
virtual knots~\cite{Kauffman:VirtualKnotTheory, Kuperberg:VirtualLink} or of a similar class
of virtually knotted objects. 
However, w-knotted objects are quite interesting in their own right,
both topologically and algebraically:
they are related to combinatorial group
theory, to groups of movies of flying rings in $\bbR^3$, and more generally, to
certain classes of knotted surfaces in $\bbR^4$. The references include
\cite{BrendleHatcher:RingsAndWickets, FennRimanyiRourke:BraidPermutation,
Goldsmith:MotionGroups, McCool:BasisConjugating, Satoh:RibbonTorusKnots}.

In \cite{Bar-NatanDancso:WKO1} we studied the universal finite type invariants of w-braids and w-knots,  
the latter of which turns out to be essentially the Alexander polynomial. 
A more thorough introduction about our ``hopes and dreams'' and the u-v-w big picture can also be
found in \cite{Bar-NatanDancso:WKO1}.

\subsection{A brief overview and large-scale explanation}
We are going to start by developing the algebraic ingredients of the paper in Section \ref{sec:generalities}.
The general notion of an {\em algebraic structure} lets us treat spaces of a topological or diagrammatic nature
in a unified algebraic manner. All of braids, w-braids, w-knots, w-tangles, etc., and their 
associated chord- or arrow-diagrammatic counterparts form algebraic structures, and so do any number of these spaces combined, with 
maps between them.

We then introduce {\em associated graded structures} with respect to a specific filtration, the machine which in our case takes an algebraic structure of 
``topological nature'' (say, braids with $n$ strands) and produces the corresponding
diagrammatic space (for braids, horizontal chord diagrams on $n$ vertical strands). This is done by taking the associated graded space 
with respect to a given filtration, namely the powers of the augmentation ideal in the algebraic structure.

An {\em expansion}, sometimes called a universal finite type invariant, is a map from an algebraic structure (in this case one of topological nature)
to its associated graded (a structure of combinatorial/diagrammatic nature), with a certain universality property. 
A {\em homomorphic expansion} is one that is in addition
``well behaved'' with respect to the {\em operations} of the algebraic structure (such as composition and strand doubling for braids, for example).

The three main results of the paper are as follows:
\begin{enumerate}

\item
As mentioned before, our goal is to provide a topological framework for the Kashiwara-Vergne (KV) problem.
One of our main results is Theorem~\ref{thm:ATEquivalence}, in which we establish a bijection between certain homomorphic expansions
of {\em w-tangled foams} (introduced in Section \ref{sec:w-foams}) and solutions of the Kashiwara-Vergne equations. More precisely, ``certain'' 
homomorphic expansions means ones that are group-like (a commonly used condition), and subject to another minor technical condition.
Section \ref{sec:w-tangles} leads up to this result by studying the simpler case of {\em w-tangles} and identifying building blocks
of its associated graded structure as the spaces which appear in the \cite{AlekseevTorossian:KashiwaraVergne} formulation of the KV equations.

\item
In Theorem \ref{thm:WenATEquivalence} we study an unoriented version of w-tangled foams, and prove that homomorphic expansions
for this space (group-like and subject to the some minor conditions) are in one-to-one correspondence with solutions to the KV problem
with {\em even Duflo function}. This sets the stage for perhaps the most interesting result of the paper:

\item Section \ref{subsec:wTFcompatibility} marries the theory above with the theory of ordinary (not w-) knotted trivalent graphs (KTGs).
For technical reasons explained in Section \ref{sec:w-foams}, we work with a signed version of KTGs (sKTG).
Roughly speaking, homomorphic expansions for sKTGs are determined by a {\em Drinfel'd associator}. Furthermore, sKTGs map naturally
into w-tangled foams. 

\noindent
In Theorem \ref{thm:ZuwCompatible} we prove that any homomorphic expansion of sKTGs coming from a {\em horizontal chord} associator
has a compatible homomorphic expansion of w-tangled foams, and furthermore, these expansions are in one-to-one correspondence
with {\em symmetric} solutions of the KV problem. This gives a topological explanation for the relationship between
Drinfel'd associators and the KV conjecture.
\end{enumerate}

We note that in \cite{Bar-NatanDancso:WKO3}
we'll further capitalize on these insights to provide a topological proof and interpretation for 
Alekseev, Enriquez and Torossian's explicit solutions
for the KV conjecture in terms of associators \cite{AlekseevEnriquezTorossian:ExplicitSolutions}.

Several of the structures of a topological nature
in this paper (w-tangles and w-foams) are introduced as {\em Reidemeister theories}. That is, the spaces are built from pictorial 
generators (such as crossings) which can be connected arbitrarily, and the resulting pictures are then factored out by certain relations
(``Reidemeister moves''). Technically speaking, this is done using the framework of {\em circuit algebras} (similar to planar algebras but 
without the planarity requirement) which are introduced in Section \ref{sec:generalities}. 

One of the fundamental theorems of classical knot theory is Reidemeister's theorem, which states that isotopy classes of knots are
in bijection with {\em knot diagrams} modulo Reidemeister moves. In our case, w-knotted objects have a Reidemeister description and
a topological interpretation in terms of ribbon knotted tubes in $\bbR^4$. However, the analogue of the Reidemeister theorem, i.e. the 
statement that these
two interpretations coincide, is only known for w-braids \cite{McCool:BasisConjugating, Dahm:GeneralBraid, BrendleHatcher:RingsAndWickets}.

For w-tangles and w-foams (and w-knots as well) there is a map $\delta$ from the Reidemeister presentation to the appropriate class of ribbon 
2-knotted objects in $\bbR^4$. In our case this means that all the generators 
have a local topological interpretation and the relations represent isotopies. The map $\delta$ is certainly a surjection,
but it is only conjectured to be injective (in other words, it is possible that some relations are missing). 

The main difficulty in proving the injectivity of $\delta$ lies
in the management of the ribbon structure. A ribbon 2-knot is a knotted sphere or long tube in $\bbR^4$ which admits a filling with only certain 
types of singularities. While there are Reidemeister theorems for general 2-knots in $\bbR^4$ \cite{CarterSaito:KnottedSurfaces},
the techniques don't translate well to ribbon 2-knots, mainly because it is not well understood how different ribbon structures (fillings) of
the same ribbon 2-knot can be obtained from each other through Reidemeister type moves. The completion of such a theorem would be of great interest.
We suspect that even if $\delta$ is not injective, the present set of generators and relations describes a set of ribbon-knotted
tubes in $\bbR^4$ with possibly some extra combinatorial information, similarly to how, say, dropping the $R1$ relation in classical knot theory
results in a Reidemeister theory for framed knots with rotation numbers.

The paper is organized as follows: we start with a discussion of general algebraic structures,
associated graded structures, expansions (universal finite type invariants) and ``circuit algebras'' in Section \ref{sec:generalities}.  
In Section \ref{sec:w-tangles} we study w-tangles and identify some of the spaces \cite{AlekseevTorossian:KashiwaraVergne}
where the KV conjecture ``lives'' as the spaces of ``arrow diagrams'' (the w-analogue of chord diagrams) 
for certain w-tangles. In Section \ref{sec:w-foams} we study w-tangled foams and we prove the main theorems discussed above.
For more detailed information
consult the ``Section Summary'' paragraphs at the beginning
of each of the sections. A glossary of notation is on
page~\pageref{sec:glossary}.

\subsection{Corrigendum}\label{subsec:Corr} Yusuke Kuno and Haruko Miyazawa pointed out an error in the topological interpretation of w-tangled foams. Following several in-depth conversations with them, we made improvements to Section~\ref{sec:w-foams}, and submitted a Corrigendum \cite{Bar-NatanDancso:WKO2Corr}. In short, the ``adjoint'' $A_e$ operations of the published version of this paper were not well-defined with respect to the Reidemeister 4 relation. The simplest fix is to omit surface (two-dimensional) orientations and the adjoint operations entirely. The adjoint operations were used to relate the two ``vertex'' generators of w-foams, however, this can be achieved just as well by composing one of the vertices by ``wens'' (cut open Klein bottles). In this version of the paper we incorporated these changes into Section~\ref{sec:w-foams} for easier reading, and made several other stylistic improvements to the exposition.

\def\summaryalg{In this section we introduce the associated graded structure 
of an ``arbitrary algebraic structure'' with respect to powers of its augmentation ideal 
(Sections~\ref{subsec:AlgebraicStructures} and \ref{subsec:Grad})
and introduce the notions of ``expansions'' and ``homomorphic expansions'' (\ref{subsec:Expansions}). 
Everything is so general that practically anything is an example,
yet our main goal is to set the language for the examples of w-tangles
and w-tangled foams, which appear later in this paper. Both of
these examples are types of ``circuit algebras'', and hence we
end this section with a general discussion of circuit algebras
(Sec.~\ref{subsec:CircuitAlgebras}).}

\def\summarytangles{In Sec.~\ref{subsec:vw-tangles} we introduce
v-tangles and w-tangles, the v- and w- counterparts of the
standard knot-theoretic notion of ``tangles'', and briefly discuss their
finite type invariants and their associated spaces of ``arrow diagrams'',
$\calA^v(\uparrow_n)$ and $\calA^w(\uparrow_n)$. We then construct a
homomorphic expansion $Z$, or a ``well-behaved'' universal finite type
invariant for w-tangles. The only algebraic tool we need to
use is $\exp(a):=\sum a^n/n!$ (Sec.~\ref{subsec:vw-tangles}
is in fact a routine extension of parts of
\cite[Section~\ref{1-sec:w-knots}]{Bar-NatanDancso:WKO1}).
In Sec.~\ref{subsec:ATSpaces} we show that
$\calA^w(\uparrow_n)\cong\calU(\fraka_n\oplus\tder_n\ltimes\attr_n)$,
where $\fraka_n$ is an Abelian algebra of rank $n$ and where
$\tder_n$ and $\attr_n$, two of the primary spaces used by Alekseev
and Torossian~\cite{AlekseevTorossian:KashiwaraVergne}, have simple
descriptions in terms of cyclic words and free Lie algebras. We also show
that some functionals studied in~\cite{AlekseevTorossian:KashiwaraVergne},
$\divop$ and $j$, have natural interpretations in our language.
In Section~\ref{subsec:sder} we discuss a subclass of w-tangles called ``special''
w-tangles, and relate them by similar means to Alekseev and Torossian's
$\sder_n$ and to ``tree level'' ordinary Vassiliev theory. A local topological interpretation is described in Sec.~\ref{subsec:TangleTopology} and the
uniqueness of $Z$ is studied in Sec.~\ref{subsec:UniquenessForTangles}.}

\def\summaryfoams{In this section we add ``foam vertices'' to w-tangles (and a few lesser
things as well) and ask the same questions we asked before;
primarily, ``is there a homomorphic expansion?''. As we shall
see, in the current context this question is equivalent to the
Alekseev-Torossian~\cite{AlekseevTorossian:KashiwaraVergne} version
of the Kashiwara-Vergne~\cite{KashiwaraVergne:Conjecture} problem and
explains the relationship between these topics and Drinfel'd's theory
of associators. This section was corrected and improved post-publication: see Section~\ref{subsec:Corr} and \cite{Bar-NatanDancso:WKO2Corr}.}

\subsection{Acknowledgement} We thank Yusuke Kuno and Haruko Miyazawa for their significant contribution which helped us improve Section 4 post-publication. We wish to thank Anton Alekseev, Jana
Archibald, Scott Carter, Karene Chu, Iva Halacheva, Joel Kamnitzer,
Lou Kauffman, Peter Lee, Louis Leung, Jean-Baptiste Meilhan, Dylan Thurston, Lucy Zhang
and the anonymous referees for comments and suggestions.

%% file: alg.tex
\draftcut
\section{Algebraic Structures, Expansions, and Circuit Algebras}
\label{sec:generalities}

\begin{quote} \small {\bf Section Summary. }
  \summaryalg
\end{quote}

\subsection{Algebraic Structures} \label{subsec:AlgebraicStructures}

An ``algebraic structure'' $\glos{\calO}$ is some collection $(\calO_\alpha)$
of sets of objects of different kinds, where the subscript
$\alpha$ denotes the ``kind'' of the objects in $\calO_\alpha$,
along with some collection of ``operations'' $\glos{\psi_\beta}$, where
each $\psi_\beta$ is an arbitrary map with domain some product
$\calO_{\alpha_1}\times\dots\times\calO_{\alpha_k}$ of sets of objects,
and range a single set $\calO_{\alpha_0}$ (so operations may be unary or
binary or multinary, but they always return a value of some fixed kind).
We also allow some named ``constants'' within some $\calO_\alpha$'s
(or equivalently, allow some 0-nary operations).\footnote{
Alternatively define ``algebraic structures'' using the theory of
``multicategories''~\cite{Leinster:Higher}. Using this language,
an algebraic structure is simply a functor from some ``structure''
multicategory $\calC$ into the multicategory {\bf Set} (or into {\bf
Vect}, if all $\calO_i$ are vector spaces and all operations are
multi-linear). A ``morphism'' between two algebraic structures over the
same multicategory $\calC$ is a natural transformation between the two
functors representing those structures.} The operations may or may not
be subject to axioms --- an ``axiom'' is an identity asserting that
some composition of operations is equal to some other composition of
operations.

\begin{figure}[h]
\parbox[m]{7cm}{\caption{An algebraic structure $\calO$ with 4 kinds of
objects and one binary, 3 unary and two 0-nary operations (the constants
$1$ and $\sigma$).}
\label{fig:AlgebraicStructure}}
\parbox[m]{9cm}{\centering
\def\ObjTypeThree{{\raisebox{3mm}{\parbox{0.7in}{\footnotesize
  $\left\{\parbox{0.5in}{\centering objects of kind 3}\right\}=$
}}}}
\input figs/AlgebraicStructure.pstex_t%
}
\end{figure}

Figure~\ref{fig:AlgebraicStructure} illustrates the general notion of an
algebraic structure. Here are a few specific examples:
\begin{itemize}
\item We will use $\langle b\rangle$, the free group on one generator $b$, as
  a running example throughout this chapter (of course $\langle b \rangle$ is isomorphic to $\bbZ$). 
  This is an algebraic structure
  with one kind of objects, a binary operation ``multiplication'',
  a unary operation ``inverse'', one constant ``the
  identity'', and the expected axioms.
\item Groups in general: one kind of objects, one binary ``multiplication'',
  one unary ``inverse'', one constant ``the identity'', and some axioms.
\item Group homomorphisms: Two kinds of objects, one for each
  group. 7 operations --- 3 for each of the two groups and the homomorphism
  itself, going between the two groups. Many axioms.
\item A group acting on a set, a group extension, a split group extension
  and many other examples from group theory.
\item A quandle is a set with an operation $\glos{\uparrow}$, satisfying
  $(x \uparrow y)\uparrow z=(x\uparrow y)\uparrow(y\uparrow z)$ and some further minor axioms. This is an algebraic 
  structure with one kind of objects and one operation. See \cite{WKO} for an analysis of quandles from the perspective of this paper.
\item Planar algebras as in~\cite{Jones:PlanarAlgebrasI} and circuit
  algebras as in Section~\ref{subsec:CircuitAlgebras}.
\item The algebra of knotted trivalent graphs as
  in~\cite{Bar-Natan:AKT-CFA, Dancso:KIforKTG}.
\item Let $\varsigma\colon B\to S$ be an arbitrary homomorphism of groups (though
  our notation suggests what we have in mind --- $B$ may well be braids,
  and $S$ may well be permutations). We can consider an algebraic structure
  $\calO$ whose kinds are the elements of $S$, for which the objects of
  kind $s\in S$ are the elements of $\calO_s:=\varsigma^{-1}(s)$, and with
  the product in $B$ defining operations
  $\calO_{s_1}\times\calO_{s_2}\to\calO_{s_1s_2}$.
\item W-tangles and w-foams, studied in the following two sections of this paper.
\item Clearly, many more examples appear throughout mathematics.
\end{itemize}

\draftcut
\subsection{Associated Graded Structures} \label{subsec:Grad}

Any algebraic structure $\calO$ has an ``especially natural'' 
associated graded structure: that is, we take the associated structure with
respect to a specific and natural filtration. This will be a repeating construction 
throughout the rest of this paper series. 

First extend
$\calO$ to allow formal linear combinations of objects of the same kind
(extending the operations in a linear or multi-linear manner), then let
$\glos{\calI}$, the ``augmentation ideal'', be the sub-structure made out of
all such combinations in which the sum of coefficients is $0$, then let
$\calI^m$ be the set of all outputs of algebraic expressions (that is,
arbitrary compositions of the operations in $\calO$) that have at least
$m$ inputs in $\calI$ (and possibly, further inputs in $\calO$), and
finally, set
\begin{equation} \label{eq:gradO}
  \glos{\grad}\calO:=\bigoplus_{m\geq 0} \calI^m/\calI^{m+1}.
\end{equation}
Clearly, with the operations inherited from $\calO$, the associated graded
$\grad\calO$ is again algebraic structure with the same multi-graph
of spaces and operations, but with new objects and with new operations
that may or may not satisfy the axioms satisfied by the operations of
$\calO$. The main new feature in $\grad\calO$ is that it is a ``graded''
structure; we denote the degree $m$ piece $\calI^m/\calI^{m+1}$ of
$\grad\calO$ by $\grads_m\calO$.

We believe that many of the most interesting graded structures that appear
in mathematics are the result of this construction (i.e., as associated graded 
structures with respect to powers of the augmentation ideal), and
that many of the interesting graded equations that appear in mathematics
arise when one tries to find ``expansions'', or ``universal finite type
invariants'', which are also morphisms\footnote{Indeed, if $\calO$ is
finitely presented then finding such a morphism $Z\colon \calO\to\grad\calO$
amounts to finding its values on the generators of $\calO$, subject to
the relations of $\calO$. Thus it is equivalent to solving a system
of equations written in some graded spaces.} $Z\colon \calO\to\grad\calO$
(see Section~\ref{subsec:Expansions}) or when one studies
``automorphisms'' of such expansions\footnote{The Drinfel'd graded
Grothendieck-Teichmuller group $\mathit{GRT}$ is an example of such an
automorphism group. See~\cite{Drinfeld:GalQQ, Bar-Natan:Associators}.}.
Indeed, the paper you are reading now is really the study of the
associated graded structures of various algebraic structures associated with
w-knotted objects. We would like to believe that much of the theory of
quantum groups (at ``generic'' $\hbar$) will eventually be shown to be a
study of the associatead graded structures of various algebraic structures associated
with v-knotted objects.

\begin{example}\label{ex:bbZ}
 We compute the associated graded structuture of the running example $\langle b \rangle$.
 Allowing formal $\bbQ$-linear combinations of elements we get $\bbQ\langle b\rangle=\bbQ[b,b^{-1}]$. 
 The augmentation ideal $\calI$ is generated by differences 
 $(b^n-1)$ as a vector space (where $1=b^0$), and generated by $(b-1)$ as an ideal. 
 
 We claim that $\grad \langle b \rangle \cong \bbQ[[c]]$, the algebra of power series
 in one variable. To show this, consider the map $\pi: \bbQ[[c]] \to \grad \langle b \rangle$ by setting $\pi(c)=[b-1]$ (mod $\calI^2$). 
 It is easy to show explicitly
 that $\pi$ is surjective. For example, in degree 1, we need to show that $b-1$ generates $\calI/\calI^2$. 
 indeed, $(b^n-1)-n(b-1)$ has a double zero at $b=1$, and hence $f=\frac{(b^n-1)-n(b-1)}{(b-1)^2}$ is a polynomial, 
 and $b^n-1=n(b-1)+f(b-1)^2$. So modulo $(b-1)^2 \in \calI^2$, $b^n-1=n(b-1)$. 
 A similar argument works to show that $(b-1)^k$ generates
 $\calI^k/\calI^{k+1}$. 
 
 Note that $\langle b \rangle$ can also be thought of as the pure braid group on two strands: $b$ would be a ``full twist'' and $c$ can be
 represented as a single ``horizontal chord''. In other knot theoretic settings,
 it is generally relatively easy to find a ``candidate associated graded'' and a map $\pi$, which can be shown to be surjective
 by explicit means. 
 
 To show that $\pi$ is injective we are going to use the machinery of ``expansions'' which is the tool we use to accomplish
 similar tasks in the later sections of this paper.
\end{example}

We end this section with two more examples of computing associated graded structures: the 
proof of Proposition \ref{prop:GradGrp} is an exercise; for the
proof of Proposition \ref{prop:GradQ} see \cite{WKO}.

\begin{proposition}\label{prop:GradGrp} If $G$ is a group, $\grad G$ is a graded associative
algebra with unit. Similarly, the associated graded structure of a group homomorphism is a homomorphism of
graded associative algebras. \qed
\end{proposition}

\begin{proposition} \label{prop:GradQ} If $Q$ is a unital quandle,
$\grads_0 Q$ is one-dimensional and $\grads_{>0} Q$ is a graded 
right Leibniz 
algebra\footnote{A Leibniz algebra is a Lie algebra without anti-commutativity, as
defined by Loday in \cite{Loday:LeibnizAlg}.}
generated by $\grads_1 Q$.
\end{proposition}

\draftcut
\subsection{Expansions and Homomorphic Expansions}
\label{subsec:Expansions}
We
start with the definition. Given an algebraic structure $\calO$ let
$\glos{\fil}\calO$ denote the filtered structure of linear combinations of
objects in $\calO$ (respecting kinds), filtered by the powers $(\calI^m)$
of the augmentation ideal $\calI$. Recall also that any graded space
$G=\bigoplus_mG_m$ is automatically filtered, by $\left(\bigoplus_{n\geq
m}G_n\right)_{m=0}^\infty$.

\begin{definition} An ``expansion'' $Z$ for $\calO$
is a map $Z\colon \calO\to\grad\calO$ that preserves the kinds of objects
and whose linear extension (also called $Z$) to $\fil\calO$ respects the
filtration of both sides, and for which $\left(\gr Z\right):
\left(\gr\fil\calO=\grad\calO\right) \to
\left(\gr\grad\calO=\grad\calO\right)$ is the identity map of
$\grad\calO$; we refer to this as the ``universality property''.
\end{definition}

In practical terms, this is equivalent to saying that $Z$ is a map
$\calO\to\grad\calO$ whose restriction to $\calI^m$ vanishes in degrees
less than $m$ (in $\grad\calO$) and whose degree $m$ piece is the
projection $\calI^m\to\calI^m/\calI^{m+1}$.

We come now to what is perhaps the most crucial definition in this paper.

\begin{definition} A ``homomorphic expansion'' is an expansion which
also commutes with all the algebraic operations defined on the algebraic
structure $\calO$.
\end{definition}

\noindent{\bf Why Bother with Homomorphic Expansions?} Primarily, for two
reasons:
\begin{itemize}
\item Often $\grad\calO$ is simpler to work with than $\calO$;
for one, it is graded and so it allows for finite ``degree by degree''
computations, whereas often times, such as in many topological examples,
anything in $\calO$ is inherently infinite. Thus it can be beneficial to
translate questions about $\calO$ to questions about $\grad\calO$. A
simplistic example would be, ``is some element $a\in\calO$ the square
(relative to some fixed operation) of an element $b\in\calO$?''. Well, if
$Z$ is a homomorphic expansion and by a finite computation it can be shown
that $Z(a)$ is not a square already in degree $7$ in $\grad\calO$, then
we've given a conclusive negative answer to the example question. Some less
simplistic and more relevant examples appear in~\cite{Bar-Natan:AKT-CFA}.
\item Often $\grad\calO$ is ``finitely presented'', meaning that it
is generated by some finitely many elements $g_1,\dots,g_k\in\calO$,
subject to some relations $R_1\dots R_{n}$ that can be written in terms
of $g_1,\dots,g_k$ and the operations of $\calO$. In this case, finding a
homomorphic expansion $Z$ is essentially equivalent to guessing the values
of $Z$ on $g_1,\dots,g_k$, in such a manner that these values
$Z(g_1),\dots,Z(g_k)$ would satisfy the $\grad\calO$ versions of the
relations $R_1\dots R_{n}$. So finding $Z$ amounts to solving equations in
graded spaces. It is often the case (as will be demonstrated in this paper;
see also~\cite{Bar-Natan:NAT, Bar-Natan:Associators})
that these equations are very interesting for their own algebraic sake, and
that viewing such equations as arising from an attempt to solve a problem
about $\calO$ sheds further light on their meaning. 
\end{itemize}

In practice, often the first difficulty in searching for an
expansion (or a homomorphic expansion) $Z\colon \calO\to\grad\calO$ is that its
would-be target space $\grad\calO$ is hard to identify. It is typically
easy to make a suggestion $\calA$ for what $\grad\calO$ could be. It
is typically easy to come up with a reasonable generating set $\calD_m$
for $\calI^m$ (keep some knot theoretic examples in mind, or $\bbZ$ in Example~\ref{ex:bbZ}). 
It is a bit harder but not
exceedingly difficult to discover some relations $\calR$ satisfied by the
elements of the image of $\calD$ in $\calI^m/\calI^{m+1}$ (4T, $\aft$, and
more in knot theory, there are no relations for $\bbZ$).
Thus we set $\calA:=\calD/\calR$; but it is often very hard to be
sure that we found everything that ought to go in $\calR$; so perhaps
our suggestion $\calA$ is still too big? Finding 4T for example 
was actually not {\em that} easy. Could we have
missed some further relations that are hiding in $\calA$?

The notion of an $\calA$-expansion, defined below, solves two problems at
once. Once we find an $\calA$-expansion we know that we've identified
$\grad\calO$ correctly, and we automatically get what we really wanted, a
($\grad\calO$)-valued expansion.

\parpic[r]{\raisebox{-9mm}{$\xymatrix{
  & \calA \ar@<-2pt>[d]_\pi \\
  \calO \ar[ur]^{Z_{\calA}} \ar[r]_<>(0.4)Z
  & \grad\calO \ar@<-2pt>[u]_{\gr Z_\calA}
}$}}
\begin{definition} \label{def:CanGrad}
A ``candidate assoctaed graded structure'' for an algebraic structure
$\calO$ is a graded structure $\glos{\calA}$ with the same operations as
$\calO$ along with a homomorphic surjective graded map $\pi\colon
\calA\to\grad\calO$. An ``$\calA$-expansion'' is a kind and filtration
respecting map $\glos{Z_\calA}\colon \calO\to\calA$ for which $(\gr
Z_\calA)\circ\pi\colon \calA\to\calA$ is the identity. One can similarly
define ``homomorphic $\calA$-expansions''.
\end{definition}

\begin{proposition} \label{prop:CanGrad}
If $\calA$ is a candidate associated graded of $\calO$
and $Z_\calA\colon \calO\to\calA$ is a homomorphic $\calA$-expansion, then
$\pi:\calA\to\grad\calO$ is an isomorphism and $Z:=\pi\circ Z_\calA$ is a
homomorphic expansion. (Often in this case, $\calA$ is identified with
$\grad\calO$ and $Z_\calA$ is identified with $Z$).
\end{proposition}

\begin{proof} Note that $\pi$ is surjective by birth. Since $(\gr Z_\calA)\circ\pi$
is the identity, $\pi$ it is also injective and hence it is an
isomorphism. The rest is immediate. \qed
\end{proof}

\begin{example}
 Back to $\langle b \rangle$, in Example~\ref{ex:bbZ} we found a candidate associated graded structure $\calA=\bbQ[[c]]$ and a map $\pi: c \mapsto [b-1]$.
 According to Proposition~\ref{prop:CanGrad}, it is enough to find a homomorphic $\calA$-expansion, that is, an algebra homomorphism
 $Z_{\calA}: \bbQ\langle b \rangle \to \bbQ[[c]]$ such that $\gr Z_\calA \circ \pi$ is the identity of $\bbQ[[c]]$. It is a straightforward calculation to
 check that any algebra map defined by $Z_{\calA}(b)=1+c+\{\text{higher order terms}\}$ satisfies this property. If one seeks a ``group-like''
 homomorphic expansion then $Z_{\calA}(b)=e^c$ is the only solution. In either case, exhibiting $Z_{\calA}$ proves that $\pi$ is injective and
 hence $\calA$ is the associated graded structure of $\langle b \rangle$.
\end{example}

\draftcut
\subsection{Circuit Algebras} \label{subsec:CircuitAlgebras}

``Circuit algebras'' are so common and everyday, and they make
such a useful language (definitely for the purposes of this paper,
but also elsewhere), we find it hard to believe they haven't made
it into the standard mathematical vocabulary\footnote{Or have they,
and we have been looking the wrong way?}. People familiar with planar
algebras~\cite{Jones:PlanarAlgebrasI} may note that circuit algebras are
just the same as planar algebras, except with the planarity requirement
dropped from the ``connection diagrams'' (and all colourings are dropped
as well). 

In our context, the main utility of circuit algebras is that they allow
for a much simpler presentation of $v$(irtual)- and $w$-tangles. 
There are planar algebra presentations of $v$- and $w$-tangles, generated by the usual crossings and the ``virtual crossing'', 
modulo the usual as well as the ``virtual'' and ``mixed'' Reidemeister moves.
Switching from planar algebras to circuit algebras however renders the extra generators
and relations unnecessary: the ``virtual crossing'' becomes merely a circuit algebra
artifact, and the new Reidemeister moves are implied by the circuit algebra structure 
(see Warning~\ref{warn:virtualxings}, Definition~\ref{def:vw-tangles}, and Remark \ref{rmk:VirtualXings}). 

The everyday intuition for circuit algebras comes from electronic circuits,
whose components can be wired together in many, not necessarily planar, 
ways, and it is not important to know how these wires are embedded in space.
For details and more motivation see Section~\ref{subsec:CAMotivation}.
We start formalizing this image by defining
``wiring diagrams'', the abstract analogs of printed circuit boards.
Let $\bbN$ denote the set of natural numbers including $0$, and for
$n\in\bbN$ let $\underline{n}$ denote some fixed set with $n$ elements,
say $\{1,2,\dots,n\}$.

{
\makeatletter\def\thm@space@setup{%
  \thm@preskip=0cm plus 0cm minus 0cm 
}\makeatother
\parpic[r]{\input{figs/WiringDiagram.pstex_t}}
\begin{definition} Let $k, n, n_1,\dots,n_k\in\bbN$ be natural
numbers. A ``wiring diagram'' $D$ with inputs $\underline{n_1},\dots
\underline{n_k}$ and outputs $\underline{n}$ is an unoriented
compact 1-manifold whose boundary is $\underline{n}\amalg
\underline{n_1}\amalg\cdots\amalg\underline{n_k}$, regarded
up to homeomorphism (on the right is an example with $k=3$, $n=6$,
and $n_1=n_2=n_3=4$). In strictly combinatorial terms, it is a
pairing\footnote{We mean ``pairing'' in the sense of combinatorics,
not in the sense of linear algebra. That is, an involution without fixed
point.} of the elements of the set $\underline{n}\amalg
\underline{n_1}\amalg\cdots\amalg\underline{n_k}$ along
with a single further natural number that counts closed
circles. If $D_1;\dots;D_m$ are wiring diagrams with inputs
$\underline{n_{11}},\dots,\underline{n_{1k_1}}; \dots;
\underline{n_{m1}},\dots,\underline{n_{mk_m}}$ and outputs
$\underline{n_1};\dots;\underline{n_m}$ and $D$ is a wiring diagram
with inputs $\underline{n_1};\dots;\underline{n_m}$ and outputs
$\underline{n}$, there is an obvious ``composition'' $D(D_1,\dots,D_m)$
(obtained by gluing the corresponding 1-manifolds, and also describable
in completely combinatorial terms) which is a wiring diagram with
inputs $(\underline{n_{ij}})_{1\leq i\leq k_j,1\leq j\leq m}$ and
outputs $\underline{n}$ (note that closed circles may be created in
$D(D_1,\dots,D_m)$ even if none existed in $D$ and in $D_1;\dots;D_m$).
\end{definition}
}

A circuit algebra is an algebraic structure (in the sense of
Section~\ref{subsec:Grad}) whose operations are parametrized by
wiring diagrams. Here's a formal definition:

\begin{definition} A circuit algebra consists of the following data:
\begin{itemize}
\item For every natural number $n\geq 0$ a set (or a $\bbZ$-module) $C_n$ ``of
circuits with $n$ legs''.
\item For any wiring diagram $D$ with inputs $\underline{n_1},\dots
\underline{n_k}$ and outputs $\underline{n}$, an operation (denoted by the
same letter) $D\colon C_{n_1}\times\dots\times C_{n_k}\to C_n$ (or linear
$D\colon C_{n_1}\otimes\dots\otimes C_{n_k}\to C_n$ if we work with
$\bbZ$-modules).
\end{itemize}
We insist that the obvious ``identity'' wiring diagrams with
$\underline{n}$ inputs and $\underline{n}$ outputs act as the identity of
$C_n$, and that the actions of wiring diagrams be compatible in the obvious
sense with the composition operation on wiring diagrams.
\end{definition}

A silly but useful example of a circuit algebra is the circuit algebra
$\glos{\calS}$ of empty circuits, or in our context, of ``skeletons''. The
circuits with $n$ legs for $\calS$ are wiring diagrams with $n$ outputs
and no inputs; namely, they are 1-manifolds with boundary $\underline{n}$
(so $n$ must be even).

More generally one may pick some collection of ``basic components''
(analogous to logic gates and junctions for electronic circuits as in
Figure~\ref{fig:FlipFlop}) and speak of the ``free circuit algebra''
generated by these components; even more generally we can speak of
circuit algebras given in terms of ``generators and relations''. (In the
case of electronics, our relations may include the likes of De Morgan's
law $\neg(p\vee q)=(\neg p)\wedge(\neg q)$ and the laws governing the
placement of resistors in parallel or in series.) We feel there is no need
to present the details here, yet many examples of circuit algebras given
in terms of generators and relations appear in this paper, starting with
the next section. We will use the notation $C=\CA\langle \, G \mid R \, \rangle$
to denote the circuit algebra generated by a collection of elements
$G$ subject to some collection $R$ of relations. 

People familiar with electric circuits know that connectors
sometimes come in ``male'' and ``female'' versions, and that you
can't plug a USB cable into a headphone jack. 
Thus one may define ``directed circuit algebras'' in which
the wiring diagrams are oriented, the circuit sets $C_n$ get replaced
by $C_{p,q}$ for ``circuits with $p$ incoming wires and $q$
outgoing wires'' and only orientation preserving connections are ever
allowed\footnote{By convention we label the boundary points of such circuits $1,\ldots,p+q$, with the first $p$
labels reserved for the incoming wires and the last $q$ for the outgoing. The inputs of wiring
diagrams must be labeled in the opposite way for the numberings to match.}. Likewise
there is a ``coloured'' version of everything, in which the wires may be
coloured by the elements of some given set $X$ (which may include among
its members the elements ``USB'' and ``audio'') and in which connections
are allowed only if the colour coding is respected. We will leave the
formal definitions of directed and coloured circuit algebras, as well
as the definitions of directed and coloured analogues of the skeletons
algebra $\calS$ and generators and relations for directed and coloured
algebras, as an exercise.

Note that there is an obvious notion of ``a morphism between
two circuit algebras'' and that circuit algebras (directed or not,
coloured or not) form a category. We feel that a precise definition
is not needed. A lovely example is the ``implementation morphism''
of logic circuits in the style of Figure~\ref{fig:FlipFlop} in Section \ref{sec:odds} into more
basic circuits made of transistors and resistors.

Perhaps the prime mathematical example of a circuit algebra is tensor
algebra. If $t_1$ is an element (a ``circuit'') in some tensor product of
vector spaces and their duals, and $t_2$ is the same except in a possibly
different tensor product of vector spaces and their duals, then once an
appropriate pairing $D$ (a ``wiring diagram'') of the relevant vector
spaces is chosen, $t_1$ and $t_2$ can be contracted (``wired together'')
to make a new tensor $D(t_1,t_2)$. The pairing $D$ must pair a vector
space with its own dual, and so this circuit algebra is coloured by the
set of vector spaces involved, and directed, by declaring (say) that
some vector spaces are of one gender and their duals are of the other. We
have in fact encountered this circuit algebra in
\cite[Section~\ref{1-subsec:LieAlgebras}]{Bar-NatanDancso:WKO1}.

Let $G$ be a group. A $G$-graded algebra $A$ is a collection $\{A_g\colon g\in
G\}$ of vector spaces, along with products $A_g\otimes A_h\to A_{gh}$ that
induce an overall structure of an algebra on $A:=\bigoplus_{g\in G}A_g$. In
a similar vein, we define the notion of an $\calS$-graded circuit algebra:

\begin{definition}\label{def:Skeleta} An $\calS$-graded circuit algebra, 
or a ``circuit algebra with skeletons'', is an algebraic structure $C$ with
spaces $C_\beta$, one for each element $\beta$ of the circuit algebra of
skeletons $\calS$, along with composition operations
$D_{\beta_1,\dots,\beta_k}\colon C_{\beta_1}\times\dots\times C_{\beta_k}\to
C_\beta$, defined whenever $D$ is a wiring diagram and
$\beta=D(\beta_1,\dots,\beta_k)$, so that with the obvious induced
structure, $\coprod_\beta C_\beta$ is a circuit algebra. A similar
definition can be made if/when the skeletons are
taken to be directed or coloured.
\end{definition}

Loosely speaking, a circuit algebra with skeletons is a circuit
algebra in which every element $T$ has a well-defined skeleton
$\varsigma(T)\in\calS$. Yet note that as an algebraic structure a circuit
algebra with skeletons has more ``spaces'' than an ordinary circuit
algebra, for its spaces are enumerated by skeleta and not merely by
integers. The prime examples for circuit algebras with skeletons appear in
the next section.

%% file: tangles.tex
\draftcut
\section{w-Tangles} \label{sec:w-tangles}

\begin{quote} \small {\bf Section Summary. }
  \summarytangles
\end{quote}

\subsection{v-Tangles and w-Tangles} \label{subsec:vw-tangles} Building on Section~\ref{subsec:CircuitAlgebras}, we define v-tangles and w-tangles combinatorially
as finitely presented circuit algebras, given by generators and relations. This facilitates the algebraic treatment we present, by which we connect w-tangle theory to Kashiwara-Vergne theory. In Section~\ref{subsec:TangleTopology} we recall Satoh's tubing map \cite[Sec.~3.1.1]{Bar-NatanDancso:WKO1}, to endow w-tangles with the expected topological meaning. 

\begin{definition}\label{def:vtanglediagrams}
The ($\calS$-graded) circuit algebra $\vD$ of v-tangle diagrams is
the $\calS$-graded directed circuit algebra freely generated by
two generators in $C_{2,2}$ called the {\em positive crossing},
$\tensor*[_1^4]{\text{\large$\overcrossing$}}{_2^3}$,
and the {\em negative crossing},
$\tensor*[_1^4]{\text{\large$\undercrossing$}}{_2^3}$. In as much
as possible we suppress the leg-numebering below; with this in
mind, $\vD:=$\raisebox{-1.5mm}{\input{figs/vDDef.pstex_t}}.
The skeleton of both crossings is the element
$\tensor*[_1^4]{\text{\large$\virtualcrossing$}}{_2^3}$
(the pairing of 1\&3 and 2\&4) in $\calS_{2,2}$. That is,
$\varsigma(\overcrossing)=\varsigma(\undercrossing)=\virtualcrossing$.
\end{definition}

\begin{figure}[b]
 \input figs/wTangleExample.pstex_t
 \caption{$V \in \vD_{3,3}$ is a $v$-tangle diagram. $V$ is the result of applying the circuit algebra operation 
 $D: C_{2,2} \times C_{2,2} \times C_{2,2} \to C_{3,3}$,
 given by the wiring diagram shown, acting on two negative crossings and one positive crossing. In other words $V=D(\undercrossing, \undercrossing, \overcrossing)$.
 The skeleton of $V$ is given by $\varsigma(V)=D(\virtualcrossing, \virtualcrossing, \virtualcrossing)$, which is equal in $\calS$ to the diagram shown here. 
 Note that we usually suppress the circuit algebra numbering of boundary points. Note also that the apparent ``virtual crossings'' of $V$ are not 
 virtual crossings but merely part of the circuit algebra structure, see Warning~\ref{warn:virtualxings}. The same is true for the  
 crossings appearing in the skeleton $\varsigma(V)$.}
 \label{fig:wTangleExample}
\end{figure}

\begin{example}\label{ex:vDiag}
 An example of a v-tangle diagram $V$ is shown the left side of Figure \ref{fig:wTangleExample}. $V$ is a circuit algebra composition of two negative crossings
 and one positive crossing by the wiring diagram $D$, as shown. The right side of the same figure shows the skeleton $\varsigma(V)$ of $V$: 
 to produce the skeleton, replace each crossing by the element
 $\virtualcrossing$ in $\calS$ and apply the same wiring diagram. The elements of $\calS$ are oriented 1-manifolds with numbered boundary points, 
 and hence the result is equal to the one shown in the figure.
\end{example}

\begin{warning}\label{warn:virtualxings}
 People familiar with the planar presentation of virtual tangles may be accustomed to 
 the notion of there being another type of crossing: the ``virtual crossing''. The main point
 of introducing circuit algebras (as opposed to working with planar algebras) is to eliminate the
 need for virtual crossings: they become part of the CA structure. This greatly simplifies the presentation of both $v$- and $w$-tangles:
 there is one less generator, as seen above, and far fewer relations, as we explain in Remark~\ref{rmk:VirtualXings}.
\end{warning}

\begin{definition}\label{def:vw-tangles} The ($\calS$-graded) circuit algebra $\glos{\vT}$ of
v-tangles is the $\calS$-graded directed circuit algebra of v-tangle diagrams $\vD$, 
modulo the \Rs, R2 and R3 moves as depicted in
Figure~\ref{fig:VKnotRels}. These relations make sense as circuit
algebra relations between the two generators, and preserve skeleta. 
To obtain the circuit algebra $\wT$ of $w$-tangles we also mod out by the OC relation of Figure~\ref{fig:VKnotRels}
(note that each side in that relation involves only two generators,
with the apparent third ``virtual'' crossing being merely a circuit algebra artifact).
In fewer words, $\vT:=$\raisebox{-1.5mm}{\input{figs/vTDef.pstex_t}}, and
$\glos{\wT}:=$\raisebox{-1.8mm}{\input{figs/wTDef.pstex_t}}.
\end{definition}

\begin{figure}[t]
 \input figs/VKnotRels.pstex_t
 \caption{The relations (``Reidemeister moves'') \Rs, R2 and R3 define $v$-tangles, adding OC to these defines $w$-tangles. VR1, VR2, VR3 and M
 are not necessary as the circuit algebra presentation eliminates the need for ``virtual crossings'' as generators. R1 is not imposed for 
 framing reasons, and not imposing UC breaks the symmetry between over and under crossings in $\wT$.}
 \label{fig:VKnotRels}
 \end{figure}

\begin{remark}\label{rmk:VirtualXings}
One may also define v-tangles and w-tangles using the
language of planar algebras, except then another generator is required
(the ``virtual crossing'') and also a number of further relations shown in Figure \ref{fig:VKnotRels} (VR1--VR3,
M), and some of the operations (non-planar wirings) become less
elegant to define.
In our context ``virtual crossings''  are automatically
present (but unimportant) as part of the circuit algebra structure, and the ``virtual Reidemeister moves'' VR1--VR3 and M are 
also automatically true. In fact, the ``rerouting move'' known in the planar presentation, which says that a purely virtual strand
of a $v$-tangle diagram can be re-routed in any other purely virtual way, is precisely the statement that virtual crossings are
unimportant, and the language of circuit algebras makes this fact manifest.
\end{remark}

\begin{remark}\label{rmk:skeleta}
 For $S\in \calS$ a given skeleton, that is, an oriented 1-manifold with numbered ends, let us denote by 
 $\vT(S)$ and $\wT(S)$, respectively, the $v$- and $w$-tangles with skeleton $S$. 
 That is, $\vT(S)$ and $\wT(S)$ are the pre-images of $S$ under the skeleton map $\varsigma$. 
 Note that in our case the skeleton map
 is ``forgetting topology'', in other words, forgetting the under/over information of crossings, resulting in empty circuits. 
 With this notation, $\wT(\uparrow)$, the set of w-tangles whose skeleton is a single line, is exactly
 the set of (long) w-knots discussed in \cite[Section 3]{Bar-NatanDancso:WKO1}. Note also that $\wT(\uparrow_n)$, the set of w-tangles whose skeleton 
 is $n$ lines, includes w-braids with $n$ strands (\cite[Section 2]{Bar-NatanDancso:WKO1}) 
 but it is more general. Neither w-knots nor w-braids are circuit algebras.  
\end{remark}

 
\begin{remark}\label{rmk:Framing}
Since we do not mod out by the R1 relation, only by its weak (or
``spun'') version \Rs, it is more appropriate to call our class of
$v$/$w$-tangles {\em framed} $v$/$w$-tangles. (Recall that framed u-tangles
are characterized as the planar algebra generated by the positive and
negative crossings modulo the \Rs, R2 and R3 relations.) However, since
we are for the most part interested in studying the framed theories
(cf. Comment \ref{com:wTFFraming}), we will reserve the unqualified name
for the framed case, and will explicitly write ``unframed v/w-tangles''
if we wish to mod out by R1. For a more detailed explanation of framings
and R1 moves, see \cite[Remark~\ref{1-rem:Framing}]{Bar-NatanDancso:WKO1}.
\end{remark}

Our next task is to study the associated graded structures $\grad\vT$ and $\grad\wT$
of $\vT$ and $\wT$, with respect to the augmentation ideal as described in Section~\ref{subsec:Grad}. These are ``arrow diagram spaces on tangle skeletons'':
directed analogues of the chord diagram spaces of ordinary finite type invariant theory,
and even more similar to the arrow diagram spaces for braids and knots discussed
in \cite{Bar-NatanDancso:WKO1}. Our convention for figures will be to show skeletons 
as thick lines with thin arrows (directed chords).
Again, the language of circuit algebras makes defining these spaces
exceedingly simple.

\parpic[r]{\raisebox{-8mm}{$\input figs/arrows.pstex_t$}}
\begin{definition} The ($\calS$-graded) circuit algebra
$\glos{\calD^v}=\glos{\calD^w}$ of arrow diagrams is the graded and
$\calS$-graded directed circuit algebra generated by a single degree
1 generator $a$ in $C_{2,2}$ called ``the arrow'' as shown on the
right, with the obvious meaning for its skeleton. There are morphisms
$\pi\colon \calD^v\to\vT$ and $\pi\colon \calD^w\to\wT$ defined by
mapping the arrow to an overcrossing minus a no-crossing. (On the
right some virtual crossings were added to make the skeleta match). Let
$\glos{\calA^v}$ be $\calD^v/6T$, let
$\glos{\calA^w}:=\calA^v/TC=\calD^w/(\aft,TC)$, and let
$\glos{\calA^{sv}}:=\calA^v/RI$ and $\glos{\calA^{sw}}:=\calA^w/RI$, with RI, $6T$, $\aft$, and $TC$ being the relations shown
in Figures~\ref{fig:ADand6T} and~\ref{fig:TCand4T}. Note that the pair of relations $(\aft, TC)$ is
equivalent to the pair $(6T,TC)$,
as discussed in \cite[Section 2.3.1]{Bar-NatanDancso:WKO1}.
\end{definition}

\begin{figure}
 \input{figs/ADand6T.pstex_t}
 \caption{Relations for v-arrow diagrams on tangle skeletons.
Skeleta parts that are not connected can lie on separate
skeleton components; and the dotted arrow that remains in the same position means
``all other arrows remain the same throughout''. }
 \label{fig:ADand6T}
\end{figure}

\begin{figure}
 \input{figs/TCand4T.pstex_t}
 \caption{Relations for w-arrow diagrams on tangle skeletons.}
 \label{fig:TCand4T}
\end{figure}

\begin{proposition} The maps $\pi$ above induce surjections
$\pi\colon \calA^{sv}\to\grad\vT$ and $\pi\colon
\calA^{sw}\to\grad\wT$. Hence in the language of
Definition~\ref{def:CanGrad}, $\calA^{sv}$ and $\calA^{sw}$ are candidate
associated graded structures of $\vT$ and $\wT$.
\end{proposition}

\begin{proof} Proving that $\pi$ is well-defined amounts to checking
directly that the RI and 6T or RI, $\aft$ and TC relations are in the
kernel of $\pi$. (Just like in the finite type theory of virtual knots and
braids.) Thanks to the circuit algebra structure, it is enough to verify
the surjectivity of $\pi$ in degree 1. We leave this as an exercise for
the reader. \qed
\end{proof}

We do not know if $\calA^{sv}$ is indeed the associated graded of $\vT$ (also
see~\cite{Bar-NatanHalachevaLeungRoukema:v-Dims}). Yet in the w case, the
picture is simple:

\begin{theorem}\label{thm:ExpansionForTangles} The assignment $\overcrossing\mapsto e^a$ (with $e^a$
denoting the exponential of a single arrow from the over strand to the
under strand, interpreted via its power series) extends to a well defined $Z\colon \wT\to\calA^{sw}$. The
resulting map $Z$ is a homomorphic $\calA^{sw}$-expansion, and in particular,
$\calA^{sw}\cong\grad\wT$ and $Z$ is a homomorphic expansion.
\end{theorem}

\begin{proof} The proof is essentially the same as the proof of
\cite[Theorem~\ref{1-thm:RInvariance}]{Bar-NatanDancso:WKO1},
and follows \cite{BerceanuPapadima:BraidPermutation,
AlekseevTorossian:KashiwaraVergne}.
One needs to check that $Z$ satisfies the
Reidemeister moves and the OC relation. \Rs~follows easily from $RI$, R2 is obvious, TC implies OC. 
For R3, let $\calA^{sw}(\uparrow_n)$ denote the space of ``arrow diagrams on $n$ vertical strands''. 
We need to verify that $R:=e^a\in \calA^{sw}(\uparrow_2)$ satisfies the Yang-Baxter equation
$$R^{12}R^{13}R^{23}=R^{23}R^{13}R^{12}, \quad \text{ in } \calA^{sw}(\up_3),$$ 
where $R^{ij}=e^{a_{ij}}$ means ``place R on strands $i$ and $j$''.
By $4T$ and $TC$ relations, both sides of the equation can be reduced to $e^{a_{12}+a_{13}+a_{23}}$,
proving the Reidemeister invariance of $Z$.

$Z$ is by definition a circuit algebra homomorphism. 
Hence to show that $Z$ is an $\calA^{sw}$-expansion we only need to check the
universality property in degree one, where it is very easy.
The rest follows from Proposition~\ref{prop:CanGrad}. \qed
\end{proof}

\begin{remark}
 Note that the restriction of $Z$ to w-knots and w-braids (in the sense of Remark~\ref{rmk:skeleta})
 recovers the expansions constructed in \cite{Bar-NatanDancso:WKO1}. Note also that the 
 filtration and associated graded structure for w-braids fits into the general algebraic 
 framework of Section~\ref{sec:generalities} by applying the machinery to the skeleton-graded group
 of w-braids instead the circuit algebra of w-tangles. (The skeleton of a w-braid is the permutation it represents.) 
 However, as w-knots do not form a finitely
 presented algebraic structure in the sense of Section~\ref{sec:generalities}, the ``finite type'' filtration
 used in \cite{Bar-NatanDancso:WKO1} does not arise as powers of any augmentation ideal.
 This captures the reason why w-knots are ``the wrong objects to study'', as we have
 mentioned at the beginning of Section 3 of \cite{Bar-NatanDancso:WKO1}.
\end{remark}

In a similar spirit to
\cite[Definition~\ref{1-def:wJac}]{Bar-NatanDancso:WKO1}, one may define a
``w-Jacobi diagram'' on an arbitrary skeleton:

\begin{definition}\label{def:wJac} A ``w-Jacobi diagram on a tangle
skeleton''\footnote{We usually short this to
``w-Jacobi diagram'', or sometimes ``arrow diagram'' or just ``diagram''.}
is a graph made of the following ingredients:
\begin{itemize}
\item An oriented ``skeleton'' consisting of long lines and circles 
(i.e., an oriented one-manifold). In figures we draw the skeleton
  lines thicker.
\item Other directed edges, usually called ``arrows''.
\item Trivalent ``skeleton vertices'' in which an arrow starts or ends on
  the skeleton line.
\item Trivalent ``internal vertices'' in which two arrows end and one arrow
  begins. The internal vertices are cyclically oriented; in figures the assumed 
  orientation is always counterclockwise unless marked otherwise. Furthermore,
  all trivalent vertices must be connected to the skeleton via arrows (but not necessarily 
  following the direction of the arrows).
\end{itemize}
\end{definition}

Note that we allow multiple and loop arrow edges, as long as trivalence and the
two-in-one-out rule is respected.

Formal linear combinations of (w-Jacobi) arrow diagrams form a circuit algebra. 
We denote by $\calA^{wt}$ the quotient of the circuit algebra of arrow diagrams modulo the $\aSTU_1$, $\aSTU_2$ relations of 
Figure \ref{fig:STU}, and the TC relation. We denote $\calA^{wt}$ modulo the RI relation by $\calA^{swt}$. We then
have the following ``bracket-rise'' theorem:

\begin{figure}
 \input{figs/aSTU.pstex_t}
 \caption{The $\protect\aSTU$ relations for arrow diagrams, with their ``central edges'' marked $e$ for easier memorization.}
 \label{fig:STU}
\end{figure}

\begin{figure}
 \input{figs/aIHX.pstex_t}
 \caption{The $\protect\aAS$ and $\protect\aIHX$ relations.}
 \label{fig:aIHX}
\end{figure}

\begin{theorem}\label{thm:BracketRise} The obvious inclusion of arrow diagrams (with no internal vertices) into 
w-Jacobi diagrams descends to a map
$\bar{\iota}: \calA^w\to\calA^{wt}$, which is a circuit
algebra isomorphism. Furthermore, the $\aAS$
and $\aIHX$ relations of Figure~\ref{fig:aIHX} hold in $\calA^{wt}$.
Consequently, it is also true that $\calA^{sw}\cong\calA^{swt}$.
\end{theorem}

\begin{proof} In the proof of
\cite[Theorem~\ref{1-thm:BracketRise}]{Bar-NatanDancso:WKO1} we showed
this for long w-knots (i.e., tangles whose skeleton is a single long
line).  That proof applies here verbatim, noting that it does not make
use of the connectivity of the skeleton.

In short, to check that $\bar{\iota}$ is well-defined, we need to show that the $\aSTU$ relations
imply the $\aft$ relation. This is shown in Figure \ref{fig:STU4T}. To show that $\bar{\iota}$ is an isomorphism,
we construct an inverse $\calA^{wt} \to \calA^w$, which ``eliminates all internal vertices'' using
a sequence of $\aSTU$ relations. Checking that this is well-defined requires some case analysis;
the fact that it is an inverse to $\bar{\iota}$ is obvious. Verifying that the $\aAS$
and $\aIHX$ relations hold in $\calA^{wt}$ is an easy exercise. 
\begin{figure}
 \input{figs/STUto4T.pstex_t}
 \caption{Applying $\protect\aSTU_1$ and $\protect\aSTU_2$ to the diagram on the left, we get the two sides of $\protect\aft$.}
 \label{fig:STU4T}
\end{figure}
\qed
\end{proof}

Given the above theorem, we no longer keep the distinction between
$\calA^w$ and $\calA^{wt}$ and between $\calA^{sw}$ and $\calA^{swt}$.

We recall from \cite{Bar-NatanDancso:WKO1} that a ``$k$-wheel'', sometimes denoted $w_k$, is a
an arrow diagram consisting of an oriented cycle of arrows with $k$
incoming ``spokes'', the tails of which rest on the skeleton. An example
is shown in Figure \ref{fig:wheels}. In this language, the RI relation
can be rephrased using the $\aSTU$ relation to say that all one-wheels are 0,
or $w_1=0$.

\begin{figure}
 \input{figs/wkl.pstex_t}
 \caption{A $4$-wheel and the RI relation re-phrased.}
 \label{fig:wheels}
\end{figure}

\begin{remark} \label{rem:HeadInvariance}
Note that if $T$ is an arbitrary $w$ tangle, then the equality on the
left side of the figure below always holds, while the one on the right
generally doesn't:
\begin{equation} \label{eq:TangleLassoMove}
  \begin{array}{c}\input{figs/TangleLassoMove.pstex_t}\end{array}
\end{equation}
The
arrow diagram version of this statement is that if $D$ is an arbitrary
arrow diagram in $\calA^w$, then the left side equality in the
figure below always holds (we will sometimes refer to this as the
``head-invariance'' of arrow diagrams), while the right side equality
(``tail-invariance'') generally fails.
\begin{equation} \label{eq:HeadInvariance}
  \begin{array}{c}\input{figs/HeadInvariance.pstex_t}\end{array}
\end{equation}
We leave it to the reader to ascertain that
Equation~\eqref{eq:TangleLassoMove} implies
Equation~\eqref{eq:HeadInvariance}. There is also a direct
proof of Equation~\eqref{eq:HeadInvariance} which we also leave
to the reader, though see an analogous statement and proof in
\cite[Lemma~3.4]{Bar-Natan:NAT}. Finally note that a restricted version of
tail-invariance does hold --- see Section~\ref{subsec:sder}.
\end{remark}

\draftcut
\subsection{$\calA^w(\uparrow_n)$ and the Alekseev-Torossian Spaces}
\label{subsec:ATSpaces}

\begin{definition} Let $\glos{\calA^v(\uparrow_n)}$ be the part of $\calA^v$ in
which the skeleton is the disjoint union of $n$ directed lines,
with similar definitions for $\glos{\calA^w(\uparrow_n)}$,
$\glos{\calA^{sv}(\uparrow_n)}$, and $\glos{\calA^{sw}(\uparrow_n)}$.
\end{definition}

\begin{theorem}
\em{(Diagrammatic PBW Theorem.)} Let  $\glos{\calB^w_n}$ denote the space of uni-trivalent diagrams\footnote{
Oriented graphs with vertex degrees either 1 or 3, where trivalent vertices must have two edges incoming and
one edge outgoing and are cyclically oriented.}
with symmetrized ends coloured with colours in some $n$-element set
(say $\{x_1,\ldots,x_n\}$), modulo the $\aAS$ and $\aIHX$ relations of Figure \ref{fig:aIHX}.
Then there is an isomorphism $\calA^w(\uparrow_n)\cong \calB^w_n$.
\end{theorem}

{\it Proof sketch.} Readers familiar with the diagrammatic PBW theorem \cite[Theorem 8]{Bar-Natan:OnVassiliev}
will note that the proof carries through almost verbatim. There is a map $\chi: \calB^w_n \to \calA^w(\uparrow_n)$,
which sends each uni-trivalent diagram to the average of all ways of attaching their univalent ends
to the skeleton of $n$ lines, so that ends of colour $x_i$ are attached to the strand numbered $i$.
I.e., a diagram with $k_i$ uni-valent vertices of colour $x_i$ is sent to a sum of $\prod_i k_i!$ terms,
divided by $\prod_i k_i!$.

The goal is to show that $\chi$ is an isomorphism by constructing an inverse for it. The image of
$\chi$ are {\em symmetric} sums of diagrams, that is, sums of diagrams that are invariant under 
permuting arrow endings on the same skeleton component. 
One can show that in fact any arrow diagram $D$ in $\calA^w(\uparrow_n)$ is equivalent via $\aSTU$ and $TC$ relations 
to a symmetric sum. The obvious candidate is its ``symmetrization'' $Sym(D)$: the average of all ways of permuting the
arrow endings on each skeleton component of $D$. It is not true that each diagram is equivalent to its symmetrization
(hence, the ``simply delete the skeleton'' map is not an inverse for $\chi$),
but it is true that $D-Sym(D)$ has fewer skeleton vertices (lower degree) than $D$, hence we can
construct $\chi^{-1}$ inductively.
The fact that this inductive procedure is well-defined requires a proof; that proof is essentially the same as the proof 
of the corresponding fact in \cite[Theorem 8]{Bar-Natan:OnVassiliev}.
\qed

Both $\calA^w(\uparrow_n)$ and $\calB^w_n$ have a natural bi-algebra structure. In $\calA^w(\uparrow_n)$
multiplication is given by stacking. For a diagram $D\in \calA^w(\uparrow_n)$, the co-product $\glos{\Delta}(D)$
is given by the sum of all ways of dividing $D$ between a ``left co-factor'' and a ``right cofactor'' so
that the connected components of $D-S$ are kept intact, where $S$ is the skeleton of $D$. In $\calB^w_n$
multiplication is given by disjoint union, and $\Delta$ is the sum of all ways of dividing the connected 
components of a diagram between two co-factors (here there is no skeleton).
Note that the isomorphism $\chi$ above is a co-algebra isomorphism, but not an algebra homomorphism.

The primitives $\glos{\calP^w_n}$ of $\calB^w_n$ are the
connected diagrams (and hence the primitives of $\calA^w(\uparrow_n)$
are the diagrams that remain connected even when the skeleton is
removed). Given the ``two in one out'' rule for internal vertices,
the diagrams in $\calP^w_n$ can only be trees (diagrams with no cycles) or wheels (a
single oriented cycle with a number of ``spokes'', or leaves, attached to it). ``Wheels of
trees'' can be reduced to simple wheels by repeatedly using $\aIHX$,
as in Figure~\ref{fig:WheelOfTreesAndPrince}.

\begin{figure}
\input{figs/WheelOfTrees.pstex_t}
\caption{A wheel of trees can be reduced to a combination of wheels, and a wheel of trees with 
a Little Prince.}\label{fig:WheelOfTreesAndPrince}
\end{figure}

Thus as a vector space $\calP^w_n$ is easy to identify. It is a direct sum
$\calP^w_n=\langle\text{trees}\rangle\oplus\langle\text{wheels}\rangle$.
The wheels part is simply the graded vector space generated by
all cyclic words in the letters $x_1,\ldots,x_n$.  Alekseev and
Torossian~\cite{AlekseevTorossian:KashiwaraVergne} denote the
space of cyclic words by $\glos{\attr_n}$, and so shall we. The trees in
$\calP^w_n$ have leafs coloured $x_1,\ldots,x_n$. Modulo $\aAS$ and
$\aIHX$, they correspond to elements of the free Lie algebra $\glos{\lie_n}$
on the generators $x_1,\ldots,x_n$. But the root of each such tree
also carries a label in $\{x_1,\ldots,x_n\}$, hence there are $n$
types of such trees as separated by their roots, and so $\calP^w_n$
is linearly isomorphic to the direct sum $\attr_n\oplus\bigoplus_{i=1}^n\lie_n$.

Note that with $\calB_n^{sw}$ and $\calP_n^{sw}$ defined in the analogous manner
(i.e., factoring out by one-wheels, as in the RI relation),
we can also conclude that there is a linear isomorphism
$\calP^{sw}_n\cong\attr_n/(\text{deg }1)\oplus\bigoplus_{i=1}^n\lie_n$.

By the Milnor-Moore theorem~\cite{MilnorMoore:Hopf}, $\calA^w(\uparrow_n)$
is isomorphic to the universal enveloping algebra $\calU(\calP^w_n)$,
with $\calP^w_n$ identified as the subspace $\glos{\calP^w(\uparrow_n)}$
of primitives of $\calA^w(\uparrow_n)$ using the PBW symmetrization
map $\chi\colon \calB^w_n\to\calA^w(\uparrow_n)$. Thus in order to
understand $\calA^w(\uparrow_n)$ as an associative algebra, it is enough
to understand the Lie algebra structure induced on $\calP^w_n$ via the
commutator bracket of $\calA^w(\uparrow_n)$.

Our goal is to identify $\calP^w(\uparrow_n)$ as the Lie algebra
$\attr_n\rtimes(\fraka_n\oplus\tder_n)$,
which in itself is a combination of the Lie algebras
$\fraka_n$, $\tder_n$ and $\attr_n$ studied by Alekseev and
Torossian~\cite{AlekseevTorossian:KashiwaraVergne}. Here are the relevant
definitions:

\begin{definition} Let $\glos{\fraka_n}$ denote the vector space with basis
$x_1,\ldots,x_n$, also regarded as an Abelian Lie algebra of dimension $n$.
As before, let $\lie_n=\lie(\fraka_n)$ denote the free Lie algebra on $n$
generators, now identified as the basis elements of $\fraka_n$. Let
$\glos{\der_n}=\der(\lie_n)$ be the (graded) Lie algebra of derivations
acting on $\lie_n$, and let
\[ \glos{\tder_n}=\left\{D\in\der_n\colon \forall i\ \exists a_i\text{ s.t.{}
  }D(x_i)=[x_i,a_i]\right\}
\]
denote the subalgebra of ``tangential derivations''. A tangential
derivation $D$ is determined by the $a_i$'s for which $D(x_i)=[x_i,a_i]$,
and determines them up to the ambiguity $a_i\mapsto a_i+\alpha_ix_i$, where
the $\alpha_i$'s are scalars. Thus as vector spaces,
$\fraka_n\oplus\tder_n\cong\bigoplus_{i=1}^n\lie_n$.
\end{definition}

\begin{definition} Let $\glos{\Ass_n}=\calU(\lie_n)$ be the free associative
algebra ``of words'', and let $\glos{\Ass_n^+}$ be the degree $>0$ part of
$\Ass_n$. As before, we let $\attr_n=\Ass^+_n/(x_{i_1}x_{i_2}\cdots
x_{i_m}=x_{i_2}\cdots x_{i_m}x_{i_1})$ denote ``cyclic words'' or
``(coloured) wheels''.  $\Ass_n$, $\Ass_n^+$, and $\attr_n$ are
$\tder_n$-modules and there is an obvious equivariant ``trace''
$\tr\colon \Ass^+_n\to\attr_n$.
\end{definition}

\begin{proposition}\label{prop:Pnses}
There
is a split short exact sequence of Lie algebras 
\[ 0 \longrightarrow \attr_n
  \stackrel{\glos{\iota}}{\longrightarrow} \calP^w(\uparrow_n) 
  \stackrel{\glos{\pi}}{\longrightarrow} \fraka_n \oplus \tder_n 
  \longrightarrow 0.
\]
\end{proposition}

\begin{proof}
The inclusion $\iota$ is defined the natural way: $\attr_n$ is
spanned by coloured ``floating'' wheels, and such a wheel is mapped
into $\calP^w(\uparrow_n)$ by attaching its ends to their assigned strands in
arbitrary order. Note that this is well-defined: wheels have only tails,
and tails commute.

 As vector spaces, the statement is already proven: $\calP^w(\uparrow_n)$ 
is generated by trees 
and wheels (with the all arrow endings fixed on $n$ strands). When factoring out by the wheels,
only trees remain. Trees have one head and many tails. All the tails commute with 
each other, and commuting a tail with a head on a strand costs a wheel (by $\aSTU$), 
thus in the quotient the head also commutes with the tails. Therefore, the quotient
is the space of coloured ``floating'' trees, which we have previously identified with
$\bigoplus_{i=1}^{n} \lie_n \cong \fraka_n\oplus\tder_n$.

It remains to show that the maps $\iota$ and $\pi$ are Lie algebra maps as well. For $\iota$ this
is easy: the Lie algebra $\attr_n$ is commutative, and is mapped to the commutative
(due to $TC$)
subalgebra of $\calP^w(\uparrow_n)$ generated by wheels. Next, we show that $\pi$ is a homomorphism.  The map $\pi$ quotients out by wheels and reads trees as lie words in $\lie_n^{\oplus n}$, as show in Figure~\ref{fig:MapPi} and below: 
\begin{align*}
           \pi:                                &\, \calP^w(\uparrow_n)  \to                        \quad \fraka_n \oplus \tder_n \quad \stackrel{\text{v.s.}}{\cong} \quad \quad \quad \lie_n^{\oplus n} \\
  \langle \text{wheels} \rangle  &    \ni w                     \quad \,  \, \mapsto      \quad \quad  \,\,  0 \\
   \langle \text{trees} \rangle    & \ni T                         \quad \, \,  \mapsto     \quad \quad  \quad \quad \quad \quad  \quad \quad \,  (0,\ldots, a_i,\ldots, 0)
\end{align*}
Here $a_i$ is the lie word corresponding to $T$ in the $i$-th component of $\lie_n^{\oplus n}$, where the head of $T$ is on strand $i$. Namely, $a_i$ is a lie word on the generators corresponding to the strand numbers to which the tails of $T$ are attached, and commutators corresponding to each of the arrow vertex of $T$, read left to right when looking at $T$ from its head.

\begin{figure}
\[
  \def\goesto{{$\overset{\pi}{\longmapsto} (0,[[x_1,x_2],[x_2,x_3]],0)$}}
  \input{figs/ArrowsToLie.pstex_t}
\]
\caption{The map $\pi$.}\label{fig:MapPi}
\end{figure}

To show that $\pi$ is a map of Lie algebras we give two proofs,
first a ``hands-on'' one, then a ``conceptual'' one.

{\bf Hands-on argument.} $\fraka_n$ is the image of single arrows on one strand.
These commute with everything in $\calP^w(\uparrow_n)$, and so does $\fraka_n$
in the direct sum $\fraka_n \oplus \tder_n$. Thus, $\pi$ respects commutators involving these local arrows. 

It remains to show that commuting trees in $\calP^w(\uparrow_n)$ maps to the bracket of $\tder_n$; or 
more accurately, the bracket of $\fraka_n \oplus \tder_n$ transferred to $\lie_n$. Let $D$ and $D'$ be elements of
$\tder_n$ represented by $(a_1,\ldots ,a_n)$ and $(a_1',\ldots ,a_n')$, meaning
that $D(x_i)=[x_i,a_i]$ and $D'(x_i)=[x_i,a_i']$ for $i=1,\ldots ,n$. We compute the commutator of these elements:
\begin{multline*}
  [D,D'](x_i)=(DD'-D'D)(x_i)=D[x_i,a_i']-D'[x_i,a_i]= \\
  =[[x_i,a_i],a_i']+[x_i,Da_i']-[[x_i,a_i'],a_i]-[x_i,D'a_i]
    = [x_i,Da_i'-D'a_i+[a_i,a_i']].
\end{multline*}

Now let $T$ and $T'$ be two trees in $\calP^w(\uparrow_n)/\attr_n$,
with heads on strands $i$ and $j$, respectively ($i$ may or may not
equal $j$).  Let us denote by $a_i$ (resp. $a_j'$) the element in $\lie_n$ corresponding to $T$ (resp. $T'$), as above. 
In $\tder_n$, let $D=\pi(T)=(0,\ldots,-a_i,\ldots,0)$ and
$D'=\pi(T')=(0,\ldots,-a_j,\ldots,0)$.  (In each case, the $i$-th, respectively
$j$-th, is the only non-zero component.) The commutator of these
elements is given by $[D,D'](x_i)=[Da_i'-D'a_i+[a_i,a_i'],x_i]$, and
$[D,D'](x_j)=[Da_j'-D'a_j+[a_j,a_j'],x_j].$ Note that unless $i=j$,
$a_j=a_i'=0$.

In $\calP^w(\uparrow_n)/\attr_n$ tails commute, as well as the head of a tree with its
own tails. Therefore, commuting two trees only incurs a cost when commuting a head of
one tree over the tails of the other on the same strand, and the two heads over each other,
if $i=j$.

If $i \neq j$, then commuting the head of $T$ over the tails of $T'$ by $\aSTU$ 
gives a sum of trees given by $-Da_j'$, with heads on strand $j$, while moving
the head of $T'$ over the tails of $T$ costs exactly $D'a_i$, with heads on strand $i$,
as needed.

If $i=j$, then everything happens on strand $i$, and the cost is 
$(-Da_i'+D'a_i-[a_i,a_i'])$, where the last term arises from commuting the two heads.

{\bf Conceptual argument.}
There is an action of $\calP^w(\uparrow_n)$ on $\lie_n$, as follows: introduce
and extra strand on the right. An element $L$ of $\lie_n$ corresponds to a tree with 
its head on the extra strand. Its commutator with an element of $\calP^w(\uparrow_n)$ 
(considered as an element of $\calP^w(\uparrow_{n+1})$ by the obvious inclusion)
is again a tree with head on strand $(n+1)$, defined to be the result of the action.

Since $L$ has only tails on the first $n$ strands, 
elements of $\attr_n$, which
also only have tails, act trivially. So do single (local) arrows on one strand
($\fraka_n$). It remains to show that trees act as $\tder_n$, and it is enough
to check this on the generators of $\lie_n$ (as the Leibniz rule is obviously
satisfied). The generators of $\lie_n$ are arrows pointing from one of the first 
$n$ strands, say strand $i$, to strand $(n+1)$. A tree $T$ with head on strand $i$
acts on this element, according $\aSTU$, by forming the commutator $[x_i,T]$, which
is exactly the action of $\tder_n$.
\end{proof}

To identify $\calP^w(\uparrow_n)$ as the semidirect product
$\attr_n\rtimes(\fraka_n\oplus\tder_n)$, it remains to show that
the short exact sequence of the Proposition splits. This is indeed the case,
although not canonically.  Two ---of the many--- splitting maps
$\glos{u},\glos{l}\colon \tder_n\oplus\fraka_n \to \calP^w(\uparrow_n)$
are described as follows: $\tder_n\oplus\fraka_n$ is identified with
$\bigoplus_{i=1}^n\lie_n$, which in turn is identified with ``floating''
coloured trees. A map to $\calP^w(\uparrow_n)$ can
be given by specifying how to place the legs on their specified strands.
A tree may have many tails but has only one head, and due to $TC$, only
the positioning of the head matters. Let $u$ (for {\it upper}) be the map
placing the head of each tree above all its tails on the same strand,
while $l$ (for {\it lower}) places the head below all the tails. It is
clear that these are both Lie algebra maps and that $\pi \circ u$ and
$\pi \circ l$ are both the identity of $\tder_n \oplus \fraka_n$. This
makes $\calP^w(\uparrow_n)$ a semidirect product. \qed

\begin{remark} Let $\glos{\attr_n^s}$ denote $\attr_n$ mod out by its
degree one part (one-wheels). Since the RI relation is in the kernel of
$\pi$, there is a similar split exact sequence
\[ 0\to \attr_n^s \stackrel{\overline{\iota}}{\rightarrow} \calP^{sw}
  \stackrel{\overline{\pi}}{\rightarrow} \fraka_n \oplus \tder_n.
\]
\end{remark}

\begin{definition}\label{div} 
For any $D \in \tder_n$, $(l-u)D$ is in the kernel of $\pi$, therefore
is in the image of $\iota$, so $\iota^{-1}(l-u)D$ makes sense. We call
this element $\glos{\divop}D$.
\end{definition}

\begin{definition}
In \cite{AlekseevTorossian:KashiwaraVergne} 
div is defined as follows: div$(a_1,\ldots,a_n):=\sum_{k=1}^n \tr((\partial_k a_k)x_k)$,
where $\partial_k$ picks out the words of a sum which end in $x_k$ and deletes their last letter
$x_k$, and deletes all other words (the ones which do not end in $x_k$).
\end{definition}

\begin{proposition}
The div of Definition \ref{div} and the div of \cite{AlekseevTorossian:KashiwaraVergne} are 
the same.  
\end{proposition}

\parpic[r]{\input{figs/combtree.pstex_t}}
{\it Proof.}
It is enough to verify the claim for the linear generators of $\tder_n$, namely, elements
of the form $(0,\ldots,a_j,\ldots,0)$, where $a_j \in \lie_n$ or equivalently, single (floating, 
coloured) trees, where the colour of
the head is $j$. By the Jacobi identity, each $a_j$ can be written 
in a form $a_j=[x_{i_1},[x_{i_2},[\ldots,x_{i_k}]\ldots]$. 
Equivalently, by $\aIHX$, each tree has a 
standard ``comb'' form, as shown on the picture on the right.

For an associative word $Y=y_1y_2\ldots y_l \in \Ass_n^+$, 
we introduce the notation $[Y]:=[y_1,[y_2,[\ldots,y_l]\ldots]$.
The div of \cite{AlekseevTorossian:KashiwaraVergne} picks out the
words that end in $x_j$, forgets the rest, and considers these as
cyclic words. Therefore, by interpreting the Lie brackets as commutators,
one can easily check that for $a_j$ written as above,
\begin{equation}\label{divformula}
{\rm div}((0,\ldots,a_j,\ldots,0))=\sum_{\alpha\colon  i_{\alpha}=x_j} 
-x_{i_1}\ldots x_{i_{\alpha-1}}[x_{i_{\alpha+1}}\ldots x_{i_k}]x_j.
\end{equation}

\parpic[r]{\input{figs/divproof.pstex_t}}
In Definition \ref{div}, div of a tree is the difference between attaching its
head on the appropriate strand (here, strand $j$) below all of its tails and above.
As shown in the figure on the right, moving the head across each of the tails on 
strand $j$ requires an $\aSTU$ relation,
which ``costs'' a wheel (of trees, which is equivalent to a sum of honest wheels). 
Namely, the head gets connected to the tail in question.
So div of the tree represented by $a_j$ is given by
\begin{center}
$\sum_{\alpha\colon  x_{i_{\alpha}}=j}$``\rm connect the head to the $\alpha$ leaf''.
\end{center}

\noindent
This in turn gets mapped to the formula above via the correspondence between 
wheels and cyclic words. \qed

\parpic[r]{\input{figs/treeactonwheel.pstex_t}}
\begin{remark}\label{rem:tderontr}
There is an action of $\tder_n$ on $\attr_n$ as
follows. Represent a cyclic word $w \in \attr_n$ as a
wheel in $\calP^w(\uparrow_n)$ via the map $\iota$. Given
an element $D \in \tder_n$, $u(D)$, as defined above, is a tree
in $\calP^w(\uparrow_n)$ whose head is above all of its tails. We
define $D \cdot w:=\iota^{-1}(u(D)\iota(w)-\iota(w)u(D))$. Note that
$u(D)\iota(w)-\iota(w)u(D)$ is in the image of $\iota$, i.e., a linear
combination of wheels, for the following reason. The wheel $\iota(w)$ has only tails. As we commute
the tree $u(D)$ across the wheel, the head of the tree is commuted
across tails of the wheel on the same strand. Each time this happens
the cost, by the $\aSTU$ relation, is a wheel with the tree attached
to it, as shown on the right, which in turn (by $\aIHX$ relations,
as Figure~\ref{fig:WheelOfTreesAndPrince} shows) is a sum of wheels.
Once the head of the tree has been moved to the top, the tails of the
tree commute up for free by $TC$. Note that the alternative definition,
$D \cdot w:=\iota^{-1}(l(D)\iota(w)-\iota(w)l(D))$ is in fact equal to
the definition above.
\end{remark}

\begin{definition}
In \cite{AlekseevTorossian:KashiwaraVergne}, the group $\glos{\TAut_n}$
is defined as $\exp(\tder_n)$. Note that $\tder_n$ is positively
graded, hence it integrates to a group. Note also that $\TAut_n$ is
the group of ``basis-conjugating'' automorphisms of $\lie_n$, i.e.,
for $g \in \TAut_n$, and any $x_i$, $i=1,\ldots ,n$ generator of
$\lie_n$, there exists an element $g_i \in \exp(\lie_n)$ such that
$g(x_i)=g_i^{-1}x_ig_i$.
\end{definition}

Note that the group multiplication in $\TAut_n$ is the one exponentiated from $\tder_n$, which is read left to right (as right actions) rather than right to left (as function composition). For example, for $f,g \in \TAut_n$, $(fg)(x_i)=g(f(x_i))=g(f_i^{-1}x_if_i)= g(f_i)^{-1}g_i^{-1}x_ig_ig(f_i)$. 

The action of $\tder_n$ on $\attr_n$ lifts to an action of $\TAut_n$ on $\attr_n$,
by interpreting exponentials formally, in other words $e^D$ acts as 
$\sum_{n=0}^\infty\frac{D^n}{n!}$. The lifted action is by conjugation:
for $w \in \attr_n$ and $e^D \in \TAut_n$, 
$e^D \cdot w=\iota^{-1}(e^{uD} \iota(w) e^{-uD})$.

Recall that in Section 5.1 of \cite{AlekseevTorossian:KashiwaraVergne}
Alekseev and Torossian construct a map $\glos{j}\colon \TAut_n \to
\attr_n$ which is characterized by two properties: the cocycle property
\begin{equation}\label{eq:jcocycle}
 j(gh)=j(g)+g\cdot j(h),
\end{equation}
where in the second term multiplication by $g$ denotes the action described above;
and the condition
\begin{equation}\label{eq:jderiv}
\frac{d}{ds}j(\exp(sD))|_{s=0}=\divop(D). 
\end{equation}

Now let us interpret $j$ in our context.
\begin{definition}\label{def:Adjoint}
The adjoint map $\glos{*}\colon \calA^w(\uparrow_n) \to
\calA^w(\uparrow_n)$ acts by ``flipping over diagrams and negating arrow
heads on the skeleton''. In other words, for an arrow diagram $D$,
\[ D^*:=(-1)^{\#\{\text{tails on skeleton}\}}S(D), \]
where $S$ denotes the map which switches the orientation of the skeleton
strands (i.e. flips the diagram over), and multiplies by $(-1)^{\#
\text{skeleton vertices}}$.
\end{definition}

Note that the number of tails on the skeleton is the same as the degree of the arrow diagram, hence $D^*=(-1)^{\deg D}S(D)$.

\begin{proposition}\label{prop:Jandj}For $D \in \tder_n$,
define a map $\glos{J}\colon \TAut_n \to \exp(\attr_n)$ by
$J(e^D):=e^{uD}(e^{uD})^*$. Then
$$\exp(j(e^D))=J(e^D).$$
\end{proposition}

\begin{proof}
Note that $(e^{uD})^*=e^{-lD}$, due to ``Tails Commute'' and the fact that a 
tree has only one head.

Let us check that $\log J$ satisfies properties \eqref{eq:jcocycle} and 
\eqref{eq:jderiv}. Namely, with $g=e^{D_1}$ and $h=e^{D_2}$, and 
using that $\attr_n$ is commutative, we need to show that
\begin{equation}
 J(e^{D_1}e^{D_2})=J(e^{D_1})\big(e^{uD_1}\cdot J(e^{D_2})\big),
\end{equation}
where $\cdot$ denotes the action of $\tder_n$ on $\attr_n$; and that
\begin{equation}
 \frac{d}{ds}J(e^{sD})|_{s=0}=\divop D.
\end{equation}

Indeed, with $\operatorname{BCH}(D_1,D_2)=\log e^{D_1}e^{D_2}$ being the 
standard Baker--Campbell--Hausdorff formula,
\begin{multline*}
  J(e^{D_1}e^{D_2})=J(e^{\operatorname{BCH}(D_1,D_2)})
  =e^{u(\operatorname{BCH}(D_1,D_2)}
  e^{-l(\operatorname{BCH}(D_1,D_2)}=
  e^{\operatorname{BCH}(uD_1,uD_2)}
  e^{-\operatorname{BCH}(lD_1,lD_2)} \\
  =e^{uD_1}e^{uD_2}e^{-lD_2}e^{-lD_1}=
  e^{uD_1}(e^{uD_2}e^{-lD_2})e^{-uD_1}e^{uD_1}e^{lD_1}
  =(e^{uD_1}\cdot J(D_2))J(D_1),
\end{multline*}
as needed.

As for condition~\eqref{eq:jderiv}, a direct computation of the derivative
yields
$$\frac{d}{ds}J(e^{sD})|_{s=0}=uD-lD=\divop D,$$
as desired. \qed
\end{proof}

\draftcut
\subsection{The Relationship with u-Tangles} \label{subsec:sder} Let
$\glos{\uT}$ be the planar algebra of classical, or ``{\it u}sual''
tangles.  There is a map $a\colon \uT \to \wT$ of $u$-tangles into
$w$-tangles: algebraically, it is defined in the obvious way on the planar
algebra generators of $\uT$. (It can also be interpreted topologically
as Satoh's tubing map, see
\cite[Section~\ref{1-subsubsec:TopTube}]{Bar-NatanDancso:WKO1},
where a u-tangle is a tangle drawn on a sphere. However, it is only
conjectured that the circuit algebra presented here is a Reidemeister
theory for ``tangled ribbon tubes in $\bbR^4$''.)  The map $a$ induces a
corresponding map $\alpha\colon  \calA^u \to \calA^{sw}$, which maps an
ordinary Jacobi diagram (i.e., unoriented chords with internal trivalent
vertices modulo the usual $AS$, $IHX$ and $STU$ relations) to the sum
of all possible orientations of its chords (many of which are zero in
$\calA^{sw}$ due to the ``two in one out'' rule).

\parpic[l]{$\xymatrix{
  \uT \ar@{.>}[r]^{Z^u} \ar[d]^a & \calA^u \ar[d]^\alpha \\
  \wT \ar[r]^{Z^w} & \calA^{sw}
}$}
It is tempting to ask whether the square on the left
commutes. Unfortunately, this question hardly makes sense, as there
is no canonical choice for the dotted line in it. Similarly to the
braid case of \cite[Section~\ref{1-subsubsec:RelWithu}]{Bar-NatanDancso:WKO1}, the definition of the
homomorphic expansion (Kontsevich integral) for $u$-tangles typically depends on various choices
of ``parenthesizations''. Choosing parenthesizations, this square becomes
commutative up to some fixed corrections. The details are in
Proposition~\ref{prop:uwBT}.

Yet already at this point we can recover something from the existence of
the map $a\colon\uT\to\wT$, namely an interpretation of the
Alekseev-Torossian~\cite{AlekseevTorossian:KashiwaraVergne} space of
special derivations, $$\glos{\sder_n}:=\{ D\in\tder_n\colon D(\sum_{i=1}^n
x_i)=0\}.$$ Recall from Remark \ref{rem:HeadInvariance} that
in general it is not possible to slide a strand under an arbitrary $w$-tangle.
However, it is possible to slide strands freely under
tangles {\em in the image of $a$}, and thus by reasoning similar to the
reasoning in Remark~\ref{rem:HeadInvariance}, diagrams $D$ in the image
of $\alpha$ respect ``tail-invariance'':
\begin{equation} \label{eq:TailInvariance}
  \begin{array}{c}\input{figs/TailInvariance.pstex_t}\end{array}
\end{equation}

Let $\calP^u(\uparrow_n)$ denote the primitives of $\calA^u(\uparrow_n)$,
that is, Jacobi diagrams that remain connected when the skeleton is
removed. Remember that $\calP^{w}(\uparrow_n)$ stands for the primitives
of $\calA^{w}(\uparrow_n)$. Equation~\eqref{eq:TailInvariance} readily
implies that the image of the composition
\[ \xymatrix{
  \calP^u(\uparrow_n) \ar[r]^(0.48){\alpha}
  & \calP^w(\uparrow_n) \ar[r]^(0.45)\pi
  & \fraka_n \oplus \tder_n
} \]
is contained in $\fraka_n \oplus \sder_n$. Even better is true.

\begin{theorem}\label{thm:sder}
The image of $\pi\alpha$ is precisely $\fraka_n \oplus \sder_n$. 
\end{theorem}

This theorem was first proven by Drinfel'd (Lemma after Proposition 6.1
in \cite{Drinfeld:GalQQ}), but the proof we give here is due to Levine
\cite{Levine:Addendum}.

\begin{proof}
Let $\lie_n^d$ denote the degree $d$ piece of $\lie_n$. Let $V_n$ be
the vector space with basis $x_1, x_2, \ldots , x_n$.  Note that
$$V_n \otimes \lie_n^d \cong \bigoplus_{i=1}^n \lie_n^d \cong
(\tder_n \oplus \fraka_n)^d,$$
where $\tder_n$ is graded by the number of tails of a tree, and $\fraka_n$ 
is contained in degree 1.  

The bracket defines a map $\beta\colon  V_n \otimes \lie_n^d \to \lie_n^{d+1}$:
for $a_i \in \lie_n^d$ where $i=1,\ldots ,n$, the ``tree'' 
$D=(a_1,a_2,\ldots ,a_n) \in (\tder_n \oplus \fraka_n)^d$ is mapped to 
$$\beta(D)=\sum_{i=1}^n[x_i,a_i]=D\left(\sum_{i=1}^n x_i\right),$$
where the first equality is by the definition of tensor product and the bracket,
and the second is by the definition of the action of $\tder_n$ on $\lie_n$.

Since $\fraka_n$ is contained in degree 1, by definition 
$\sder_n^d=(\operatorname{ker}\beta)^d$ for $d\geq2$. In degree 
1, $\fraka_n$ is obviously in the kernel, hence 
$(\operatorname{ker}\beta)^1= \fraka_n \oplus \sder_n^1$. So overall,
$\operatorname{ker}\beta=\fraka_n\oplus\sder_n$.

We want to study the image of the map $\calP^u(\uparrow^n)
\stackrel{\pi\alpha}{\longrightarrow} \fraka_n \oplus \tder_n$.
Under $\alpha$, all connected Jacobi diagrams that are not trees or
wheels go to zero, and under $\pi$ so do all wheels. Furthermore, $\pi$
maps trees that live on $n$ strands to ``floating'' trees with univalent
vertices coloured by the strand they used to end on. So for determining
the image, we may replace $\calP^u(\uparrow^n)$ by the space $\calT_n$
of connected {\em un}oriented ``floating trees'' (uni-trivalent graphs), the ends (univalent vertices)
of which are coloured by the $\{x_i\}_{i=1,..,n}$. We denote the degree
$d$ piece of $\calT_n$, i.e., the space of trees with $d+1$ ends,
by $\calT_n^{d}$. Abusing notation, we shall denote the map induced by
$\pi\alpha$ on $\calT_n$ by $\alpha\colon  \calT_n \to \fraka_n \oplus
\tder_n$. Since choosing a ``head'' determines the entire orientation of
a tree by the two-in-one-out rule, $\alpha$ maps a tree in $\calT_n^d$
to the sum of $d+1$ ways of choosing one of the ends to be the ``head''.

We want to show that $\operatorname{ker}\beta=\operatorname{im}\alpha$.
This is equivalent to saying that $\bar{\beta}$ is injective, where
$\bar{\beta}\colon V_n\otimes\lie_n/\operatorname{im}\alpha
\to \lie_n$ is map induced by $\beta$ on the quotient by
$\operatorname{im}\alpha$.

\parpic[r]{\input{figs/beta.pstex_t}}
The degree $d$ piece of $V_n \otimes \lie_n$, in the pictorial
description, is generated by floating trees with $d$ tails and one head,
all coloured by $x_i$, $i=1,\ldots ,n$. This is mapped to $\lie_n^{d+1}$,
which is isomorphic to the space of floating trees with $d+1$ tails and
one head, where only the tails are coloured by the $x_i$. The map $\beta$
acts as shown on the picture on the right.

\parpic[r]{\input{figs/taudef.pstex_t}}
We show that $\bar{\beta}$ is injective by exhibiting a map $\tau\colon
\lie_n^{d+1} \to V_n\otimes\lie_n^d/\operatorname{im}\alpha$ so that
$\tau\bar{\beta}=I$. The map $\tau$ is defined as follows: given a tree with
one head and $d+1$ tails $\tau$ acts by deleting the head and the
arc connecting it to the rest of the tree and summing over all ways of
choosing a new head from one of the tails on the left half of the tree relative to the
original placement of the head (see the
picture on the right). As long as we show that $\tau$ is well-defined,
it follows from the definition and the pictorial description of $\beta$
that $\tau\bar{\beta}=I$.

For well-definedness we need to check that the images of $\aAS$ and
$\aIHX$ relations under $\tau$ are in the image of $\alpha$. This we do 
in the picture below. In both cases it is enough to check the
case when the ``head'' of the relation is the head of the tree
itself, as otherwise an $\aAS$ or $\aIHX$ relation in the domain is mapped
to an $\aAS$ or $\aIHX$ relation, thus zero, in the image.
\[ \input figs/tauproof.pstex_t \]
\[ \input figs/tauproof2.pstex_t \]
In the $\aIHX$ picture, in higher degrees $A$, $B$ and $C$ may denote
an entire tree. In this case, the arrow at $A$ (for example) means the
sum of all head choices from the tree $A$.
\qed
\end{proof}

\begin{comment} In view of the relation between the right half of
Equation~\eqref{eq:TailInvariance} and the special derivations $\sder$,
it makes sense to call w-tangles that satisfy the condition in the left
half of Equation~\eqref{eq:TailInvariance} ``special''. The $a$ images
of u-tangles are thus special. We do not know if the global version of
Theorem~\ref{thm:sder} holds true. Namely, we do not know whether every
special w-tangle is the $a$-image of a u-tangle.
\end{comment}

\draftcut
\subsection{The local topology of w-tangles}\label{subsec:TangleTopology}
So far throughout this section we have presented $w$-tangles as
a Reidemeister theory: a circuit algebra given by generators and
relations. There is a topological intuition behind this definition:
we can interpret the strands of a w-tangle diagram as embedded
tubes in $\bbR^4$ (with oriented ``cores'': 1D curves that
run along them), as shown in Figure \ref{fig:CrossingTubes}.
For each tube there exists a 3-dimensional ``filling'', and each crossing
represents a ribbon intersection between the filled tubes where the
one corresponding to the under-strand intersects the filling of
the over-strand. (For an explanation of ribbon singularities see
\cite[Section~\ref{1-subsubsec:ribbon}]{Bar-NatanDancso:WKO1}.)
In Figure~\ref{fig:CrossingTubes} we use the drawing conventions of
\cite{CarterSaito:KnottedSurfaces}: we draw surfaces as if projected
from $\bbR^4$ to $\bbR^3$, and cut them open when they are ``hidden''
by something with a higher 4-th coordinate.

\begin{figure}
 \input{figs/RibbonTubes.pstex_t}
 \caption{
  Strands correspond to tubes with cores, a virtual crossing
  corresponds to non-interacting tubes, while a crossing means that
  the tube corresponding to the under strand ``goes through'' the tube
  corresponding to the over strand.
} \label{fig:CrossingTubes}
\end{figure}

Note that w-braids can also be thought of in terms of flying circles,
with ``time'' being the fourth dimension; this is equivalent to
the tube interpretation in the obvious way. In this language a
crossing represents a circle (the under strand), flying through
another (the over strand). This is described in detail in
\cite[Section~\ref{1-subsubsec:FlyingRings}]{Bar-NatanDancso:WKO1}.

The assignment of tangled ribbon tubes in $\bbR^4$ to w-tangles is
well-defined (the Reidemeister and OC relations are satisfied), and after
Satoh \cite{Satoh:RibbonTorusKnots} we call it the tubing map and denote
it by $\glos{\delta}\colon\{\text{w-tangles}\} \to \{\text{Ribbon tubes in
} \bbR^4\}$. It is natural to expect that $\delta$ is an isomorphism, and indeed it is a 
surjection. However,
the injectivity of $\delta$ remains unproven even for long w-knots.  Nonetheless, ribbon
tubes in $\bbR^4$ will serve as the topological motivation and local
topological interpretation behind the circuit algebras presented in
this paper. 

\parpic[r]{\input{figs/TubeOrientation_2.pstex_t}}
Ribbon tubes in the image of $\delta$ are directed: 
where the direction comes from the direction of each tube as a strand of the tangle. In other 
words, each tube has a ``core''\footnote{The core of Lord Voldemort's wand was made of a phoenix feather.}: 
a distinguished line along the tube,  
which is oriented as a 1-dimensional manifold. An example is shown on the right.

Note that there are in fact four types (non-virtual) of crossings, given
by whether the core of tube A intersects the filling of B or vice versa,
and two possible directions in each case.  In the flying ring interpretation,
the orientation of the tube is the direction of the flow of time, and the 
four types of crossings represent: ring A flies through ring B from below or
from above; and ring B flies through ring A from ``below'' or from ``above''
(cf.~\cite[Exercise~\ref{1-ex:swBn}]{Bar-NatanDancso:WKO1}). The two crossings
not shown as generators are given by conjugating the classical crossing with virtual crossings, as 
in the bottom row of Figure~\ref{fig:BandCrossings}.

\parpic[r]{\input{figs/PushMembranes_2.pstex_t}}
We take the opportunity here to introduce another notation, to be called
the ``band notation'', which is more suggestive of the 4D topology than
the strand notation we have been using so far.  We represent a tube in
$\bbR^4$ by a picture of a directed band in $\bbR^3$.  By ``directed
band'' we mean that the band has a 1D direction (for example
an orientation of one of the edges). To
interpret the 3D picture of a band as an tube in $\bbR^4$, we add an extra
coordinate. Let us refer to the $\bbR^3$ coordinates as $x, y$ and $t$,
and to the extra coordinate as $z$. Think of $\bbR^3$ as being embedded
in $\bbR^4$ as the hyperplane $z=0$, and think of the band as being made
of a thin double membrane. Push the membrane up and down in the $z$
direction at each point as far as the distance of that point from the
boundary of the band, as shown on the right. This produces a tube embedded in $\bbR^4$.
In band notation, the four types of crossings appear as in Figure~\ref{fig:BandCrossings}, where underneath each crossing we indicate the corresponding
strand picture.

\begin{figure}
\input{figs/BandCrossings_2.pstex_t}
\caption{Crossings and crossing signs in band notation.}\label{fig:BandCrossings}
\end{figure}

\draftcut
\subsection{Good properties and uniqueness of the homomorphic expansion}
\label{subsec:UniquenessForTangles}

In much the same way as in the case of braids
\cite[Section~\ref{1-subsubsec:BraidCompatibility}]{Bar-NatanDancso:WKO1},
$Z$ has a number of good properties with respect to various tangle
operations: it is group-like\footnote{In practice this simply means that
the value of the crossing is an exponential.}; it commutes with adding
an inert strand (note that this is a circuit algebra operation, hence it
doesn't add anything beyond homomorphicity); and it commutes with deleting
a strand and with strand orientation reversals. All but the last of
these were explained in the context of braids and the explanations still
hold. Orientation reversal $\glos{S_k}\colon\wT\to\wT$ is the operation
which reverses the orientation of the $k$-th component. In the image of Satoh's tubing map 
this translates to reversing both the tube directions.  The induced diagrammatic operation $S_k\colon
\calA^w(T) \to \calA^w(S_k(T))$, where $T$ denotes the skeleton of a
given w-tangle, acts by multiplying each arrow diagram by $(-1)$ raised
to the power the number of arrow endings (both heads and tails) on the
$k$-th strand, as well as reversing the strand orientation.  Saying that
``$Z$ commutes with $S_k$'' means that the appropriate square commutes.

The following theorem asserts that a well-behaved homomorphic expansion of 
$w$-tangles is unique:
\begin{theorem}\label{thm:Tangleuniqueness}
The only homomorphic expansion satisfying the good properties described
above is the $Z$ defined in Section \ref{subsec:vw-tangles}.
\end{theorem}

\parpic[r]{\input{figs/rho.pstex_t}}
\begin{proof}
We first prove the following claim: Assume, by contradiction, that $Z'$ is a different 
homomorphic expansion
of $w$-tangles with the good properties described above. Let $R'=Z'(\overcrossing)$ and
$R=Z(\overcrossing)$, and denote by $\rho$ the lowest degree homogeneous
non-vanishing term of $R'-R$. (Note that $R'$ determines $Z'$, so if $Z'\neq Z$, then
$R' \neq R$.) Suppose $\rho$ is of degree $k$. 
Then we claim that $\rho=\alpha_1 w_k^1+\alpha_2 w_k^2$ is a linear combination of $w_k^1$ and $w_k^2$, 
where $w_k^i$ denotes a $k$-wheel 
living on strand $i$, as shown on the right.

Before proving the claim, note that it leads to a contradiction.
Let $d_i$ denote the operation ``delete strand $i$''.
Then up to degree $k$, we have $d_1(R')=\alpha_2 w_k^1$ and $d_2(R')=\alpha_1 w_k^2$, but
$Z'$ is compatible with strand deletions, so $\alpha_1=\alpha_2=0$. Hence
$Z$ is unique, as stated.

On to the proof of the claim, note that $Z'$ being an expansion determines the degree 1 term of $R'$ 
(namely, the single arrow 
$a^{12}$ from strand 1 to strand 2, with coefficient 1). So we can assume that $k \geq 2$. Note also that since both $R'$ and $R$ are 
group-like, $\rho$ is primitive. Hence $\rho$ is a linear combination of connected diagrams,
namely trees and wheels. 

Both $R$ and $R'$ satisfy the Reidemeister 3 relation:
$$R^{12}R^{13}R^{23}=R^{23}R^{13}R^{12}, \qquad R'^{12}R'^{13}R'^{23}=R'^{23}R'^{13}R'^{12}$$
where the superscripts denote the strands on which $R$ is placed
(compare with the proof of Theorem \ref{thm:ExpansionForTangles}).
We focus our attention on the degree $k+1$ part of the equation for $R'$,
and use that up to degree $k+1$. We can write $R'=R+\rho+\mu$, where $\mu$ denotes the degree
$k+1$ homogeneous part of $R'-R$. Thus, up to degree $k+1$, we have
$$(R^{12}\!+\!\rho^{12}\!+\!\mu^{12})(R^{13}\!+\!\rho^{13}\!+\!\mu^{13})(R^{23}\!+\!\rho^{23}\!+\!\mu^{23})=
(R^{23}\!+\!\rho^{23}\!+\!\mu^{23})(R^{13}\!+\!\rho^{13}\!+\!\mu^{13})(R^{12}\!+\!\rho^{12}\!+\!\mu^{12}).$$
The homogeneous degree $k+1$ part of this equation is a sum of some terms which contain $\rho$
and some which don't. The diligent reader can check that those which don't involve $\rho$ 
cancel on both sides, either due to the
fact that $R$ satisfies the Reidemeister 3 relation, or by simple degree counting. 
Rearranging all the terms which do involve $\rho$ to the left side, we get the following equation,
where $a^{ij}$ denotes an arrow pointing from strand $i$ to strand $j$:
\begin{equation}\label{eq:Reid3forrho}
[a^{12}, \rho^{13}]+[\rho^{12},a^{13}]+[a^{12},\rho^{23}]+[\rho^{12},a^{23}]+
[a^{13},\rho^{23}]+[\rho^{13},a^{23}]=0. 
\end{equation}

The third and fifth terms sum to $[a^{12}+a^{13},\rho^{23}]$,
which is zero due to the ``head-invariance'' of diagrams, as in Remark
\ref{rem:HeadInvariance}.

We treat the tree and wheel components of $\rho$ separately.
Let us first assume that $\rho$ is a linear combination of trees. Recall that the
space of trees on two strands is isomorphic to $\lie_2 \oplus \lie_2$, the
first component given by trees whose head is on the first strand, and the second 
component by trees with their head on the second strand.
Let $\rho=\rho_1 +\rho_2$, where $\rho_i$ is the projection to the $i$-th component
for $i=1,2$.

Note that due to $TC$, we have $[a^{12}, \rho^{13}_2]=[\rho^{12}_2,a^{13}]=
[\rho^{12}_1,a^{23}]=0$. So Equation (\ref{eq:Reid3forrho}) reduces to
$$[a^{12},\rho^{13}_1]+[\rho^{12}_1,a^{13}]+[\rho^{12}_2,a^{23}]+[\rho^{13}_1,a^{23}]+[\rho^{13}_2,a^{23}]=0$$
The left side of this equation lives in $\bigoplus_{i=1}^3 \lie_3$. Notice that only the
first term lies in the second direct sum component, while the second, third and last terms live in the third one,
and the fourth term lives in the first.
This in particular means that the first term is itself zero. By $\aSTU$, this implies 
$$0=[a^{12},\rho^{13}_1]=-[\rho_1, x_1]^{13}_2,$$
where $[\rho_1, x_1]^{13}_2$ means the tree defined by the element $[\rho_1,x_1] \in \lie_2$,
with its tails on strands 1 and 3, and head on strand 2. Hence, $[\rho_1, x_1]=0$, so $\rho_1$
is a multiple of $x_1$. The tree given by $\rho_1=x_1$ is a degree 1 element, a possibility we have eliminated, so
$\rho_1=0$.

Equation (\ref{eq:Reid3forrho}) is now reduced to  
$$[\rho^{12}_2,a^{23}]+[\rho^{13}_2,a^{23}]=0.$$
Both terms are words in $\lie_3$, but notice that the first term does not involve
the letter $x_3$. This means that if the second term involves $x_3$ at all, i.e., if
$\rho_2$ has tails on the second strand, then both terms have to be zero individually.
Assuming this and looking at the first term, $\rho^{12}_2$ is a Lie word in $x_1$ and $x_2$,
which does involve $x_2$ by assumption. We have
$[\rho^{12}_2,a^{23}]=[x_2, \rho^{12}_2]=0$, which implies $\rho^{12}_2$ is a multiple of $x_2$, in
other words, $\rho$ is a single arrow on the second strand. This is ruled out by the 
assumption that $k \geq 2$.

On the other hand if the second term does not involve $x_3$ at all, then $\rho_2$ has no tails on the second
strand, hence it is of degree 1, but again $k \geq 2$. We have proven that the ``tree part''
of $\rho$ is zero.

So $\rho$ is a linear combination of wheels. 
Wheels have only tails, so the
first, second and fourth terms of (\ref{eq:Reid3forrho}) are zero due to the tails commute relation.
What remains is $[\rho^{13}, a^{23}]=0$. We assert that this is true if and only if each
linear component of $\rho$ has all of its tails on one strand. 

To prove this, recall each wheel of $\rho^{13}$ represents a cyclic word in letters $x_1$ and $x_3$.
The map $r\colon  \rho^{13} \mapsto [\rho^{13}, a^{23}]$ is a map $\attr_2 \to \attr_3$, which sends each
cyclic word in letters $x_1$ and $x_3$ to the sum of all ways of substituting $[x_2,x_3]$ for one 
of the $x_3$'s in the word.
Note that if we expand the commutators, then all terms that have $x_2$
between two $x_3$'s cancel. Hence all remaining terms will be cyclic words in $x_1$ and $x_3$ with
a single occurrence of $x_2$ in between an $x_1$ and an $x_3$. 

We construct an almost-inverse $r'$ to $r$: for a cyclic word $w$ in $\attr_3$ with one occurrence of $x_2$,
let $r'$ be the map that deletes $x_2$ from $w$ and maps it to the resulting word in 
$\attr_2$ if $x_2$ is followed by $x_3$ in $w$, and maps it to 0 otherwise. On the rest of $\attr_3$
the map $r'$ may be defined to be 0.

The composition $r'r$ takes a cyclic word in $x_1$ and $x_3$ to itself multiplied by the number of times
a letter $x_3$ follows a letter $x_1$ in it. The kernel of this map can consist only of cyclic words 
that do not contain the sub-word $x_3x_1$, namely, these are the words of the form $x_3^k$ or $x_1^k$.
Such words are indeed in the kernel of $r$, so these make up exactly the kernel of $r$. This is exactly what 
needed to be proven: all wheels in $\rho$ have all their tails on one strand.

This concludes the proof of the claim, and the proof of the theorem. \qed
\end{proof}

%% file: figs/vDDef.pstex_t
\begin{picture}(0,0)%
\includegraphics{figs/vDDef.pstex}%
\end{picture}%
%
%
\setlength{\unitlength}{3158sp}%
\begingroup\makeatletter\ifx\SetFigFont\undefined%
\gdef\SetFigFont#1#2#3#4#5{%
  \reset@font\fontsize{#1}{#2pt}%
  \fontfamily{#3}\fontseries{#4}\fontshape{#5}%
  \selectfont}%
\fi\endgroup%
\begin{picture}(1058,292)(961,-1339)
\put(1636,-1276){\makebox(0,0)[lb]{\smash{{\SetFigFont{11}{13.2}{\rmdefault}{\mddefault}{\updefault}{\color[rgb]{0,0,0},}%
}}}}
\put(976,-1241){\makebox(0,0)[lb]{\smash{{\SetFigFont{10}{12.0}{\rmdefault}{\mddefault}{\updefault}{\color[rgb]{0,0,0}$\CA$}%
}}}}
\end{picture}%

%% file: figs/wTDef.pstex_t
\begin{picture}(0,0)%
\includegraphics{figs/wTDef.pstex}%
\end{picture}%
%
%
\setlength{\unitlength}{3158sp}%
\begingroup\makeatletter\ifx\SetFigFont\undefined%
\gdef\SetFigFont#1#2#3#4#5{%
  \reset@font\fontsize{#1}{#2pt}%
  \fontfamily{#3}\fontseries{#4}\fontshape{#5}%
  \selectfont}%
\fi\endgroup%
\begin{picture}(1121,271)(2011,-1670)
\put(2708,-1577){\makebox(0,0)[lb]{\smash{{\SetFigFont{8}{9.6}{\rmdefault}{\mddefault}{\updefault}{\color[rgb]{0,0,0}$=$}%
}}}}
\put(2026,-1561){\makebox(0,0)[lb]{\smash{{\SetFigFont{10}{12.0}{\rmdefault}{\mddefault}{\updefault}{\color[rgb]{0,0,0}$\vT$}%
}}}}
\end{picture}%

%% file: figs/taudef.pstex_t
\begin{picture}(0,0)%
\includegraphics{figs/taudef.pstex}%
\end{picture}%
%
%
\setlength{\unitlength}{3158sp}%
\begingroup\makeatletter\ifx\SetFigFont\undefined%
\gdef\SetFigFont#1#2#3#4#5{%
  \reset@font\fontsize{#1}{#2pt}%
  \fontfamily{#3}\fontseries{#4}\fontshape{#5}%
  \selectfont}%
\fi\endgroup%
\begin{picture}(2799,849)(1489,-448)
\put(3376,-361){\makebox(0,0)[b]{\smash{{\SetFigFont{10}{12.0}{\rmdefault}{\mddefault}{\updefault}{\color[rgb]{0,0,0}$+$}%
}}}}
\put(2626,239){\makebox(0,0)[b]{\smash{{\SetFigFont{10}{12.0}{\rmdefault}{\mddefault}{\updefault}{\color[rgb]{0,0,0}$\tau$}%
}}}}
\end{picture}%

%% file: figs/TubeOrientation_2.pstex_t
\begin{picture}(0,0)%
\includegraphics{figs/TubeOrientation_2.pstex}%
\end{picture}%
%
%
\setlength{\unitlength}{4934sp}%
\begingroup\makeatletter\ifx\SetFigFont\undefined%
\gdef\SetFigFont#1#2#3#4#5{%
  \reset@font\fontsize{#1}{#2pt}%
  \fontfamily{#3}\fontseries{#4}\fontshape{#5}%
  \selectfont}%
\fi\endgroup%
\begin{picture}(557,1046)(3731,-217)
\end{picture}%

%% file: figs/rho.pstex_t
\begin{picture}(0,0)%
\includegraphics{figs/rho.pstex}%
\end{picture}%
%
%
\setlength{\unitlength}{2368sp}%
\begingroup\makeatletter\ifx\SetFigFont\undefined%
\gdef\SetFigFont#1#2#3#4#5{%
  \reset@font\fontsize{#1}{#2pt}%
  \fontfamily{#3}\fontseries{#4}\fontshape{#5}%
  \selectfont}%
\fi\endgroup%
\begin{picture}(2389,944)(214,-83)
\put(676,314){\makebox(0,0)[rb]{\smash{{\SetFigFont{10}{12.0}{\rmdefault}{\mddefault}{\updefault}{\color[rgb]{0,0,0}$\rho=$}%
}}}}
\put(1726,314){\makebox(0,0)[b]{\smash{{\SetFigFont{10}{12.0}{\rmdefault}{\mddefault}{\updefault}{\color[rgb]{0,0,0}$+$}%
}}}}
\end{picture}%

%% file: foams.tex
\draftcut
\section{w-Tangled Foams} \label{sec:w-foams}

\begin{quote} \small {\bf Section Summary. }
  \summaryfoams
\end{quote}

\subsection{The Circuit Algebra of w-Tangled Foams} \label{subsec:wTF}
In the same manner as we did for tangles, we present the circuit algebra
of w-tangled foams via its Reidemeister-style diagrammatic description
accompanied by a local topological interpretation. To give a finite presentation for a circuit algebra with auxiliary (additional) operations, we use the notation
$$\CA\left\langle  \parbox{1.5in}{\centering Circuit algebra generators} \left| \parbox{1.5in}{\centering Circuit algebra relations} \right|  \parbox{1.2in}{\centering Auxiliary operations} \right\rangle.$$

\begin{definition}\label{def:wTF} Let $\glos{\wTF}$ denote the circuit algebra given by the following generators, relations and auxiliary operations:
 \[
  \wTF=\CA\!\left\langle\left.\left.
  \raisebox{-2mm}{\input{figs/wTFgensWen.pstex_t}}
  \right|
  \parbox{1.9in}{\centering \Rs, R2, R3, R4, OC, CP, FR, $W^2$, CW, TV}
  \right|
  \parbox{0.7in}{\centering $S_e, u_e, d_e $}
  \right\rangle.
\]

\end{definition}

\parpic[r]{\input{figs/VertexExamples.pstex_t}}
The generators consist of crossings, caps, wens, and foam vertices. Note that the foam vertices, where three strands meet, also come in all possible combinations of strand directions. Some additional examples are shown on the right. These generators are related by the orientation switch operation $S_e$, whose topological interpretation is explained in Section~\ref{subsubsec:wops}.

The relations \Rs, R2, R3 and $OC$ are as in Section~\ref{sec:w-tangles}. The other relations are shown and explained in the context of their local topological meaning in Section~\ref{subsubsec:wrels}. 

An {\em edge} of a w-tangled foam is a line between two vertices, tangle ends (boundary points), or caps; edges may go over and under multiple crossings. The auxiliary operations of $\wTF$ are edge orientation switches $S_e$, edge unzips $u_e$, and deletions $d_e$ of {\em long strands} which end in two tangle ends. These are described, along with their topological interpretations, in Section~\ref{subsubsec:wops}.

The circuit algebra $\wTF$ is skeleton-graded
where the circuit algebra of skeleta $\calS$ is a version of the skeleton algebra $\calS$ introduced
in Section~\ref{subsec:CircuitAlgebras}, but with vertices, caps and wens included:
 \[
  \calS=\CA\!\left.\left\langle
  \raisebox{-2mm}{\input{figs/SkelGenWen_2.pstex_t}}
  \right|
  \parbox{1.2in}{\centering $W^2$, CW, TV}
  \right\rangle.
\]
Denote by $\sigma: \wTF \to \calS$ the skeleton map, given by $\sigma(\overcrossing)=\sigma(\undercrossing)=\virtualcrossing$, and where all other generators are mapped to themselves.


\subsection{The local topology of w-tangled foams}
In this section we present the local topological meaning of  $\wTF$ generators, present the relations and show that they represent local isotopies for a space of ribbon-embedded tubes in $\bbR^4$ with caps, wens (that is, open Klein bottles), and foam vertices. We interpret the auxiliary operations as topological operations on this space.

\begin{comment}\label{com:MissingTopology}
We conjecture that the generators and relations of $\wTF$ provide a Reidemeister theory for this topological interpretation of w-tangled foams. However, there is no complete Reidemeister theorem even for w-knots (see \cite[Section 3]{Bar-NatanDancso:WKO1}). For any rigorous purposes below,
$\wTF$ is studied as a circuit algebra given by generators and relations,
with topology serving only as intuition.
\end{comment}


\subsubsection{The generators of $\wTF$}\label{subsubsec:wTFgens}
There is topological meaning to each of the generators
of $\wTF$: via a generalization of the Satoh tubing map $\delta$
of Section~\ref{sec:w-tangles} they each stand for
local features of framed knotted ribbon tubes in $\bbR^4$.

The map $\delta$ treats the strands in the same way as in Section~\ref{subsec:TangleTopology}.
The crossings are also as explained in Section
\ref{subsec:TangleTopology}: the under-strand denotes the small circle flying through a larger one,
or, equivalently, a ``thin'' tube braided through a thicker one. Recall that for tangles there are four kinds of crossings (left or right circle flying through from below or from above the other). Two of these are the generators shown, and the other two are obtained from the generators by adding virtual crossings (see Figures~\ref{fig:CrossingTubes} and~\ref{fig:BandCrossings}).

\parpic[r]{\input{figs/Cap.pstex_t}}
A bulleted end denotes a cap on the tube, or a flying circle that shrinks to a point,  
as in the figure on the right.

\vskip 5mm

The $\glos{w}$ marking on a strand indicates a
{\em wen}. A wen is a Klein bottle cut apart (see
\cite[Section~\ref{1-subsubsec:NonHorRings}]{Bar-NatanDancso:WKO1}). In the flying circle perspective, a wen represents a circle which flips over and at the same time changes its orientation; in the band persepctive, a twisted band.
\begin{center}
 \input{figs/Wen2.pstex_t}
\end{center}

\parpic[l]{\includegraphics[height=5cm]{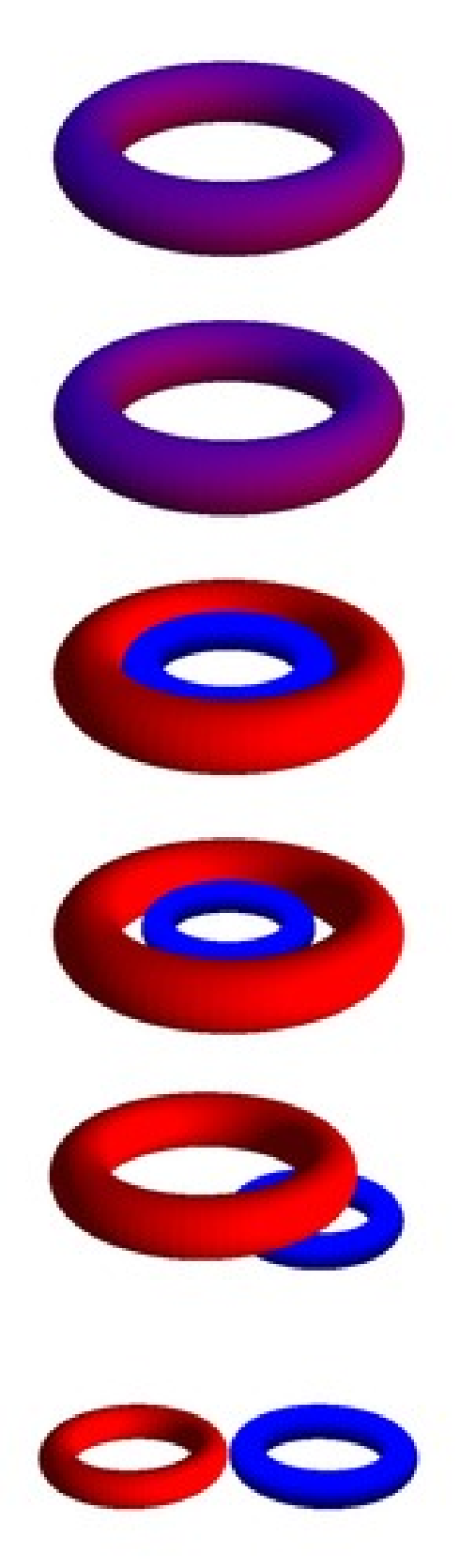}}
The final generators denote {\em singular foam vertices}. As the notation
suggests, a vertex can be thought of as a crossing with either the bottom or the top half tubes identified. To make
this precise using the flying circles interpretation, the
vertex $\raisebox{-1mm}{\input{figs/PlusVertex.pstex_t}}$ represents the movie shown on the left: the circle corresponding
to the right strand approaches the ring represented by the left strand 
from below, flies inside it, and then the two rings fuse (as opposed to
a crossing where the ring coming from the right would continue to fly
out to above and to the left of the other one).   
The second vertex $\raisebox{-1mm}{\input{figs/MinusVertex.pstex_t}}$ is
the movie where a ring splits radially into a smaller and a larger ring,
and the small one flies out to the right and below the big one. The edge corresponding to two rings identified (i.e., the top edge of $\raisebox{-1mm}{\input{figs/PlusVertex.pstex_t}}$ and the bottom edge of  $\raisebox{-1mm}{\input{figs/MinusVertex.pstex_t}}$) is called the {\em stem}. 

Vertices are rigid: the three tubes meeting at a vertex play different roles. Combinatorially, this means that the three edges meeting at a vertex are labelled (stem, inner and outer), and carry a cyclic orientation. In practice we use asymmetric pictures to avoid the clutter of labels.

As with crossings, we obtain the vertices with opposite fly-in directions by composing the generating vertices with virtual crossings, as shown
in Figure~\ref{fig:VertexTypes}. In the figure the band notation for vertices is used
the same way as it is for crossings: the fully coloured band stands for
the thin (inner) ring. 

\begin{figure}[h!]
\input{figs/VertexTypes_2.pstex_t}
\caption{Vertex types in $\wTFo$.  }\label{fig:VertexTypes}
\end{figure}


\subsubsection{The relations of $\wTF$} \label{subsubsec:wrels} 
Next, we discuss the relations of $\wTF$ and show that they represent local isotopies of w-foams.
The usual \Rs, R2, R3, and OC relations of
Figure~\ref{fig:VKnotRels} continue to apply.

The Reidemeister 4 (\glost{R4}) relations assert that a strand can be moved under
or over a vertex, as shown below. The ambiguously
drawn vertices in the figure denote a vertex of any sign with any strand directions (as in Section~\ref{subsubsec:wTFgens}). The local isotopies can be read from the band pictures in the bottom row. 
\begin{center}
 \input{figs/R4.pstex_t}
\end{center}

Recall that topologically, a cap represents a capped tube or equivalently, 
flying circle shrinking to a point. Hence, a cap
on the thin (or under) strand can be ``pulled out'' from a crossing,
but the same is not true for a cap on the thick (or over) strand, as
shown below. We
denote this relation by \glost{CP}, for Cap Pull-out. This is the case for any strands directions. 
\[ \input{figs/CapRel.pstex_t} \]

The FR, $W^2$, CW and TV relations describe the behaviour of the wens, and together we refer to them as the {\em wen relations}.

The interaction of a wen and a crossing is described by the following
Flip Relations (\glost{FR}):
\begin{center}
 \input figs/FlipRels2.pstex_t
\end{center}
To explain this in the flying circle interpretation, recall that a wen represents a circle that flips over. It does not matter whether ring B flips first and
then flies through ring A or vice versa. However, the movies in which ring A first flips and then ring B flies through it, 
or B flies through A first and then A flips differ in the fly-through direction of B through A, hence the virtual crossings. 

A double flip is isotopic to no flip, in other words two consecutive wens
are isotopic to no wen. We denote this relation by $\glos{W^2}$. 

\parpic[r]{$\input{figs/CapWen.pstex_t}$}
A cap can slide through a wen, hence a capped wen disappears, as shown
on the right, to be denoted \glost{CW}.

\vspace{3mm}

\parpic[l]{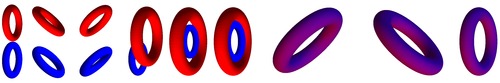}
The last wen relation describes the interaction of wens and
vertices, as illustrated on the left. In the band notation the non-filled band represents the larger circle, and the band the inner/smaller circle, as usual. Conjugating a vertex by three wens
switches the top and bottom bands, as shown in the figure on the left. Alternatively in the flying circle interpretation,
if both rings flip, then merge, and then the merged ring flips again,
this is isotopic to no flips, except the fly-in direction (from below
or from above) has changed. We denote the diagrammatic relation arising from this isotopy -- shown in the bottom right -- by
\glost{TV}, for {\em Twisted Vertex}. 

\vspace{5mm}


\subsubsection{The auxiliary operations of $\wTF$} \label{subsubsec:wops}
The circuit algebra $\wTF$ is equipped with several auxiliary operations.

The first of these is the familiar orientation switch: given an edge $e$ of a w-foam, $\glos{S_e}$ switched the direction of the edge $e$.

\parpic[r]{\input{figs/StringUnzip.pstex_t}}
The most interesting operation on w-foams is the {\em edge unzip} $\glos{u_e}$, which doubles the edge $e$ using the
blackboard framing, then attaches the ends of the doubled edge to the
connecting ones, as shown on the right. Unzip is only defined when the directions of the edges involved match, as shown on the left.  We restrict unzip to edges that are the {\em stem} of each of the two vertices they connect, and whose two vertices have their inner and outer edges aligned, as shown (in other words, edges which connect  two different generating vertices).

\begin{comment}\label{com:wTFFraming}
In order to understand the local topological meaning of the unzip operation, we need to discuss framings in more depth. Recall that
framings were mentioned in Section \ref{subsec:TangleTopology}, but have not played
a significant role so far, except to explain the lack of a Reidemeister 1
relation. 

In the local topological interpretation of $\wTF$, edges represent ribbon-knotted tubes with foam vertices, 
which are also equipped with a framing, arising from the blackboard framing of the strand diagrams via Satoh's tubing map. 
Topologically, unzip is the operation of doubling a tube by ``pushing it off itself slightly'' in the framing direction,
as shown in Figure \ref{fig:BandUnzip}.

Recall that ribbon knotted tubes have a ``filling'', with only ``ribbon''
self-intersections
\cite[Section~\ref{1-subsubsec:ribbon}]{Bar-NatanDancso:WKO1}. When we
double a tube, we want this ribbon property to be preserved.  This is
equivalent to saying that the circle obtained by pushing off any given girth
of the tube in the framing direction is not linked with the original tube,
which is indeed the case.

\parpic[r]{$\input{figs/DetailedVertex_2.pstex_t}$} 
Framings arising from the blackboard framing of strand diagrams via Satoh's tubing map 
always match at the vertices, with the normal vectors pointing either directly towards or
away from the center of the singular ring. Note that the directions of the three tubes may or may not match.
An example of a vertex with the orientations and framings shown is on the right. Note that the
framings on the two sides of each band are always mirror images of each other.
\end{comment}

\begin{figure}
\input{figs/BandUnzip_2.pstex_t}
\caption{Unzipping a tube, in band notation with orientations and framing marked.}\label{fig:BandUnzip}
\end{figure}

When a tube is unzipped, at each of the vertices at the two ends of the doubled tube there are two tubes
to be attached to the doubled tube. At each end, the normal vectors pointed either directly towards or away from the center,
so there is an ``inside'' and an ``outside'' boundary circle. The two tubes to be attached also come as an ``inside'' and an
``outside'' one, which defines which one to attach to which. An example is shown in Figure \ref{fig:BandUnzip}. 

A related operation, {\it disk unzip}, is unzip done on a capped edge, pushing the edge off in the direction of the blackboard framing, as 
before. An example in the line and band notations (with the framing suppressed) is shown below.

\begin{center}
\input{figs/CapUnzip.pstex_t}
\end{center}

Note that edges which contain wens may be unzipped by first relocating the wens to other edges or removing them, using the wen relations.

The edge deletion (denoted $d_e$) operation is restricted to ``long linear'' edges, meaning edges that do not end in a vertex or cap.


\draftcut
\subsection{The Associated Graded Structure} \label{subsec:fgrad}
Mirroring the previous section, we describe the associated graded structure
$\glos{\calA^{sw}}$ of $\wTF$ and its ``full version'' $\glos{\calA^w}$
as circuit algebras on certain generators modulo a number of
relations. From now on we will write $\glos{A^{(s)w}}$ to mean ``$\calA^{w}$
and/or $\calA^{sw}$''.
\[ \calA^{(s)w}=\CA\!\left.\left.\left\langle
  \raisebox{-2mm}{\input{figs/wTFprojgensWen.pstex_t}}
  \right|
  \parbox{2in}{\centering $\aft$, TC, VI, CP, $W^2$, TW, CW, FR, (RI for $\calA^{sw}$)}
  \right|
  \parbox{0.8in}{\centering $S_e, u_e, d_e$}
  \right\rangle.
\]

In other words, $\calA^{(s)w}$ is the circuit algebra of arrow diagrams on trivalent (foam) skeleta with 
caps and wens. That is, the skeleta are elements of $\calS$ as in Section~\ref{subsec:wTF}.  
With the exception of the first generator (the {\em arrow}), all generators are skeleton features (of degree 0). The arrow is of degree $1$. As for the generating vertices, the same remark applies as in Definition \ref{def:wTFo},
that is, vertices come in all possible combinations of edge directions.

\subsubsection{The relations of $\calA^{(s)w}$}\label{subsubsec:wTFGradRels}
In addition to the usual $\aft$ and TC 
relations (see Figure \ref{fig:TCand4T}), as well as RI in the case of $\calA^{sw}=\calA^w/RI$, arrow
diagrams in $\calA^{(s)w}$ satisfy the following additional relations:

{\it Vertex invariance}, denoted by \glost{VI}, are relations which arise
from the same principle as the classical $\aft$ relation, but with a vertex in place of a crossing:
\begin{center}
 \input figs/VI_2.pstex_t
\end{center}
The other end of the arrow is in the same place throughout the relation, somewhere outside the picture
shown. The signs are positive whenever the edge on which the arrow ends is directed towards the vertex,
and negative when directed away. The ambiguously drawn vertex means any kind of vertex, with edges oriented in any direction, as long as it is the same one throughout the relation.

\parpic[r]{\input{figs/CapHeads.pstex_t}}
The CP relation (a cap can be pulled out from under a edge but not from
over, Section \ref{subsubsec:wrels}) implies that arrow heads vanish next to a cap, as shown on the right. We denote this relation also by
\glost{CP}. (Note that an arrow tail near a cap may not vanish.)

The $W^2$, TV, and CW relations are skeleton relations introduced in Section~\ref{subsubsec:wrels} which describe the interactions of wens with each other, vertices and caps, and they continue to apply to the skeleta of arrow diagrams.

The Flip Relations FR imply that wens ``commute'' with
arrow heads, but ``anti-commute'' with arrow tails. We denote the associated graded Flip Relations also by \glost{FR}. 
\begin{center}
  \input figs/WenRel.pstex_t
\end{center}

\medskip
As in Definition~\ref{def:wJac}, we define a
``w-Jacobi diagram'' (or just ``arrow diagram'') on a foam
skeleton by allowing trivalent arrow vertices. Denote the circuit algebra of formal 
linear combinations of arrow diagrams, modulo the relations of $\calA^{(s)w}$ and the $\aSTU$ relations of Figure~\ref{fig:STU}, by $\calA^{(s)wt}$. We have the following bracket-rise theorem:

\begin{theorem}\label{thm:FoamBracketRise} The natural inclusion of diagrams induces a circuit
algebra isomorphism $\calA^{(s)w}\cong\calA^{(s)wt}$. Furthermore, the $\aAS$
and $\aIHX$ relations of Figure~\ref{fig:aIHX} hold in $\calA^{(s)wt}$.
\end{theorem}

\begin{proof} Same as the proof of Theorem~\ref{thm:BracketRise}. \qed
\end{proof}

\medskip
As in Section~\ref{subsec:vw-tangles}, the primitive elements of
$\calA^{(s)w}$ are connected diagrams (that is, connected with the skeleton removed), which are linearly generated by trees and wheels. Before
moving on to the auxiliary operations of $\calA^{(s)w}$, we make
two useful observations:

\begin{lemma}\label{lem:CapIsWheels}
$\calA^w(\raisebox{-1mm}{\input{figs/SmallCap.pstex_t}})$, the part of
$\calA^w$ with skeleton $\raisebox{-1mm}{\input{figs/SmallCap.pstex_t}}$,
is isomorphic as a vector space to the completed polynomial 
algebra freely generated by wheels $w_k$ with $k \geq 1$. Likewise
$\calA^{sw}(\raisebox{-1mm}{\input{figs/SmallCap.pstex_t}})$, except here
$k \geq 2$.
\end{lemma}

\begin{proof}
 Any arrow diagram with an arrow head at its top is zero by the Cap Pull-out (CP) relation. If $D$ is an arrow
diagram that has a head somewhere on the skeleton but not at the top, then one can use repeated $\aSTU$ relations 
to commute the head to the top at the cost of diagrams with one fewer skeleton head. 

Iterating
this procedure, we can get rid of all arrow heads, and hence write $D$ as a linear combination of 
diagrams having no heads on the skeleton. All connected components of such diagrams are wheels. 

To prove that there are no relations between wheels in $\calA^{(s)w}(\raisebox{-1mm}{\input{figs/SmallCap.pstex_t}})$, 
let $S_L\colon \calA^{(s)w}(\uparrow_1) \to \calA^{(s)w}(\uparrow_1)$ 
(resp. $S_R$) be the map that sends an arrow diagram to the sum of all ways of dropping one left (resp. right) arrow 
(on a vertical edge, left means down and right means up). Define
$$\glos{F}:=\sum_{k=0}^{\infty}\frac{(-1)^k}{k!}D_R^k(S_L+S_R)^k,$$
where $D_R$ is a short right arrow.
We leave it as an exercise for the reader to check that $F$ is a bi-algebra homomorphism that kills diagrams with an arrow head at the top
(i.e., CP is in the kernel of $F$), and $F$ is injective on wheels. This concludes the proof.
\qed
\end{proof}

\begin{lemma}
$\calA^{(s)w}(\raisebox{-1mm}{\input{figs/PlusVertex.pstex_t}})=\calA^{(s)w}(\uparrow_2)$, where $\calA^{(s)w}(\raisebox{-1mm}{\input{figs/PlusVertex.pstex_t}})$
stands for the space of arrow diagrams whose skeleton is a vertex of any type, with
any orientation of the edges, and $\calA^{(s)w}(\uparrow_2)$
denotes the space of arrow diagrams on two strands.
\end{lemma}

\begin{proof}
 Use the vertex invariance (VI) relation to push all arrow heads and tails from the ``trunk'' of the vertex to the other two edges.
\qed
\end{proof}

\subsubsection{The auxiliary operations of $\calA^{(s)w}$}\label{subsubsec:wTFGradOps}
Recall from Section \ref{subsec:UniquenessForTangles} that the associated graded orientation switch  operation  $S_e\colon \calA^{(s)w}(s) \to \calA^{(s)w}(S_e(s))$ acts by reversing the direction of the skeleton edge $e$, and multiplying  each arrow diagram by $(-1)$ raised to the number of arrow endings on $e$ (counting both heads and tails).

\parpic[r]{\input{figs/Unzip.pstex_t}}
The arrow diagram operations induced by unzip and disc unzip $u_e\colon  \calA^{(s)w}(s) \to \calA^{(s)w}(u_e(s))$ are both denoted $u_e$, and interpreted appropriately according to whether the 
edge $e$ is capped. They both map each arrow ending (head or tail) on $e$ to
a sum of two arrows, one ending on each of the new edges, as shown on the right. In other words, if in a primitive arrow diagram $D$ there are $k$ arrow
ends on $e$, then $u_e(D)$ is a sum of $2^k$ primitive arrow diagrams.

The operation induced by deleting the long linear strand $e$ is the map $d_e\colon  \calA^{(s)w}(s) \to \calA^{(s)w}(d_e(s))$ which kills arrow diagrams with
any arrow ending (head or tail) on $e$, and leaves all else unchanged, except with $e$ removed. 


\draftcut
\subsection{The homomorphic expansion}\label{subsec:wTFExpansion}
If a homomorphic expansion for $\wTF$ exists, it is determined by the values of the generators, as $\wTF$ has a finite presetantion. We are interested in particular in group-like homomorphic expansions, where the values of the generators are exponentials of (infinite series of) primitive arrow diagrams. For more detail see \cite[Section 2.5.1.2]{Bar-NatanDancso:WKO1}. 

We will see that the value of the vertex $V\in \calA^{sw}(\raisebox{-1mm}{\input{figs/PlusVertex.pstex_t}})\cong\calA^{sw}(\uparrow_2)$ plays a particularly important role. 
It will also  become clear that the short arrows of $V$ can be ignored: here a short arrow means an arrow on a single edge of the vertex, whose head and tail are adjacent on the skeleton. The reason is that, given a homomorphic expansion $Z$ with $Z(\raisebox{-1mm}{\input{figs/PlusVertex.pstex_t}})=V=e^v$, and $a$ is a short arrow on the vertex, then changing $V$ to $V'=e^{v+a}$ defines another homomorphic expansion $Z'$. We explain this in more detail later; for now we make the following definition:

\begin{definition} 
A homomorphic expansion is {\em v-small} if $Z(\raisebox{-1mm}{\input{figs/PlusVertex.pstex_t}})=e^v$ where $v$ is a (possibly infinite) linear combination of primitive arrow diagrams which does not include short arrows in degree one.
\end{definition}

Given a homomorphic expansion $Z:\wTF \to \calA^{sw}$, denote by $W$ the $Z$-value of the wen. We are now able to state one of the main theorems of this paper:

\begin{theorem}\label{thm:WenATEquivalence}
Group-like\footnote{The
formal definition of the group-like property is along the lines of
\cite[Section~\ref{1-par:Delta}]{Bar-NatanDancso:WKO1}. In practice, it means
that the $Z$-values of the vertices, crossings, and cap (denoted $V$,
$R$ and $C$ below) are exponentials of linear combinations of
connected diagrams.}
 homomorphic expansions $Z: \wTF \to \calA^{sw}$ exist, and these which are v-small and satisfy $W=1$ are in one-to one correspondence 
with solutions to the Kashiwara-Vergne equations (defined in Section~\ref{subsec:EqWithAT}) with even Duflo function. \qed
\end{theorem}

Our goal is to explain and prove this theorem. To begin, observe that
finding a homomorphic expansion $Z: \wTF \to \calA^{sw}$ is equivalent to finding values for the generators
of $\wTF$ in $\calA^{sw}$, so that these values satisfy the equations
which arise from the relations in $\wTF$ and the homomorphicity
with respect to the auxiliary operations. 
In this subsection we derive these equations; in Section~\ref{subsec:EqWithAT} we show
that they are equivalent to the Alekseev-Torossian version of the
Kashiwara-Vergne equations \cite{AlekseevTorossian:KashiwaraVergne} with even Duflo function. In
\cite{AlekseevEnriquezTorossian:ExplicitSolutions} Alekseev Enriquez
and Torossian construct explicit solutions to these equations using
associators. In~\cite{Bar-NatanDancso:WKO3} we will interpret and independently prove this
result in the context of homomorphic expansions for w-tangled foams.

First we set notation for the images of the most important generators. Assume that $Z$ is a homomorphic expansion. Let $\glos{R}:=Z(\overcrossing) \in \calA^{sw}(\uparrow_2)$. 
Let
$\glos{C}:=Z(\raisebox{-1mm}{\input{figs/SmallCap.pstex_t}})\in
\calA^{sw}(\raisebox{-1mm}{\input{figs/SmallCap.pstex_t}})$.
By Lemma \ref{lem:CapIsWheels}, we know that
$C$ is made up of wheels only.  Let $\glos{W}\in \calA^{sw}(\uparrow)$ denote the $Z$-value of the wen, and we adopt the convention that $W$ is always placed on the skeleton edge after the wen. Finally, let
$\glos{V}=\glos{V^+}:=Z(\raisebox{-1mm}{\input{figs/PlusVertex.pstex_t}})\in
\calA^{sw}(\raisebox{-1mm}{\input{figs/PlusVertex.pstex_t}})\cong
\calA^{sw}(\uparrow_2)$, and
$\glos{V^-}:=Z(\raisebox{-1mm}{\input{figs/MinusVertex.pstex_t}})\in
\calA^{sw}(\raisebox{-1mm}{\input{figs/MinusVertex.pstex_t}})\cong
\calA^{sw}(\uparrow_2)$.

We first address the value of the wen. 
Recall that the FR relation in $\calA^{sw}$ states that skeleton wens commute with arrow heads and anti-commute with arrow tails. For a primitive arrow diagram $D$ we denote
$$\overline{D}:= (-1)^{\#\{ \text{arrow tails in }D\}} D.$$
Then we have that $wD=\overline{D}w$, where $w$ denotes a skeleton wen and $D\in \calA^{sw}(\uparrow)$. 
The same equality holds in $\calA^{sw}(\uparrow_n)$ if all strands of $D$ are commuted with wens.

\begin{lemma}\label{lem:ValueW}
Under any group-like homomorphic expansion $Z$ the value of the wen can be expressed as $W=\exp({\sum_{k=1}^\infty c_{2k+1}w_{2k+1}})$, where $c_{2k+1}$ are constants, and $w_{2k+1}$ are odd wheels, and any such $W$ satisies the equations induced by the $W^2$ relation and homomorphicity with respect to $S_e$.
\end{lemma}

\begin{proof}
From the $W^2$ relation we obtain that $\overline{W}W=1$, see Figure~\ref{fig:ZofWen2}. 
Since $Z$ is group-like, we have $W=e^\omega$ for some primitive $\omega$. Since $\omega$ is a primitive element of $\calA^{sw}(\uparrow)$, by the description of primitive arrow diagrams in Section~\ref{subsec:ATSpaces} it can be written a sum of wheels in degrees 2 and above, with possibly a multiple of a single arrow in degree 1. (Higher degree trees on a single strand reduce to wheels by the AS and STU relations.) Write $\omega= p_1a + \sum_{k=1}^{\infty} p_{2k}w_{2k} + \sum_{l=1}^\infty p_{2l+1}w_{2l+1}$, where $a$ denotes the degree 1 arrow, $w_i$ are $i$-wheels, and $p_i$ are constants.  

Then $\overline{W}=e^{\overline{\omega}}$, and $\overline{\omega}= - p_1a + \sum_{k=1}^{\infty} p_{2k}w_{2k} - \sum_{l=1}^\infty p_{2l+1}w_{2l+1}$. Thus, $\overline{W}=W^{-1}$ means that $p_{2k}=0$ for all $k\geq 1$, in other words, $W$ is contained in odd degrees only.

From the homomorphicity of $Z$ with respect to orientation switches we further see that $S(W)=\overline{W'}$, where $S(W)$ denotes the orientation switch of $W$. Combining this with the previous result, we have $S(W)=\overline{W}$, which further implies that $p_1=0$, completing the proof. \qed
\end{proof}

\begin{figure}
\input{figs/ZofWen2_2.pstex_t}
\caption{The implication of the $W^2$ relation.}\label{fig:ZofWen2}
\end{figure}

\medskip
Recall that by convention we number strands at the bottom of each diagram from left to right, and for an arrow diagram $D\in \calA^{sw}(\uparrow_k)$, $D^{i_1i_2\ldots i_k}$ means ``$D$ placed on strands $i_1,\ldots,i_k$. For instance, $R^{23}$ means ``$R$ placed on strands 2 and 3''. In this section 
we also need to use co-simplicial notation, for example $R^{(23)1}$ means ``$R$ with its first strand doubled (unzipped), then placed on strands 2, 3 and 1''.

We recall the following result Sections
\ref{subsec:vw-tangles} and \ref{subsec:UniquenessForTangles}:

\begin{lemma}\label{lem:ValueR}
For any homomorphic expansion $Z$, the values of the crossings are as follows: $R=Z(\overcrossing)=e^a$ where $a$ denotes a single arrow from the over to the under strand, $Z(\undercrossing)=(R^{-1})^{21}=e^{-a^{21}}$, where again $a^{21}$ points from the over strand to the under strand. These values satisfy the equations induced by $R1^s$, $R2$, $R3$ and $OC$ in $\calA^{sw}$. \qed
\end{lemma}

One of the important restrictions on the value $V$ arises from the R4 relations:

\begin{lemma}\label{lem:R4}
The R4 relations induce the following single equation on $V$ and $R$:
\begin{equation}\tag{R4}\label{eq:HardR4}
 V^{12}R^{(12)3}=R^{23}R^{13}V^{12}.
\end{equation}
\end{lemma}

\begin{proof}
The Reidemeister 4
relation with a strand over a vertex induces an equation that is automatically satisfied, as follows:
\begin{center}
 \input{figs/R4ToEquation.pstex_t} 
\end{center}
In other words, the over strand $R4$ relation induces the equation
$$V^{12}R^{3(12)}=R^{32}R^{31}V^{12}.$$
However, observe that by the ``head-invariance'' property of arrow diagrams (Remark \ref{rem:HeadInvariance})
$V^{12}$ and $R^{3(12)}$ commute on the left hand side. Hence the left hand side equals $R^{3(12)}V^{12}=R^{32}R^{31}V^{12}$.
Also, $R^{3(12)}=e^{a^{31}+a^{32}}=e^{a^{32}}e^{a^{31}}=R^{32}R^{31}$, where the second step is an application of the TC relation ($a^{31}$ and $a^{32}$ commute). Therefore, this equation is true regardless of the choice of $V$.

We have no such luck with the second Reidemeister 4 relation, which, in the same manner as above, translates to the \eqref{eq:HardR4} equation $V^{12}R^{(12)3}=R^{23}R^{13}V^{12}$.
There is no ``tail invariance'' of arrow diagrams, so $V$ and $R$ do not commute on the left hand side; also, heads do not commute and so $R^{(12)3}\neq R^{23}R^{13}$. 
Thus, this equation places a genuine restriction on the choice of $V$. \qed
\end{proof}

\medskip

\begin{lemma}\label{lem:CP}
The equation induced by the CP relation is automatically satisfied for any choice of $C$.
\end{lemma}
\begin{proof}
The Cap Pull-out (CP) relation translates to the equation $R^{12}C^2=C^2$. By head invariance,
$R^{12}C^2=C^2R^{12}$. Now $R^{12}$ is just below the cap on strand $2$, and thus by the CP relation in $\calA^{sw}$, every term of $R^{12}$
with an arrow head at the top of strand $2$ is zero. Hence, the only surviving term of $R^{12}$ is $1$ (the empty diagram), which makes the
equation true. \qed
\end{proof}

\begin{lemma}\label{lem:VV-}
If $V=Z(\raisebox{-1mm}{\input{figs/PlusVertex.pstex_t}})$, and $V^-=Z(\raisebox{-1mm}{\input{figs/MinusVertex.pstex_t}})$, then $V^-=V^{-1}$.
\end{lemma}

\begin{proof}
This is an immediate consequence of the homomorphicity of $Z$ with respect to the unzip operation. \qed
\end{proof}
\medskip

For the value of the cap denote $C=e^c$, where $c=\sum_{j=1}^\infty r_j w_j$, with $r_j$ constants and $w_j$ the $j$-wheel. The value of the cap is the product of even and odd parts, that is, $C=C_{eve}C_{odd}$, where $C_{eve}=e^{c_{eve}}$ with $c_{eve}=\sum_{k=1}^\infty r_{2k} w_{2k}$, and  $C_{odd}=e^{c_{odd}}$ with $c_{odd}=\sum_{l=1}^\infty r_{2k} w_{2k}$

\begin{lemma}\label{lem:CW}
The equation induced by the CW relation is $c_{odd}=-\frac{1}{2}\omega$, or equivalently $C_{odd}=W^{-1/2}$.
\end{lemma}

\begin{proof}
Applying $Z$ to each side of the CW relation, we obtain $\overline{W}\overline{C}=C$ in $\calA^{sw}(\upcap)$. Substituting the formulas for $W$ and $C$, the statement follows.\qed
\end{proof}

\medskip

Possibly the most interesting equation is the one induced by the twisted vertex relation. For this we introduce one additional piece of notation. Given $D\in \calA^{sw}(\uparrow_n)$, denote $\glos{D^*}:=S_1S_2\cdots S_n(\overline{D})$, and call this the {\em adjoint} of $D$. In other words, the operation 
$*: \calA^{sw}(\uparrow_n) \to \calA^{sw}(\uparrow_n)$ 
reverses the edge directions and multiplies an arrow diagram $D$ by 
$(-1)^{\#{\{\text{arrow heads on the skeleton}}\}}$ .

\begin{lemma}\label{lem:WenUnitarity}
The TV relation is induces the ``Wen-Unitarity'' equation   
\begin{equation}\tag{WU}\label{eq:WU}
(W^{-1})^{(12)} V^*W^1W^2V=1.
\end{equation}
\end{lemma}

\begin{proof}
We apply $Z$ to each side of the TV relation, as shown in Figure~\ref{fig:TVUnitarity}.
On the right hand side of the relation is a vertex $\raisebox{-1mm}{\input{figs/MinusVertex.pstex_t}}$ with the edge orientations reversed, upside down and the edges numbered $(2,1)$ as the vertex follows a virtual crossing. Therefore, the value of this vertex is $S_1S_2(V^{-1})$ by Lemma~\ref{lem:VV-}. 

The top wen value $W$ on the far left side can be ``pulled down'' to the bottom two edges using the VI relation. Therefore, we obtain the following equation in $\calA^{sw}(\uparrow_2)$:
$$\overline{W}^1\overline{W}^2 \overline{V}W^{(12)}=S_1S_2(V^{-1})$$

Applying $S_1S_2$ to both sides, and using from Lemma~\ref{lem:ValueW} that $\bar{W}=W^{-1}=S(W)$, we obtain the equation \eqref{eq:WU}. \qed
\end{proof}
\medskip

\begin{figure}
{
  \def\SSV{$S_1S_2(V^{-1})$}
  \def\OV{$\overline{V}$}
  \def\OW{$\overline{W}$}
  \begin{center} \input{figs/TVUnitarity_2.pstex_t} \end{center}
}
\caption{Applying $Z$ to each side of the TV relation.}\label{fig:TVUnitarity}
\end{figure}

\begin{lemma}\label{lem:CapEq} 
Homomorphicity of $Z$ with respect to the disc unzip operation is equavialent to the Cap Equation:
\begin{equation}\tag{C}\label{eq:CapEqn}
  V^{12}C^{(12)}=C^1C^2 \qquad\text{in}\quad
  \calA^{sw}(\raisebox{-1mm}{\input{figs/SmallCap.pstex_t}}_2)
\end{equation}
\end{lemma}

\parpic[r]{\input{figs/CapEqn.pstex_t}}\picskip{5}
\noindent {\em Proof.} We need to apply $Z$ and the cap unzip $u$ in either order to the w-foam shown in the figure on the right. On the left hand side, the value of the cap is unzipped and gives $C^{(12)}$. Note that \eqref{eq:CapEqn} is an equation in $\calA^{sw}(\upcap_2)$. \qed

\vspace{1cm}

To summarize, we have proven the following theorem:

\begin{theorem}\label{thm:ZWenEquations}
 $Z:\wTF \to \calA^{sw}$ is a group-like homomorphic expansion if and only if the values of the key generators $R=e^a$, $V$, $W$ and $C$ are group-like, and satisfy the equations
\begin{enumerate}
\item[(R4)]  $V^{12}R^{(12)3}=R^{23}R^{13}V^{12}$ in $\calA^{sw}(\uparrow_3)$,
\item[(WU)] $(W^{-1})^{(12)} V^*W^1W^2V=1$ in $\calA^{sw}(\uparrow_2)$,
\item[(C)] $V^{12}C^{(12)}=C^1C^2$ in
  $\calA^{sw}(\upcap_2)$,
\item[(CW)] $C_{odd}=W^{-1/2}$ in $\calA^{sw}(\upcap_1)$, or equivalently, as power series in odd wheels.
\end{enumerate}
\end{theorem}


\draftcut
\subsection{The equivalence with the Alekseev-Torossian equations}
\label{subsec:EqWithAT}
First let us recall Alekseev and Torossian's
formulation of the generalized Kashiwara-Vergne problem
(see~\cite[Section~5.3]{AlekseevTorossian:KashiwaraVergne}):

{\bf Generalized KV problem:} Find an element $\glos{F}\in \TAut_2$ with the properties
\begin{equation}\label{eq:ATKVEqns}
 F(x+y)=\log(e^xe^y), \text{ and } j(F)\in \im(\tilde{\delta}).
\end{equation}
Here $\tilde{\delta}\colon  \attr_1 \to \attr_2$ is defined by $(\tilde{\delta}a)(x,y)=a(x)+a(y)-a(\log(e^{x}e^{y}))$,
where $\attr_2$ is generated by cyclic words in the letters $x$ and $y$. (See
\cite{AlekseevTorossian:KashiwaraVergne}, Equation (8)). Note that an element of $\attr_1$ is a power series in one variable 
with no constant term, called the {\em Duflo function}. 
In other words, the second condition says that there exists 
$a \in \attr_1$ such that $jF=a(x)+a(y)-a(\log(e^{x}e^{y}))$.

\medskip

\noindent{\em Proof of Theorem~\ref{thm:WenATEquivalence}.}
We need to translate the equations of Theorem~\ref{thm:ZWenEquations} to equations in the Alekseev-Torossian spaces, using the identifications of Proposition~\ref{prop:Pnses} and the identification of wheels with cyclic words. Note the condition in Theorem~\ref{thm:WenATEquivalence} that $W=1$. With this simplification the (CW) equation simply asserts that the value $C$ is an even power series in wheels. The \eqref{eq:WU} equation simplifies to the following, which we call the {\em Unitarity} of the vertex:
\begin{equation}\tag{U}\label{eq:unitarity}
V^*V=1.
\end{equation}

Recall from Section \ref{subsec:ATSpaces} that the map $u\colon  \tder_2 \to \calA^{sw}(\uparrow_2)$ plants the head of a tree
above all of its tails.
Suppose that the values $V$ and $C$ satisfy the simplified equations of Theorem~\ref{thm:ZWenEquations} with $W=1$. Write  $V=e^be^{uD}$, where $b \in \tr_2^s$, $D \in  \tder_2\oplus \fraka_2$, 
and where  $V$ can be written in this form without loss of generality because wheels can always be commuted to the bottom of a diagram (at the possible cost of more wheels). Furthermore, $V$ is group-like and hence it can be written in exponential form. Similarly, write $C=e^c$ with $c \in \attr_1^s$.

Note that $u(\fraka_2)$ is central in $\calA^{sw}(\uparrow_2)$ and
that replacing a solution $(V,C)$ by $(e^{u(a)}V, C)$ for any $a \in
\fraka_2$ does not interfere with any of the equations (\ref{eq:HardR4}),
(\ref{eq:unitarity}) or (\ref{eq:CapEqn}). Hence we may assume that $D$
does not contain any single arrows, that is, $Z$ is v-small and $D \in \tder_2$. Also, a
solution $(V,C)$ in $\calA^{sw}$ can be lifted to a solution in $\calA^w$
by simply setting the degree one terms of $b$ and $c$ to be zero. It is
easy to check that this $b \in \attr_2$ and $c \in \attr_1$ along with $D$
still satisfy the equations. (In fact, in $\calA^w$ (\ref{eq:unitarity})
and (\ref{eq:CapEqn}) respectively imply that $b$ is zero in degree 1,
and $c$ is already assumed to be even.)  In light of this we declare that $b\in \attr_2$ and $c
\in \attr_1$.

The hard Reidemeister 4 equation (\ref{eq:HardR4}) reads $V^{12}R^{(12)3}=R^{23}R^{13}V^{12}$. Denote the arrow from strand 1 to strand 3 by $x$, and the
arrow from strand 2 to strand 3 by $y$. Substituting the known value for $R$ and rearranging, we get 
$$e^be^{uD}e^{x+y}e^{-uD}e^{-b}=e^ye^x.$$ Equivalently, $e^{uD}e^{x+y}e^{-uD}=e^{-b}e^ye^xe^b.$ Now on the right side there are only tails on the
first two strands, hence $e^b$ commutes with $e^ye^x$, so $e^{-b}e^b$ cancels. Taking logarithm of both sides we obtain 
$e^{uD}(x+y)e^{-uD}=\log e^ye^x$. Now for notational alignment with \cite{AlekseevTorossian:KashiwaraVergne} we switch strands 1 and 2, which exchanges 
$x$ and $y$ so we obtain:
\begin{equation}\label{eq:HardR4Translated}
e^{uD^{21}}(x+y)e^{-uD^{21}}=\log e^xe^y.
\end{equation}

The unitarity of $V$ (Equation (\ref{eq:unitarity})) translates to $1=e^be^{uD}(e^be^{uD})^*,$ where $*$ denotes the adjoint map (Definition \ref{def:Adjoint}). Note that the adjoint switches 
the order of a product and acts trivially on wheels. Also, $e^{uD}(e^{uD})^*=J(e^D)=e^{j(e^D)}$, by Proposition \ref{prop:Jandj}. 
So we have $1=e^be^{j(e^D)}e^b$. Multiplying by $e^{-b}$ on the right and by $e^b$ on the left, we get $1=e^{2b}e^{j(e^D)}$, and again by switching strand 1 and 2 we arrive at
\begin{equation}\label{eq:UnitarityTranslated}
1=e^{2b^{21}}e^{j(e^{D^{21}})}.
\end{equation}

As for the cap equation, if $C^1=e^{c(x)}$ and $C^2=e^{c(y)}$, then $C^{12}=e^{c(x+y)}$. Note that wheels
on different strands commute, hence $e^{c(x)}e^{c(y)}=e^{c(x)+c(y)}$, so the cap equation reads $$e^be^{uD}e^{c(x+y)}=e^{c(x)+c(y)}.$$ As this equation lives
in the space of arrow diagrams on two \emph{capped} strands, it remains unchanged if we multiply the left side on the right by $e^{-uD}$: $uD$ has its head at the top, so it
is 0 by the Cap relation, hence $e^{uD}=1$ near the cap. Hence, $$e^be^{uD}e^{c(x+y)}e^{-uD}=e^{c(x)+c(y)}.$$ 

\parpic[r]{\input{figs/Sigma.pstex_t}}
On the right side of the equation above \linebreak $e^{uD}e^{c(x+y)}e^{-uD}$ reminds us of Equation (\ref{eq:HardR4Translated}), however we cannot use (\ref{eq:HardR4Translated})
directly as we are working in a different space now. In particular, $x$ there meant an arrow from strand 1 to strand 3, while here it means a one-wheel on (capped) 
strand 1, and similarly for $y$. Fortunately, there is a linear map $\sigma\colon  \calA^{sw}(\uparrow_3) \to \calA^{sw}(\raisebox{-1mm}{\input{figs/SmallCap.pstex_t}}_2)$,
where $\sigma$ ``closes the third strand and turns it into a chord (or internal) strand, and caps the first two strands'', as shown on the right. This map is
well defined (in fact, it kills almost all relations, and turns one $\aSTU$ into an $\aIHX$). Under this map, using our abusive notation, $\sigma(x)=x$ and 
$\sigma(y)=y$.

Now we can apply Equation (\ref{eq:HardR4Translated}) to get $e^{uD}e^{c(x+y)}e^{-uD}=e^{c(\log e^y e^x)}$. 
Substituting this into the cap equation we obtain 
$e^be^{c(\log e^y e^x)}=e^{c(x)+c(y)}$, which, using that tails commute, implies
$b=c(x)+c(y)-c(\log e^y e^x)$. Switching strands 1 and 2, we obtain
\begin{equation}\label{eq:CapEqnTranslated}
b^{21}=c(x)+c(y)-c(\log e^x e^y)
\end{equation}

In summary, we can use $(V,C)$ to produce $F:=e^{D^{21}}$ (sorry\footnote{%
  We apologize for the annoying $2\leftrightarrow 1$ transposition in this equation,
  which makes some later equations, especially~\eqref{eq:ATPhiandV},
  uglier than they could have been. There is no depth here, just
  mis-matching conventions between us and Alekseev-Torossian.
})
and $a:=-2c$ which satisfy the Alekseev-Torossian equations
(\ref{eq:ATKVEqns}), as follows: $e^{D^{21}}$ acts on $\lie_2$ by conjugation
by $e^{uD^{21}}$, so the first part of (\ref{eq:ATKVEqns})
is implied by (\ref{eq:HardR4Translated}). The second half of
(\ref{eq:ATKVEqns}) is true due to (\ref{eq:UnitarityTranslated}) and
(\ref{eq:CapEqnTranslated}).

On the other hand, suppose that we have found $F\in \TAut_2$ and even Duflo function $a \in \tr_1$ satisfying (\ref{eq:ATKVEqns}). 
Then set $D^{21}:=\log F$, $b^{21}:=\frac{-j(e^{D^{21}})}{2}$,
and $c \in \tilde{\delta}^{-1}(b^{21})$, in particular $c=-\frac{a}{2}$ works. Then $V=e^be^{uD}$ and the even cap value $C=e^c$ satisfy the equations for 
homomorphic expansions (\ref{eq:HardR4}), (\ref{eq:unitarity})
and (\ref{eq:CapEqn}), and hence define a homomorphic expansion of $\wTF$ with $W=1$.

Furthermore, these maps between solutions of the KV problem and nv-small homomorphic expansions for $\wTF$ with $W=1$ are
obviously inverses of each other, and hence they provide a bijection between these sets as stated.\qed

\begin{remark}
The fact that $Z$ can be chosen to have $W=1$ and $C$ even
follows from Proposition 6.2 of \cite{AlekseevTorossian:KashiwaraVergne}. In Proposition 6.2 Alekseev and Torossian 
show that the even part of $f$ is $\frac{1}{2}\frac{\log(e^{x/2}-e^{-x/2})}{x}$, and that for any $f$ with this even part  (and any odd part) there exists a corresponding solution $F$
of the generalized $KV$ problem. In particular, $f$ can be assumed to be even, and hence it can be guaranteed that $C$ consists of even wheels only.
\end{remark}

\subsection{Orientable w-tangled foams}\label{subsec:OriFoams}

There is a sub-circuit algebra of $\wTF$ consisting of the w-tangled foams which contain no wens. We call this the circuit algebra of orientable w-foams, and denote it by $\wTFo$. (These foams can be equipped with a global surface orientation, which induces crossing and vertex signs consistent with the signs suggested by the diagrams. However, this is not necessary.)

\begin{lemma}\label{lem:CancelWens}
Let $F\in \wTF$ be a w-foam with the property that there are an 
even number of wens along any path connecting two tangle ends, and
along any cycle in $F$. Then all of the wens in $F$ cancel by the wen
relations. Furthermore, the process of cancelling all wens can be made
canonical by a choice -- for each connected component of the skeleton
of $F$ -- of a spanning tree $T$, and a basepoint on $T$, which is a
tangle end if there are any. 
\end{lemma}

\begin{proof}
First note that the statement of the lemma concerns only the skeleton 
$\sigma(F)$: by the FR relations wens slide through crossings, at
the possible cost of more virtual crossings. The skeleton of $F$ is
a uni-trivalent graph whose univalent ends are either caps of tangle
ends. Due to the CW relation, capped edges can be ignored, that is,
deleted without loss of generality. Thus, assume that $\sigma(F)$ is a
uni-trivalent graph in which all univalent vertices are tangle ends.

Given a choice of spanning tree $T$ for $\sigma(F)$ and a base point 
on it, there is a unique way to ``clear $T$ of wens''. Namely, use the
TV relation to push wens off of $T$ away from the base point. The $TV$
relation does not change the parity of the number of wens along any cycle, 
or any path connecting two tangle ends. At the end of this process, 
all wens will end up either on an edge of $\sigma(F)$ not in $T$, or at
a tangle end (which are all necessarily in $T$).

At the end of the process, there is still an even number of wens 
on the path from any given tangle end to the base point (which is
also necessarily a tangle end in this case), and so there is an even
number of wens at each tangle end, therefore they cancel by the $W^2$
relation. For any non-$T$ edge $e$ of $\sigma(F)$, there is a unique
path $\gamma$ in $T$ which connects the two ends of $e$. Since there
originally was an even number of wens along the cycle $e \cup \gamma$,
there is an even number of wens on $e$ at the end of the process, which
therefore cancel. \qed
\end{proof}
\medskip

We derive a generators - relations - operations presenation for $\wTFo$. Since the wen is no longer a generator, there are no wen relations. The operations $S_e$, $u_e$ and $d_e$ restrict to $\wTFo$. 
The composition with wens in $\wTF$ induces an involution on $\wTF^o$: while wens are not included in $\wTFo$, composition with a wen at {\em every} tangle end is well-defined:

\begin{definition}\label{def:wenjugation}
For $F\in \wTFo$, consider $F$ as an element of $\wTF$, and let $\bar{F}$ denote $F$ composed with a wen at every tangle end. Then by Lemma~\ref{lem:CancelWens}, $\bar{F}\in \wTFo$. We call this operation {\em wenjugation}, and denote it by $\glos{-}:\wTFo \to \wTFo$. 
\end{definition}

\begin{definition}\label{def:wTFo} The circuit algebra of oriented w-foams is defined by the presentation
\[
  \wTFo=\CA\!\left.\left.\left\langle
  \raisebox{-2mm}{\input{figs/wTFgens.pstex_t}}
  \right|
  \parbox{1.9in}{$R1^s$, R2, R3, R4, OC, CP} 
  \right|
  \parbox{0.9in}{$S_e,u_e,d_e, -$}
  \right\rangle.
\]
\end{definition}

Next, we verify that $\wTFo$ -- as defined by the presentation above -- injects into $\wTF$. In other words, the generators and relations description above is indeed a description of sub-circuit algebra of $\wTF$ generated by all orientable (non-wen) $\wTF$ generators. 

\begin{proposition}
The circuit algebra of oriented w-foams $\wTFo$ injects into $\wTF$.
\end{proposition}

\begin{proof}
We need to show that given $F, F' \in \wTFo$ for which $F~F'$ via a sequence of $\wTF$ relations, then $F~F'$ also in $\wTFo$. This can be verified explicitly, as follows. Choose a spanning tree and base point for each connected component of $\sigma(F)=\sigma(F')$. Let $F=F_0\sim F_1\sim \cdots \sim F_n=F'$ be a sequence of $\wTF$ moves. Since all $F_i$ are $\wTF$-equivalent to $F\in \wTFo$, they all satisfy the conditions of Lemma~\ref{lem:CancelWens}. Via the process of Lemma~\ref{lem:CancelWens}, each $F_i$ (i=0,\ldots, n) is canonically equivalent to an element of $\wTFo$, call this element $\Omega(F_i)$. Hence, we only need to show that $\Omega(F_i)\sim \Omega(F_{i+1})$ in $\wTFo$, where $F_i$ and $F_{i+1}$ differ in a single relation in $\wTF$. This is obvious if that relation is not a wen relation, easy for the $W^2$ and $CW$ relations, and directly verified with some effort for the $FR$ and $TV$ relations.\qed
\end{proof}
\medskip

The circuit algebra $\wTFo$ is again skeleton graded, with skeleton circuit algebra given by
\[
  \calS^0 = \CA\!\left\langle\raisebox{-2mm}{\input{figs/SkelGen.pstex_t}}\right\rangle
\]

The associated graded structure -- which we continue to denote by $\calA^{sw}$ to avoid too many superscripts -- consists of arrow diagrams on uni-coloured skeleta (elements of $\calS^o$), given by the presentation
 \[ \calA^{sw}=\CA\!\left.\left.\left\langle
  \raisebox{-2mm}{\input{figs/wTFprojgens.pstex_t}}
  \right|
  \parbox{1.5in}{\centering $\aft$, TC, VI, CP, RI}
  \right|
  \parbox{0.9in}{\centering $S_e, u_e, d_e, -$}
  \right\rangle.
\]
We denote by $-: \calA^{sw}\to \calA^{sw}$ the associated graded operation of wenjugation. It is also an involution, and coincides with the operation $D\mapsto \overline{D}$ defined in Section~\ref{subsec:wTFExpansion}. Namely, $\overline{D}$ is the arrow diagram $D$ multiplied with $(-1)^{\#\{\text{arrow tails}\}}$.

As before, arrow diagrams have an alternative, equivalent description in terms of Jacobi diagrams, as in Theorem~\ref{thm:FoamBracketRise}.

The main theorem of this section states that homomorphic expansions for $\wTFo$ are in bijection with Kashiwara-Vergne solutions, without restriction on the Duflo function:

\begin{theorem}
\label{thm:ATEquivalence} There exist a group-like homomorphic expansions for $\wTFo$, and
there is a bijection
between the set of solutions $(F,a)$ of the generalized KV equations \eqref{eq:ATKVEqns} and the set of v-small group-like homomorphic expansions
for $\wTFo$.
\end{theorem}

\begin{proof}
Since there are no wens, a homomorphic expansion is determined by the values $R$, $C$, and $V$, with  $Z(\undercrossing)=(R^{-1})^{21}$ by the R2 relation, and $Z(\raisebox{-1mm}{\input{figs/MinusVertex.pstex_t}})=V^{-1}$ by the homomorphicity with respect to edge unzip.

We derive $R=e^a$, and the \eqref{eq:HardR4} and \eqref{eq:CapEqn} equations as before: from the R3, R4 and CP relations, and the homomorphicity with respect to $S_e$, $u_e$ and $d_e$. There are no wen relations, hence no restriction on the odd part of $C$, nor a Unitarity equation.

Recall that for $\wTF$, the TV relation gave rise to the unitarity equation. Since one side of the TV relation is the wenjugate of the vertex $\raisebox{-1mm}{\input{figs/PlusVertex.pstex_t}}$. Thus, homomorphicity with respect to wenjugation is equivalent to the Unitarity equation \eqref{eq:unitarity}.

We showed in the proof or Theorem~\ref{thm:WenATEquivalence} that the equations \eqref{eq:HardR4}, \eqref{eq:CapEqn} and \eqref{eq:unitarity}, given the v-small condition, translate exactly to the Kashiwara--Vergne equations. This completes the proof.
\qed
\end{proof}


\draftcut
\subsection{Interlude: $u$-Knotted Trivalent Graphs}
\label{subsec:KTG}
The ``$u$sual'', or classical knot-theoretical objects corresponding to
$\wTF$ are loosely speaking Knotted Trivalent Graphs, or \glost{KTGs}.
We give a brief introduction/review of this structure before studying the 
relationship between their homomorphic expansions and homomorphic expansions for 
$\wTF$. The last goal of this paper is to show that the topological relationship between the
two spaces explains the relationship between the KV problem and Drinfel'd associators.

A trivalent graph is a graph with three edges meeting at each vertex,
equipped with a cyclic orientation of the three half-edges at each
vertex. KTGs are framed embeddings of trivalent graphs into $\bbR^3$,
regarded up to isotopies. The skeleton of a KTG is the trivalent
graph (as a combinatorial object) behind it.  For a detailed
introduction to KTGs see for example \cite{Bar-NatanDancso:KTG}.
Here we only recall the most important facts. The reader might
recall that in Section~\ref{1-sec:w-knots}, the w-knot section,
of \cite{Bar-NatanDancso:WKO1}  we only dealt with long $w$-knots,
as the $w$-theory of round knots is essentially trivial (see
\cite[Theorem~\ref{1-prop:AwCirc}]{Bar-NatanDancso:WKO1}). A similar issue
arises with ``$w$-knotted trivalent graphs''. Hence, the space we are
really interested in is ``long KTGs'', meaning, trivalent tangles with
1 or 2 ends.

\parpic[r]{\input{figs/UnzipAndInsertion.pstex_t}}
Long KTGs form an algebraic structure with operations as follows. {\em Orientation
switch} reverses the orientation of a specified edge. {\em Edge unzip} doubles a 
specified edge as shown on the right. {\em Tangle
insertion} is inserting a small copy of a $(1,1)$-tangle $S$ into
the middle of some specified edge of a tangle $T$, as shown in the second row on
the right (tangle composition is a special case of this). The {\em stick-on} operation ``sticks a 1-tangle $S$ onto a specified edge of another
tangle $T$'', as shown. (In the figures $T$ is a 2-tangle, but this is irrelevant.) {\em Disjoint union} of
two 1-tangles produces a 2-tangle.
Insertion, disjoint union and stick-on are a slightly weaker set of operations than the connected sum
of~\cite{Bar-NatanDancso:KTG}.  

The associated graded structure of the algebraic structure of long KTGs is the 
graded space $\glos{\calA^u}$ of chord diagrams on
trivalent graph skeleta, modulo the $\glos{4T}$ and vertex invariance
(VI) relations. The induced operations on $\calA^u$ are as expected:
orientation switch multiplies a chord diagram by $(-1)$ to the number
of chord endings on the edge.  The edge unzip $u_e$ maps a chord diagram
with $k$ chord endings on the edge $e$ to a sum of $2^k$ diagrams where
each chord ending has a choice between the two daughter edges. Finally,
tangle insertion, stick-on and disjoint union 
induces the insertion, sticking on and disjoint union 
of chord diagrams, respectively.

\parpic[r]{\input{figs/glitch.pstex_t}}
In \cite{Bar-NatanDancso:KTG} the authors prove that there is no
\emph{homomorphic} expansion for KTGs. This theorem, as well as the proof,
applies to long KTGs with slight modifications. However there are well-known --- and nearly homomorphic ---
expansions constructed by extending the Kontsevich integral to KTGs, 
or from Drinfel'd associators. There are several such constructions 
(\cite{MurakamiOhtsuki:KTGs}, \cite{ChepteaLe:EvenAssociator},
\cite{Dancso:KIforKTG}). For now, let us denote any one of these expansions by
$Z^{old}$. All $Z^{old}$ are ``almost homomorphic'': they intertwine every operation
except for edge unzip with their
chord-diagrammatic counterparts; but commutativity with unzip fails by
a controlled amount, as shown on the right. Here $\glos{\nu}$ denotes
the ``invariant of the unknot'', the value of which was conjectured in
\cite{Bar-NatanGaroufalidisRozanskyThurston:WheelsWheeling} and proven
in \cite{Bar-NatanLeThurston:TwoApplications}.

In \cite{Bar-NatanDancso:KTG} the authors fix this anomaly by slightly
changing the space of KTGs and adding some extra combinatorics (``dots''
on the edges), and construct a homomorphic expansion for this new space by
a slight adjustment of $Z^{old}$. Here we are going to use a similar but
different adjustment of the space of trivalent 1- and 2-tangles. Namely
we break the symmetry of the vertices and restrict the set of allowed
unzips.

\begin{definition}\label{def:sKTG} A ``signed KTG'', denoted $\glos{\sKTG}$, is 
a trivalent oriented 1- or 2-tangle embedded in $\bbR^3$ with a cyclic orientation of edges meeting at each vertex,
and in addition each vertex is equipped with a sign and one of the three incident edges is marked as distinguished (sometimes denoted
by a thicker line). Our pictorial convention will
be that a vertex drawn in a ``\inverted{$Y$}'' shape
with all edges oriented up and the top edge distinguished is always
positive and a vertex drawn in a ``$Y$'' shape with edges oriented
up and the bottom edge distinguished is always negative (see Figure
\ref{fig:ZatVertices}).

\parpic[r]{\input{figs/StickOnSigns.pstex_t}}
The algebraic structure $\sKTG$ has one kind of objects for each skeleton (a skeleton is a uni-trivalent graph with signed vertices but
no embedding), as well as several operations: orientation switch, edge unzip, tangle insertion, disjoint union of 1-tangles, and stick-on.
Orientation switch of either of the non-distinguished edges changes the sign of the vertex, switching the orientation of 
the distinguished edge does not. Unzip of an edge
is only allowed if the edge is distinguished at both of its ends and the vertices at either end are of opposite signs. 
The stick-on operation can be done in either one of the two ways shown on the right (i.e., the stuck-on edge can be attached at 
a vertex of either sign, but it can not become the distinguished edge of that vertex).
\end{definition}

To consider expansions of $\sKTG$, and ultimately the compatibility of
these with $Z^w$, we first note that $\sKTG$ is finitely generated (and
therefore any expansion $Z^u$ is determined by its values on finitely
many generators). The proof of this is not hard but somewhat lengthy,
so we postpone it to Section~\ref{subsec:sKTGgensProof}.

\begin{proposition}\label{prop:sKTGgens}
The algebraic structure $\sKTG$ is finitely generated by the following
list of elements:
\begin{center} \input{figs/sKTGgens.pstex_t} \end{center}
\end{proposition}

Note that we ignore edge orientations for simplicity in the statement of this proposition; this is not a problem as orientation switches are 
allowed in $\sKTG$ without restriction.


\subsubsection{Homomorphic expansions for $\sKTG$}\label{subsubsec:Zu}
Suppose that $Z^u:\sKTG \to \calA^u$ is a homomorphic expansion. We hope to determine the value of $Z^u$ on each of the generators. 

\parpic[r]{\input{figs/BubbleSquared.pstex_t}}
The value of $Z^u$ on the single strand is an element of $\calA^u(\uparrow)$ whose square is itself, hence it is 1. 
The value of 
the bubble is an element $x \in \calA^u(\uparrow_2)$, as all chords can be pushed to the ``bubble'' part using the VI relation. Two bubbles
can be composed and unzipped to produce a single bubble (see on the right), hence we have $x^2=x$, which implies $x=1$ 
in $\calA^u(\uparrow_2)$.

Recall that a Drinfel'd associator is a group-like
element $\Phi \in \calA^u(\uparrow_3)$ along with a group-like element $R^u \in \calA^u(\up_2)$ satisfying the so-called pentagon and
positive and negative hexagon equations, as well as a non-degeneracy
and mirror skew-symmetry property.  For a detailed explanation see
Section 4 of \cite{Bar-NatanDancso:KTG}; associators were first defined
in \cite{Drinfeld:QuasiHopf}. We claim that the $Z^u$-value $\glos{\Phi}$ of the right
associator, along with the value $\glos{R^u}$ of the right twist forms a Drinfel'd associator pair. The proof of this statement
is the same as the proof of Theorem 4.2 of \cite{Bar-NatanDancso:KTG},
with minor modifications (making heavy use of the assumption 
that $Z^u$ is homomorphic).
It is easy to check by composition and unzips that the value of the 
left associator and the left twist are $\Phi^{-1}$ and $(R^u)^{-1}$, respectively.
Note that if $\Phi$ is a {\em horizontal chord} associator (i.e., 
all the chords of $\Phi$ are horizontal on three strands) then $R^u$ is forced to
be $e^{c/2}$ where $c$ denotes a single chord. Note that the reverse is not true:
there exist non-horizontal chord associators $\Phi$ that satisfy the hexagon equations with $R^u=e^{c/2}$.

\parpic[r]{\input{figs/NooseBalloonPhi.pstex_t}}
Let $b$ and $n$ denote the $Z^u$-values of the balloon and the noose, respectively. Note that using the $VI$ relation all chord endings
can be pushed to the ``looped'' strands, so $b$ and $n$ live in $\calA^u(\uparrow)$, as seen in Figure \ref{fig:NooseBalloonProof}. The argument in that figure
shows that $n\cdot b$ is the inverse in $\calA^u(\uparrow)$ of ``an associator on a squiggly strand'', as shown on the right.
In Figure \ref{fig:NooseBalloonProof} we start with the $\sKTG$ on the top left and either apply $Z^u$ followed by unzipping the 
edges marked 
by stars, or first unzip the same edges and then apply $Z^u$. 
Since $Z^u$ is homomorphic, the two results in the bottom right corner must agree. 
(Note that two of the four unzips we perform are ``illegal'',
as the strand directions can't match. However, it is easy to get around this issue by inserting small bubbles at the top of the balloon and the bottom 
of the noose, and switching the appropriate edge orientations before and after the unzips. The $Z^u$-value of a bubble is 1, hence this will not effect 
the computation and so we ignore the issue for simplicity.)

\begin{figure}[h]
  \input{figs/NooseBalloonProof.pstex_t}
\caption{Unzipping a noose and a balloon to a squiggle.}
\label{fig:NooseBalloonProof}
\end{figure}

In addition, it follows from Theorem 4.2 of \cite{Bar-NatanDancso:KTG} via deleting two edges 
that the inverse of an ``associator on a squiggly strand'' is $\nu$, the invariant of the unknot.
To summarize, we have proven the following:

\begin{lemma}\label{lem:nb} If $Z^u$ is a homomorphic expansion then the $Z^u$ values of
the strand and the bubble are 1, the values of the right associator and right twist form an associator pair $(\Phi,R^u)$,
and the values of the left twist and left associator are inverses of these.
With $n$ and $b$ denoting the value of the noose and the balloon, respectively, and $\nu$
being the invariant of the unknot, we have $n \cdot b =\nu$ in $\calA^u(\uparrow)$.
\end{lemma}

The natural question to ask is whether any triple $(\Phi, R^u, n)$ gives rise to a homomorphic expansion. We don't know
whether this is true, but we do know that any pair $(\Phi, R^u)$ gives rise to a ``nearly homomorphic'' expansion of KTGs 
\cite{MurakamiOhtsuki:KTGs, ChepteaLe:EvenAssociator, Dancso:KIforKTG},
and we can construct a homomorphic expansion for $\sKTG$ from any of these (as shown below). However, all of these expansions
take the same specific value on the noose and the balloon (also see below). We don't know whether there really is a one 
parameter family of homomorphic expansions $Z^u$ for each choice of $(\Phi, R^u)$ or if we are simply missing 
a relation. 

\parpic[r]{\input{figs/ZoldOfTangle.pstex_t}}
We now construct explicit homomorphic expansions $Z^u \colon  \sKTG \to \calA^u$ from any $Z^{old}$ (where $Z^{old}$ stands for
an ``almost homomorphic'' expansion of KTGs) as follows. First of all we need to interpret
$Z^{old}$ as an invariant of 2-tangles. This can be done by connecting the top and bottom ends by a non-interacting long
strand followed by a normalization, as shown on the right. By ``multiplying by $\nu^{-1}$'' we mean that after computing $Z^{old}$
we insert $\nu^{-1}$ on the long strand (recall that $\nu$ is the ``invariant of the unknot''). We interpret $Z^{old}$ of a 1-tangle
as follows: stick the 1-tangle onto a single strand to obtain a 2-tangle, then proceed as above. The result will only have chords on the 
1-tangle (using that the extensions of the Kontsevich Integral are homomorphic with respect to ``connected sums''), 
so we define the result to be the value of $Z^{old}$ on the 1-tangle.
As an example, we compute the value of $Z^{old}$ for the noose in Figure \ref{fig:uValueNoose} 
(note that the computation for the balloon is the same).

\begin{figure}
 \input{figs/uValueOfTheNoose.pstex_t}
 \caption{Computing the $Z^{old}$ value of the noose. The third step uses that the Kontsevich integral of KTGs is homomorphic with
 respect to the ``connected sum'' operation and that the value of the unknot is $\nu$ (see \cite{Bar-NatanDancso:KTG} for an explanation of
 both of these facts).}
 \label{fig:uValueNoose}
\end{figure}

\begin{figure}[h]
\input{figs/ZatVertices.pstex_t}
\caption{Normalizations for $Z^u$ at the vertices.}\label{fig:ZatVertices}
\end{figure}

\parpic[r]{\input{figs/NooseBalloonValue.pstex_t}}
Now to construct a homomorphic $Z^u$ from $Z^{old}$ we add normalizations
near the vertices, 
as in Figure~\ref{fig:ZatVertices}, where $c$ denotes a single chord. 
Checking that $Z^u$ is a homomorphic expansion is a simple calculation
using the almost homomorphicity of $Z^{old}$, which we leave to the
reader. The reader can also verify that $Z^u$ of the strand and the bubble 
is 1 as it should be. $Z^u$ of the right twist is $e^{c/2}$ and $Z^u$ of the right associator is a Drinfel'd associator
$\Phi$ (note that $\Phi$ depends on which $Z^{old}$ was used).
From the calculation of Figure \ref{fig:uValueNoose} it follows that the $Z^u$
value of the balloon and the noose (for any $Z^{old}$) are as shown on the right,
and indeed $n\cdot b = \nu$.


\subsection{The relationship between $\sKTG$ and $\wTFo$}\label{subsec:wTFcompatibility}
We move on to the question of compatibility between the homomorphic expansions
$Z^u$ and $Z^w$ (from now on we are going to refer to the homomorphic
expansion of $\wTFo$ --- called $Z$ in the previous section --- as $Z^w$
to avoid confusion).

There is a map $a\colon  \sKTG \to \wTFo$, given by interpreting $\sKTG$
diagrams as $\wTFo$ diagrams. In particular, positive vertices (of edge
orientations as shown in Figure \ref{fig:ZatVertices}) are interpreted as the $\wTFo$ vertex
\input{figs/PlusVertex.pstex_t} and negative vertices as the $\wTFo$ vertex
\input{figs/NegVertex.pstex_t}. (The map $a$ can also be interpreted topologically as 
Satoh's tubing map.) The induced map $\alpha\colon  \calA^u
\to \calA^{sw}$ is as defined in Section \ref{subsec:sder}, that is,
$\alpha$ maps each chord to the sum of its two possible orientations.
Hence we can ask whether the two expansions are compatible (or can 
be chosen to be compatible), which takes us to the main result of 
this section:

\def\uwsquare{{\xymatrix{
  \sKTG \ar[d]^{Z^u} \ar[r]^a & \wTFo \ar[d]^{Z^w} \\
  \calA^u \ar[r]^\alpha & \calA^{sw}
}}}

\parpic(2in,1in)[r]{\null\raisebox{-5mm}{\begin{minipage}{2in}
  \begin{equation} \label{eq:uwcompatibility}
    \begin{array}{c} \uwsquare \end{array}
  \end{equation}
\end{minipage}}}
\begin{theorem}\label{thm:ZuwCompatible}
\picskip{4}
Let $Z^u$ be a homomorphic expansion for $\sKTG$ with the properties that $\Phi$
is a horizontal chord associator and $n=e^{-c/4}\nu^{1/2}$ in the sense of Section \ref{subsubsec:Zu}.\footnote{It will become apparent that in the
proof we only use slightly weaker but less aesthetic conditions on $Z^u$.}
Then there exists a homomorphic expansion $Z^w$ for $\wTFo$ compatible with $Z^u$ in the sense
that the square on the right commutes. 

Furthermore, such $Z^w$ are in one to one correspondence\footnote{An even nicer theorem would be a classification of 
homomorphic expansions for the combined algebraic structure $\left(\sKTG\overset{a}{\longrightarrow}\wTFo\right)$ in terms of solutions of 
the KV problem. The two obstacles to this are clarifying whether there is a free choice of $n$ for $Z^u$, and --- probably much harder --- how much of the 
horizontal chord condition is necessary for a compatible $Z^w$ to exist.} with
``symmetric solutions of the KV problem'' satisfying the KV equations \eqref{eq:ATKVEqns}, the ``twist equation''
\eqref{eq:TwistWithF} and the associator equation \eqref{eq:ATPhiandV}.
\end{theorem}

\picskip{0}Before moving on to the proof let us state and prove the following Lemma,
to be used repeatedly in the proof of the theorem.

\begin{lemma}\label{lem:TreesAndUnitarity}
If $a$ and $b$ are group-like elements in $\calA^{sw}(\uparrow_n)$, then $a=b$ if and only if $\pi(a)=\pi(b)$ and $aa^*=bb^*$. Here $\pi$
is the projection induced by $\pi\colon  \calP^w(\uparrow_n) \to \tder_n \oplus \fraka_n$ (see Section \ref{subsec:ATSpaces}),
and $*$ refers to the adjoint map of Definition \ref{def:Adjoint}. \end{lemma}

\begin{proof}
Write $a=e^we^{uD}$ and $b=e^{w'}e^{uD'}$, where $w\in \attr_n$, $D\in \tder_n\oplus \fraka_n$ 
and $u\colon  \tder_n\oplus \fraka_n \to \calP_n$ is the ``upper'' map of 
Section \ref{subsec:ATSpaces}. Assume that $\pi(a)=\pi(b)$ and $aa^*=bb^*$. Since $\pi(a)=e^D$ and $\pi(b)=e^{D'}$, we conclude
that $D=D'$. Now we compute $aa^*=e^we^{uD}e^{-lD}e^w=e^we^{j(D)}e^w,$ where $j\colon  \tder_n \to \attr_n$ is the map defined in Section 5.1 of 
\cite{AlekseevTorossian:KashiwaraVergne} and discussed in \ref{prop:Jandj} of this paper. Now note that both $w$ and $j(D)$ are elements of 
$\attr_n$, hence they commute, so $aa^*=e^{2w+j(D)}$. Thus, $aa^*=bb^*$ means that $e^{2w+j(D)}=e^{2w'+j(D)}$, which implies that $w=w'$ and
$a=b$. \qed
\end{proof}

\medskip

\noindent{\em Proof of Theorem \ref{thm:ZuwCompatible}.}
In addition to being a homomorphic expansion for $\wTFo$, $Z^w$ has to satisfy an
the added condition of being compatible with $Z^u$. Since $\sKTG$ is finitely generated, this translates
to one additional equation for each generator of $\sKTG$, some of which are
automatically satisfied. To deal with the others, we use the machinery established in the previous sections
to translate these equations to conditions on $F$, and they turn out to be the properties studied in \cite{AlekseevTorossian:KashiwaraVergne}
which link solutions of the KV problem with Drinfel'd associators. 

To start, note that for the single strand and the bubble the commutativity of the square \eqref{eq:uwcompatibility} is
satisfied with any $Z^w$: both the $Z^u$ and $Z^w$ values are 1 (note that the $Z^w$ value of the bubble 
is 1 due to the unitarity \eqref{eq:unitarity} of $Z^w$). Each of the other generators will require more study.

{\em Commutativity of~\eqref{eq:uwcompatibility}
for the twists.} Recall that the $Z^u$-value of the right twist (for a $Z^u$ with horizontal chord $\Phi$)
is $R^u=e^{c/2}$; and note that its $Z^w$-value is $V^{-1}RV^{21}$,
where $R=e^{a_{12}}$ is the $Z^w$-value of the crossing (and $a_{12}$ is a single arrow pointing from strand 1 to strand 2). Hence the
commutativity of \eqref{eq:uwcompatibility} for the right twist is equivalent to the
``Twist Equation'' $\alpha(R^u)=V^{-1}RV^{21}$.
By definition
of $\alpha$, $\alpha(R^u)=e^{\frac{1}{2}(a_{12}+a_{21})}$, where $a_{12}$ and $a_{21}$ are single arrows pointing from strand 1 to 2 and 2 to 1, 
respectively. Hence we have
\begin{equation}\label{eq:twist}
 e^{\frac{1}{2}(a_{12}+a_{21})}=V^{-1}RV^{21}.
\end{equation}
To translate this to the language of \cite{AlekseevTorossian:KashiwaraVergne}, we use
Lemma \ref{lem:TreesAndUnitarity}, which implies that it is enough
for $V$ to satisfy the Twist Equation ``on tree level'' (i.e., after applying $\pi$), and for which the adjoint condition of the Lemma holds.

We first prove that the adjoint condition holds for any homomorphic expansion of $\wTFo$. Multiplying the left hand side of the Twist Equation by its adjoint, we get 
$$e^{\frac{1}{2}(a_{12}+a_{21})}(e^{\frac{1}{2}(a_{12}+a_{21})})^*=e^{\frac{1}{2}(a_{12}+a_{21})}e^{-\frac{1}{2}(a_{12}+a_{21})}=1.$$
As for the right hand side, we have to compute $V^{-1}RV^{21}(V^{21})^*R^*(V^{-1})^*$. Since $V$ is unitary (Equation (\ref{eq:unitarity})), $VV^*=1$.
Now $R=e^{a_{12}}$, so $R^*=e^{-a_{12}}=R^{-1}$, hence the expression on the right hand side also simplifies to 1, as needed.  

As for the ``tree level'' of the Twist Equation, recall that in Section \ref{subsec:wTFExpansion} we used 
Alekseev and Torossian's solution $F\in \TAut_2$ to the Kashiwara--Vergne
equations \cite{AlekseevTorossian:KashiwaraVergne} to find solutions $V$ to equations (\ref{eq:HardR4}),(\ref{eq:unitarity}) and (\ref{eq:CapEqn}). 
We produced $V$ from $F$ by setting $F=e^{D^{21}}$ with $D\in \tder_2^s$, $b:=\frac{-j(F)}{2}\in \attr_2$ and $V:=e^be^{uD}$, so $F$ is ``the tree part'' of $V$, 
up to re-numbering strands. Hence, the tree level Twist Equation translates to a new equation for $F$.
Substituting $V=e^be^{uD}$ into the Twist Equation we obtain
$ e^{\frac{1}{2}(a_{12}+a_{21})}=e^{-uD}e^{-b}e^{a_{12}}e^{b^{21}}e^{uD^{21}},$
and applying $\pi$, we get
\begin{equation}\label{eq:TwistWithF}
 e^{\frac{1}{2}(a_{12}+a_{21})}=(F^{21})^{-1}e^{a_{12}}F.
\end{equation}
In \cite{AlekseevTorossian:KashiwaraVergne} the solutions $F$ of the KV equations which also satisfy
this equation are called ``symmetric solutions of
the Kashiwara-Vergne problem'' discussed in Sections 8.2
and 8.3. (Note that in
\cite{AlekseevTorossian:KashiwaraVergne} $R$ denotes $e^{a_{21}}$). 

{\em Commutativity of~\eqref{eq:uwcompatibility} for the
associators.} 
Recall that the $Z^u$ value of the right associator is a Drinfel'd associator $\Phi \in \calA^u(\uparrow_3)$;
for the $Z^w$ value see Figure \ref{fig:wAssociator}.
Hence the new condition on $V$ is the following:
\begin{equation}\label{eq:AssociatorAndV}
  \alpha(\Phi)=V_-^{(12)3}V_-^{12}V^{23}V^{1(23)}
  \qquad\text{in}\qquad
  \calA^{sw}(\uparrow_3)
\end{equation}

\begin{figure}
 \input{figs/Associator.pstex_t}
 \caption{The $Z^w$-value of the right associator.}
 \label{fig:wAssociator}
\end{figure}

Again we treat the ``tree and wheel parts'' separately
using Lemma \ref{lem:TreesAndUnitarity}. As $\Phi$ is by
definition group-like, let us denote $\Phi=:\glos{e^\phi}$. 
We first verify that the ``wheel part'' or adjoint condition of the Lemma holds. Starting
with the right hand side of Equation~\eqref{eq:AssociatorAndV}, the
unitarity $VV^*=1$
of $V$ implies that
\[ V_-^{(12)3}V_-^{12}V^{23}V^{1(23)}
  (V^{1(23)})^*(V^{23})^*(V_-^{12})^*(V_-^{(12)3})^*=1.
\]
For the left hand side of~\eqref{eq:AssociatorAndV} we need to show that
$e^{\alpha(\phi)}(e^{\alpha(\phi)})^*=1$ as well, and this is true for 
any {\em horizontal chord} associator. 
Indeed, restricted to the $\alpha$-images of horizontal
chords $*$ is multiplication by $-1$, and as it is an anti-Lie morphism,
this fact extends to the Lie algebra generated by $\alpha$-images
of horizontal chords. Hence $e^{\alpha(\phi)}(e^{\alpha(\phi)})^*
=e^{\alpha(\phi)}e^{\alpha(\phi)^*}=e^{\alpha(\phi)}e^{-\alpha(\phi)}=1$.

On to the tree part, applying $\pi$ to Equation (\ref{eq:AssociatorAndV})
and keeping in mind that $V_-=V^{-1}$ by the unitarity of $V$, we obtain
\begin{multline}\label{eq:ATPhiandV}
  e^{\pi\alpha(\phi)}
  = (F^{3(12)})^{-1} (F^{21})^{-1} F^{32} F^{(23)1}
  =e^{-D^{(12)3}}e^{-D^{12}}e^{D^{23}}e^{D^{1(23)}} \\
  \text{in }
  \glos{\SAut_3}:=\exp(\sder_3)\subset\TAut_3.
\end{multline}
This is Equation (26) of \cite{AlekseevTorossian:KashiwaraVergne},
up to re-numbering strands 1 and 2 as 2 and 1\footnote{Note that
in \cite{AlekseevTorossian:KashiwaraVergne} ``$\Phi'$ is an
associator'' means that $\Phi'$ satisfies the pentagon equation,
mirror skew-symmetry, and positive and negative hexagon equations
in the space $\SAut_3$. These equations are stated in
\cite{AlekseevTorossian:KashiwaraVergne} as equations (25), (29),
(30), and (31), and the hexagon equations are stated with strands 1
and 2 re-named to 2 and 1 as compared to \cite{Drinfeld:QuasiHopf}
and \cite{Bar-NatanDancso:KTG}. This is consistent with $F=e^{D^{21}}$.}.
The following fact from
\cite{AlekseevTorossian:KashiwaraVergne} (their Theorem 7.5, Propositions
9.2 and 9.3 combined) implies that there is a solution $F$ to the KV equations \eqref{eq:ATKVEqns} 
which also satisfies \eqref{eq:TwistWithF} and \eqref{eq:ATPhiandV}.

\begin{fact}
If $\Phi'=e^{\phi'}$ is an associator in $\SAut_3$ so that $j(\Phi')=0$\footnote{The condition
$j(\phi')=0$ is equivalent to the condition $\Phi\in KRV^0_3$ in \cite{AlekseevTorossian:KashiwaraVergne}.
The relevant definitions in \cite{AlekseevTorossian:KashiwaraVergne} can be found in Remark 4.2 and at the bottom of 
page 434 (before Section 5.2).}
then Equation~(\ref{eq:ATPhiandV}) has a solution $F=e^{D^{21}}$ which is
also a solution to the KV equations, and all such solutions are symmetric
(i.e. verify the Twist Equation (\ref{eq:TwistWithF})). \qed
\end{fact}

To use this Fact, we need to show that $\Phi':=\pi\alpha(\Phi)$ is an
associator in $\SAut_3$ and that $j(\Phi')=j(\pi\alpha(\Phi))=0$. The
latter is the unitarity of $\Phi$ which is already proven. The
former follows from the fact that $\Phi$ is an associator and the fact
(Theorem~\ref{thm:sder}) that the image of $\pi\alpha$ is contained in
$\sder$ (ignoring degree 1 terms, which are not present in an associator
anyway).

In summary, the condition of the Fact are satisfied and so there exists a solution $F$ which 
in turn induces a $Z^w$ which is compatible with $Z^u$ for the strand, the bubble, the twists and 
the associators. That is, all
generators of $\sKTG$ except possibly the balloon and the noose. As the last step of the proof 
of Theorem \ref{thm:ZuwCompatible} we show that any such $Z^w$ also automatically make
\eqref{eq:uwcompatibility} commutative for the balloon and the noose.

\parpic[r]{\input{figs/wValueOfTheNoose.pstex_t}}
{\em Commutativity of~\eqref{eq:uwcompatibility}
for the balloon and the noose.} Since we know the $Z^u$-values $B$ and $n$ of the balloon and the noose,
we start by computing $Z^w$ of the noose. $Z^w$ assigns a $V$ value to the vertex with the
first strand orientation switched as shown in the figure on the
right. The balloon is the same, except with the inverse vertex and the second strand reversed.
Hence what we need to show is that the two equations 
below hold:
\begin{center}
 \input{figs/NooseEquations.pstex_t}
\end{center}

Let us denote the left hand side of the first and second equation above by $n^w$ and $b^w$, respectively (that is, the $Z^w$ value of the noose
and the balloon, respectively).
We start by proving that the product of these two equations holds,
namely that $n^wb^w=\alpha(\nu)$.
(We used that any local (small) arrow diagram on a single strand is central in
$\calA^{sw}(\uparrow_n)$, hence the cancellations.)
This product equation is satisfied due to an argument identical to that of
Figure~\ref{fig:NooseBalloonProof}, but carried out in $\wTFo$, and using 
that by the compatibility with associators, $Z^w$ of an associator is $\alpha(\Phi)$. 

What remains is to show that the noose and balloon equations hold individually. In light of the results so far, it is sufficient to show that
\begin{equation}\label{eq:NooseSymmetry}
n^w=b^w\cdot e^{-D_A}, 
\end{equation}
where $D_A$ stands for a single arrow on one strand (whose direction
doesn't matter due to the $RI$ relation.  As stated in
\cite[Theorem~\ref{1-thm:Aw}]{Bar-NatanDancso:WKO1},
$\calA^{sw}(\uparrow_1)$ is the polynomial algebra freely generated by the
arrow $D_A$ and wheels of degrees 2 and higher. Since $V$ is group-like,
$n^w$ (resp. $b^w$) is an exponential $e^{A_1}$ (resp. $e^{A_2}$)
with $A_1, A_2 \in \calA^{sw}(\uparrow_1)$. We want to show that
$e^{A_1}=e^{A_2}\cdot e^{-D_A}$, equivalently that $A_1=A_2-D_A$.

\begin{figure}
  \input figs/NooseCappedProof.pstex_t
\caption{The proof of Equation (\ref{eq:NooseCapped}). Note that the unzips are ``illegal'', as the strand directions don't match. This can be fixed
by inserting a small bubble at the bottom of the noose and doing a number of orientation switches. As this doesn't change the result or the main argument,
we suppress the issue for simplicity. Equation (\ref{eq:NooseCapped}) is obtained from this result by multiplying by $S(C)^{-1}$ on the bottom and by $C^{-1}$
on the top.}\label{fig:NooseCappedProof}
\end{figure}

In degree 1, this can be done by explicit verification. Let $A_1^{\geq
2}$ and $A_2^{\geq 2}$ denote the degree 2 and higher parts of $A_1$ and
$A_2$, respectively. We claim that capping the strand at both its top 
and its bottom takes $e^{A_1}$ to $e^{A_1^{\geq 2}}$, and similarly $e^{A_2}$
to $e^{A_2^{\geq 2}}$. (In other words, capping kills arrows but leaves
wheels un-changed.) This can be proven similarly to the proof of
Lemma~\ref{lem:CapIsWheels}, but using
\[
  F' := \sum_{k_1,k_2=0}^{\infty}
    \frac{(-1)^{k_1+k_2}}{k_1!k_2!}D_A^{k_1+k_2}S_L^{k_1}S_R^{k_2}
\]
in place of $F$ in the proof. What we need to prove, then, is the following equality, and
the proof is shown in Figure~\ref{fig:NooseCappedProof}.
\begin{equation}\label{eq:NooseCapped}
  \raisebox{-6mm}{\input figs/NooseCapped.pstex_t}
\end{equation}
This concludes the proof of Theorem~\ref{thm:ZuwCompatible}. \qed

\vspace{2mm}
Recall from Section~\ref{subsec:sder} that there is no commutative
square linking $Z^u\colon\uT\to\calA^u$ and $Z^w\colon\wT\to\calA^{sw}$,
for the simple reason that the Kontsevich integral for tangles $Z^u$
is not canonical, but depends on a choice of parenthesizations for
the ``bottom'' and the ``top'' strands of a tangle $T$. Yet given
such choices, a tangle $T$ can be ``closed up with trees'' as within the proof of
Proposition~\ref{prop:sKTGgens} (see Section \ref{sec:odds}) into an $\sKTG$ which we will denote
$G$. For $G$ a commutativity statement does hold as we have just
proven. The $Z^u$ and $Z^w$ invariants of $T$ and of $G$ differ only
by a number of vertex-normalizations and vertex-values on skeleton-trees
at the bottom or at the top of $G$, and using VI, these values can slide
so they are placed on the original skeleton of $T$. This is summarized
as the following proposition:

\begin{proposition} \label{prop:uwBT} Let $n$ and $n'$ be natural numbers.
Given choices $c$ and and $c'$ of parenthesizations of $n$ and $n'$
strands respectively, there exists invertible elements
$C\in\calA^{sw}(\uparrow_n)$ and $C'\in\calA^{sw}(\uparrow_{n'})$ so
that for any u-tangle $T$ with $n$ ``bottom'' ends and  $n'$ ``top'' ends
we have
\[ \alpha Z^u_{c,c'}(T)=C^{-1}Z^w(aT)C', \]
where $Z^u_{c,c'}$ denotes the usual Kontsevich integral of $T$ with
bottom and top parenthesizations $c$ and $c'$.
\end{proposition}

For u-braids the above proposition may be stated with $c=c'$ and then $C$
and $C'$ are the same.

%% file: figs/VertexExamples.pstex_t
\begin{picture}(0,0)%
\includegraphics{figs/VertexExamples.pstex}%
\end{picture}%
%
%
\setlength{\unitlength}{3158sp}%
\begingroup\makeatletter\ifx\SetFigFont\undefined%
\gdef\SetFigFont#1#2#3#4#5{%
  \reset@font\fontsize{#1}{#2pt}%
  \fontfamily{#3}\fontseries{#4}\fontshape{#5}%
  \selectfont}%
\fi\endgroup%
\begin{picture}(1228,325)(3028,-1274)
\end{picture}%

%% file: figs/Cap.pstex_t
\begin{picture}(0,0)%
\includegraphics{figs/Cap.pstex}%
\end{picture}%
%
%
\setlength{\unitlength}{3158sp}%
\begingroup\makeatletter\ifx\SetFigFont\undefined%
\gdef\SetFigFont#1#2#3#4#5{%
  \reset@font\fontsize{#1}{#2pt}%
  \fontfamily{#3}\fontseries{#4}\fontshape{#5}%
  \selectfont}%
\fi\endgroup%
\begin{picture}(1187,771)(176,-1123)
\end{picture}%

%% file: figs/Wen2.pstex_t
\begin{picture}(0,0)%
\includegraphics{figs/Wen2.pstex}%
\end{picture}%
%
%
\setlength{\unitlength}{4934sp}%
\begingroup\makeatletter\ifx\SetFigFont\undefined%
\gdef\SetFigFont#1#2#3#4#5{%
  \reset@font\fontsize{#1}{#2pt}%
  \fontfamily{#3}\fontseries{#4}\fontshape{#5}%
  \selectfont}%
\fi\endgroup%
\begin{picture}(4737,1339)(-2861,-1744)
\put(-1488,-1180){\makebox(0,0)[lb]{\smash{{\SetFigFont{14}{16.8}{\rmdefault}{\mddefault}{\updefault}{\color[rgb]{0,0,0}$w$}%
}}}}
\put(301,-1186){\makebox(0,0)[lb]{\smash{{\SetFigFont{14}{16.8}{\rmdefault}{\mddefault}{\updefault}{\color[rgb]{0,0,0}$=$}%
}}}}
\end{picture}%

%% file: figs/R4.pstex_t
\begin{picture}(0,0)%
\includegraphics{figs/R4.pstex}%
\end{picture}%
%
%
\setlength{\unitlength}{3158sp}%
\begingroup\makeatletter\ifx\SetFigFont\undefined%
\gdef\SetFigFont#1#2#3#4#5{%
  \reset@font\fontsize{#1}{#2pt}%
  \fontfamily{#3}\fontseries{#4}\fontshape{#5}%
  \selectfont}%
\fi\endgroup%
\begin{picture}(7077,2124)(7561,-1273)
\put(7576,-286){\makebox(0,0)[lb]{\smash{{\SetFigFont{12}{14.4}{\rmdefault}{\mddefault}{\updefault}{\color[rgb]{0,0,0}$R4:$}%
}}}}
\end{picture}%

%% file: figs/CapWen.pstex_t
\begin{picture}(0,0)%
\includegraphics{figs/CapWen.pstex}%
\end{picture}%
%
%
\setlength{\unitlength}{3158sp}%
\begingroup\makeatletter\ifx\SetFigFont\undefined%
\gdef\SetFigFont#1#2#3#4#5{%
  \reset@font\fontsize{#1}{#2pt}%
  \fontfamily{#3}\fontseries{#4}\fontshape{#5}%
  \selectfont}%
\fi\endgroup%
\begin{picture}(2666,771)(347,-1123)
\put(526,-811){\makebox(0,0)[lb]{\smash{{\SetFigFont{10}{12.0}{\rmdefault}{\mddefault}{\updefault}{\color[rgb]{0,0,0}$w$}%
}}}}
\end{picture}%

%% file: figs/DetailedVertex_2.pstex_t
\begin{picture}(0,0)%
\includegraphics{figs/DetailedVertex_2.pstex}%
\end{picture}%
\setlength{\unitlength}{3947sp}%
\begingroup\makeatletter\ifx\SetFigFont\undefined%
\gdef\SetFigFont#1#2#3#4#5{%
  \reset@font\fontsize{#1}{#2pt}%
  \fontfamily{#3}\fontseries{#4}\fontshape{#5}%
  \selectfont}%
\fi\endgroup%
\begin{picture}(1265,1555)(5561,-3576)
\end{picture}%

%% file: figs/CapHeads.pstex_t
\begin{picture}(0,0)%
\includegraphics{figs/CapHeads.pstex}%
\end{picture}%
%
%
\setlength{\unitlength}{3158sp}%
\begingroup\makeatletter\ifx\SetFigFont\undefined%
\gdef\SetFigFont#1#2#3#4#5{%
  \reset@font\fontsize{#1}{#2pt}%
  \fontfamily{#3}\fontseries{#4}\fontshape{#5}%
  \selectfont}%
\fi\endgroup%
\begin{picture}(925,632)(288,-983)
\put(826,-736){\makebox(0,0)[lb]{\smash{{\SetFigFont{10}{12.0}{\rmdefault}{\mddefault}{\updefault}{\color[rgb]{0,0,0}$=0$}%
}}}}
\end{picture}%

%% file: figs/SmallCap.pstex_t
\begin{picture}(0,0)%
\includegraphics{figs/SmallCap.pstex}%
\end{picture}%
%
%
\setlength{\unitlength}{3158sp}%
\begingroup\makeatletter\ifx\SetFigFont\undefined%
\gdef\SetFigFont#1#2#3#4#5{%
  \reset@font\fontsize{#1}{#2pt}%
  \fontfamily{#3}\fontseries{#4}\fontshape{#5}%
  \selectfont}%
\fi\endgroup%
\begin{picture}(64,201)(151,-1123)
\end{picture}%

%% file: figs/Sigma.pstex_t
\begin{picture}(0,0)%
\includegraphics{figs/Sigma.pstex}%
\end{picture}%
%
%
\setlength{\unitlength}{3158sp}%
\begingroup\makeatletter\ifx\SetFigFont\undefined%
\gdef\SetFigFont#1#2#3#4#5{%
  \reset@font\fontsize{#1}{#2pt}%
  \fontfamily{#3}\fontseries{#4}\fontshape{#5}%
  \selectfont}%
\fi\endgroup%
\begin{picture}(3207,1265)(10299,-683)
\put(11688,-30){\makebox(0,0)[lb]{\smash{{\SetFigFont{10}{12.0}{\rmdefault}{\mddefault}{\updefault}{\color[rgb]{0,0,0}$\sigma$}%
}}}}
\end{picture}%

%% file: figs/SkelGen.pstex_t
\begin{picture}(0,0)%
\includegraphics{figs/SkelGen.pstex}%
\end{picture}%
%
%
\setlength{\unitlength}{3158sp}%
\begingroup\makeatletter\ifx\SetFigFont\undefined%
\gdef\SetFigFont#1#2#3#4#5{%
  \reset@font\fontsize{#1}{#2pt}%
  \fontfamily{#3}\fontseries{#4}\fontshape{#5}%
  \selectfont}%
\fi\endgroup%
\begin{picture}(1193,366)(41,431)
\put(156,468){\makebox(0,0)[lb]{\smash{{\SetFigFont{14}{16.8}{\rmdefault}{\mddefault}{\updefault}{\color[rgb]{0,0,0},}%
}}}}
\put(694,468){\makebox(0,0)[lb]{\smash{{\SetFigFont{14}{16.8}{\rmdefault}{\mddefault}{\updefault}{\color[rgb]{0,0,0},}%
}}}}
\end{picture}%

%% file: figs/BubbleSquared.pstex_t
\begin{picture}(0,0)%
\includegraphics{figs/BubbleSquared.pstex}%
\end{picture}%
\setlength{\unitlength}{4144sp}%
\begingroup\makeatletter\ifx\SetFigFont\undefined%
\gdef\SetFigFont#1#2#3#4#5{%
  \reset@font\fontsize{#1}{#2pt}%
  \fontfamily{#3}\fontseries{#4}\fontshape{#5}%
  \selectfont}%
\fi\endgroup%
\begin{picture}(1096,924)(4808,-3268)
\put(5311,-2761){\makebox(0,0)[b]{\smash{{\SetFigFont{12}{14.4}{\rmdefault}{\mddefault}{\updefault}{\color[rgb]{0,0,0}$u$}%
}}}}
\end{picture}%

%% file: figs/NegVertex.pstex_t
\begin{picture}(0,0)%
\includegraphics{figs/NegVertex.pstex}%
\end{picture}%
%
%
\setlength{\unitlength}{3158sp}%
\begingroup\makeatletter\ifx\SetFigFont\undefined%
\gdef\SetFigFont#1#2#3#4#5{%
  \reset@font\fontsize{#1}{#2pt}%
  \fontfamily{#3}\fontseries{#4}\fontshape{#5}%
  \selectfont}%
\fi\endgroup%
\begin{picture}(249,249)(3514,-1498)
\end{picture}%

%% file: odds.tex
\draftcut
\section{Odds and Ends}\label{sec:odds}

\subsection{Motivation for circuit algebras: electronic circuits}
\label{subsec:CAMotivation}
\begin{figure}[h!]
\parpic[r]{\raisebox{-28mm}{\includegraphics[width=50mm]{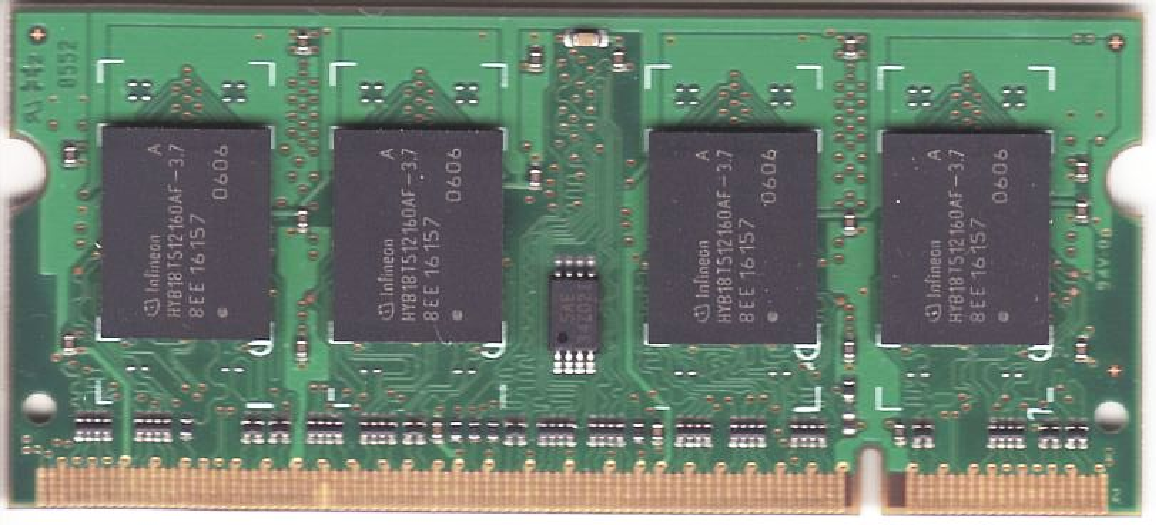}}}
\caption{
  The circuit algebra product of 4 big black components and 1 small black
  component carried out using a green wiring diagram, is an even bigger
  component that has many golden connections (at bottom). When plugged into
  a yet bigger circuit, the CPU board of a laptop, our circuit functions
  as 4,294,967,296 binary memory cells.
} \label{fig:Circuit}
\end{figure}

Electronic circuits are made of ``components'' that can
be wired together in many ways. On a logical level, we only care to
know which pin of which component is connected with which other pin of
the same or other component. On a logical level, we don't really need
to know how the wires between those pins are embedded in space (see
Figures~\ref{fig:FlipFlop} and~\ref{fig:Circuit}). ``Printed Circuit
Boards'' (PCBs) are operators that make smaller components (``chips'')
into bigger ones (``circuits'') --- logically speaking, a PCB is simply a
set of ``wiring instructions'', telling us which pins on which components
are made to connect (and again, we never care precisely how the wires
are routed provided they reach their intended destinations, and ever
since the invention of multi-layered PCBs, all conceivable topologies for
wiring are actually realizable). PCBs can be composed (think ``plugging
a graphics card onto a motherboard''); the result of a composition of
PCBs, logically speaking, is simply a larger PCB which takes a larger
number of components as inputs and outputs a larger circuit. Finally,
it doesn't matter if several PCB are connected together and then the
chips are placed on them, or if the chips are placed first and the PCBs
are connected later; the resulting overall circuit remains the same.

\begin{figure}[h!]
\parpic[r]{\hspace{-1cm}\raisebox{-18mm}{$\input figs/FlipFlop.pstex_t$}}
\caption{
  The J-K flip flop, a very basic memory cell, is an electronic circuit
  that can be realized using 9 components --- two triple-input ``and''
  gates, two standard ``nor'' gates, and 5 ``junctions'' in which 3 wires
  connect (many engineers would not consider the junctions to be real
  components, but we do). Note that the ``crossing'' in the middle of
  the figure is merely a projection artifact and does not indicate an
  electrical connection, and that electronically speaking, we need not
  specify how this crossing may be implemented in $\bbR^3$. The J-K flip
  flop has 5 external connections (labelled J, K, CP, Q, and Q') and hence
  in the circuit algebra of computer parts, it lives in $C_5$. In the
  directed circuit algebra of computer parts it would be in $C_{3,2}$
  as it has 3 incoming wires (J, CP, and K) and two outgoing wires
  (Q and Q').
} \label{fig:FlipFlop}
\end{figure}

\subsection{Proof of Proposition \ref{prop:sKTGgens}}
\label{subsec:sKTGgensProof}

We are going to ignore strand orientations throughout this proof for 
simplicity. This is not an issue as orientation switches are allowed in
$\sKTG$ without restriction. We are also going to omit vertex signs from
the pictures given the pictorial convention stated in Section \ref{subsec:KTG}.

We need to prove that any $\sKTG$ (call it $G$) can be built from the
generators listed in the statement of the proposition, using $\sKTG$ operations. To show this, consider a Morse 
drawing of $G$, that is, a planar projection of $G$ with a height function
so that all singularities along the strands are Morse and so that every
``feature'' of the projection (local minima and maxima, crossings and
vertices) occurs at a different height.

The idea in short is to decompose $G$ into levels of this Morse drawing where at each level only one ``feature'' occurs. The levels themselves
are not $\sKTG$'s, but we show that the composition of the levels can be achieved by composing their ``closed-up'' $\sKTG$ versions followed
by some unzips. Each feature gives rise to a generator by ``closing up'' extra ends at its top and bottom. We then show that we can
construct each level using the generators and the tangle insert operation.

So let us decompose $G$ into a composition of trivalent tangles (``levels''), each of which has one ``feature'' and (possibly) some 
straight vertical strands. Note that by isotopy we can make sure that every level has strands ending at both its bottom and top,
except for the first or the last level in the case of 1-tangles. An example of level decomposition is shown in the figure below. Note that the levels are generally not
elements of $\sKTG$ (have too many ends). However, we can turn each of them into a $(1,1)$-tangle (or a 1-tangle in case of the aforementioned
top first or last levels) by ``closing up'' their tops and bottoms
by arbitrary trees. In the example below we show this for one level of the Morse-drawn $\sKTG$ containing a crossing and two vertical strands.
\begin{center}
\input{figs/MorseTangle.pstex_t}
\end{center}

Now we can compose the $\sKTG$'s obtained from closing up each level. Each tree
that we used to close up the tops and bottoms of levels determines a ``parenthesization'' of the strand endings. If these parenthesizations
match on the top of each level with the bottom of the next, then we can recreate tangle composition of the levels by composing their closed
versions followed by a number of unzips performed on the connecting trees. This is illustrated in the example below, for two consecutive levels
of the $\sKTG$ of the previous example.
\begin{center}
 \input{figs/CombineLevels.pstex_t}
\end{center}

If the trees used to close up consecutive levels correspond to different parenthesizations, then we can use insertion of the left and right associators 
(the 5th and 6th pictures of the list of generators in the statement of the theorem) to change one parenthesization to match the other. This is 
illustrated in the figure below.
\begin{center}
 \input{figs/Reassociate.pstex_t}
\end{center}

So far we have shown that $G$ can be assembled from closed versions of the levels in its Morse drawing. The closed versions of the levels of $G$ are
simpler $\sKTG$'s, and it remains to show that these can be obtained from the generators using $\sKTG$ operations. 

\parpic[r]{\input{figs/SrtandsToBubble.pstex_t}}
Let us examine what each level
might look like. First of all, in the absence of any ``features'' a level might be a single strand, in which case it is the first generator itself. 
Two parallel strands when closed up become the ``bubble'', as shown on the right.

Now suppose that a level consists of $n$ parallel strands, and that the trees used to close it up on the top and bottom are horizontal mirror images 
of each other, as shown below (if not, then this can be achieved by associator insertions and unzips). We want to show that this $\sKTG$ can
be obtained from the generators using $\sKTG$ operations. Indeed, this can be achieved by repeatedly inserting bubbles into a bubble, as shown:
\begin{center}
\input{figs/BubbleInsertions.pstex_t}
\end{center}

A level consisting of a single crossing becomes a left or right twist
when closed up (depending on the sign of the crossing). Similarly, a
single vertex becomes a bubble. A single minimum or maximum becomes a 
noose or a balloon, respectively.

It remains to see that the $\sKTG$'s obtained when closing up simple features accompanied by more through strands can be built from the generators.
A minimum accompanied by an extra strand gives rise to the $\sKTG$ obtained by sticking a noose onto a vertical strand (similarly, a balloon for a maximum). 
In the case of all the other simple features and for minima and maxima accompanied by more strands, we inserting the already generated elements into 
nested bubbles (bubbles inserted into bubbles), as in the example shown below.
This completes the proof.
\begin{center}
 \input{figs/MoreStrands.pstex_t}
\end{center}
\qed

%% file: figs/Reassociate.pstex_t
\begin{picture}(0,0)%
\includegraphics{figs/Reassociate.pstex}%
\end{picture}%
\setlength{\unitlength}{4144sp}%
\begingroup\makeatletter\ifx\SetFigFont\undefined%
\gdef\SetFigFont#1#2#3#4#5{%
  \reset@font\fontsize{#1}{#2pt}%
  \fontfamily{#3}\fontseries{#4}\fontshape{#5}%
  \selectfont}%
\fi\endgroup%
\begin{picture}(6729,2364)(7144,-2863)
\put(12016,-1636){\makebox(0,0)[lb]{\smash{{\SetFigFont{12}{14.4}{\rmdefault}{\mddefault}{\updefault}{\color[rgb]{0,0,0}unzips}%
}}}}
\end{picture}%

%% file: glossary.tex
\draftcut \section{Glossary of notation} \label{sec:glossary} Greek letters, then Latin, then symbols:

\noindent
{\small \begin{multicols}{2}
\begin{list}{}{
  \renewcommand{\makelabel}[1]{#1\hfil}
}

\item[{$\delta$}] Satoh's tube map\hfill\ \ref{subsec:TangleTopology}
\item[{$\Delta$}] co-product\hfill\ \ref{subsec:ATSpaces}
\item[{$\iota$}] inclusion $\attr_n\to\calP^w(\uparrow_n)$\hfill\ 
  \ref{subsec:ATSpaces}
\item[{$\nu$}] the invariant of the unknot\hfill\ \ref{subsec:KTG}
\item[{$\pi$}] the projection $\calP^w(\uparrow_n)\to\fraka_n\oplus\tder_n$
  \hfill\ \ref{subsec:ATSpaces}
\item[{$\phi$}] log of an associator\hfill\ \ref{subsec:KTG}
\item[{$\Phi$}] an associator\hfill\ \ref{subsec:KTG}
\item[{$\psi_\beta$}] ``operations''\hfill\ \ref{subsec:AlgebraicStructures}

\item
\item[{$\fraka_n$}] $n$-dimensional Abelian Lie algebra\hfill\ 
  \ref{subsec:ATSpaces}
\item[{$\calA$}] a candidate associated graded structure\hfill\ \ref{subsec:Expansions}
\item[{$\calA^{sv}$}] $\calD^v$ mod 6T, RI\hfill\ \ref{subsec:vw-tangles}
\item[{$\calA^{sw}$}] $\calD^w$ mod $\aft$, TC, RI\hfill\ 
  \ref{subsec:vw-tangles}
\item[{$\calA^{sw}$}] $\grad\wTF$\hfill\ \ref{subsec:fgrad}
\item[{$\calA^{sw}$}] $\grad\wTFo$\hfill\ \ref{subsec:OriFoams}
\item[{$\calA^{(s)w}$}] $\calA^{w}$ and/or $\calA^{sw}$\hfill\ 
  \ref{subsec:fgrad}
\item[{$\calA^u$}] chord diagrams mod rels for KTGs\hfill\ \ref{subsec:KTG}
\item[{$\calA^v$}] $\calD^v$ mod 6T\hfill\ \ref{subsec:vw-tangles}
\item[{$\calA^w$}] $\calD^w$ mod $\aft$, TC\hfill\ \ref{subsec:vw-tangles}
\item[{$\calA^w$}] $\grad\wTFo$ without RI\hfill\ \ref{subsec:fgrad}
\item[{$\calA^-(\uparrow_n)$}] $\calA^-$ for pure $n$-tangles\hfill\ 
  \ref{subsec:ATSpaces}
\item[{$A_e$}] 1D orientation reversal\hfill\ \ref{subsubsec:wops}
\item[{$\Ass$}] associative words\hfill\ \ref{subsec:ATSpaces}
\item[{$\Ass^+$}] non-empty associative words\hfill\ \ref{subsec:ATSpaces}
\item[{$\calB^w_n$}] $n$-coloured unitrivalent arrow \newline
  diagrams\hfill\ \ref{subsec:ATSpaces}
\item[{$C$}] the invariant of a cap\hfill\ \ref{subsec:wTFExpansion}
\item[{CP}] the Cap-Pull relation\hfill\ \ref{subsubsec:wrels},
  \ref{subsec:fgrad}
\item[{CW}] Cap-Wen relations\hfill\ \ref{subsubsec:wrels}
\item[{$c$}] a chord in $\calA^u$\hfill\ \ref{subsec:KTG}
\item[{$\der$}] Lie-algebra derivations\hfill\ \ref{subsec:ATSpaces}
\item[{$\calD^v$, $\calD^w$}] arrow diagrams for v/w-tangles\hfill\
  \ref{subsec:vw-tangles}
\item[{$\divop$}] the ``divergence''\hfill\ \ref{subsec:ATSpaces}
\item[{$F$}] a map $\calA^w\to\calA^w$\hfill\ \ref{subsec:fgrad}
\item[{$F$}] the main~\cite{AlekseevTorossian:KashiwaraVergne} unknown
  \hfill\ \ref{subsec:EqWithAT}
\item[{FR}] Flip Relations\hfill\ \ref{subsubsec:wrels},
  \ref{subsec:fgrad}
\item[{$\fil$}] a filtered structure\hfill\ \ref{subsec:Expansions}
\item[{$\calI$}] augmentation ideal\hfill\ \ref{subsec:Grad}
\item[{$J$}] a map $\TAut_n\to\exp(\attr_n)$\hfill\ \ref{subsec:ATSpaces}
\item[{$j$}] a map $\TAut_n\to\attr_n$\hfill\ \ref{subsec:ATSpaces}
\item[{KTG}] Knotted Trivalent Graphs\hfill\ \ref{subsec:KTG}
\item[{$\lie_n$}] free Lie algebra\hfill\ \ref{subsec:ATSpaces}
\item[{$l$}] a map $\tder_n\to\calP^w(\uparrow_n)$\hfill\ \ref{subsec:ATSpaces}
\item[{$\calO$}] an ``algebraic structure''\hfill\ 
  \ref{subsec:AlgebraicStructures}
\item[{$\calP^w_n$}] primitives of $\calB^w_n$\hfill\ \ref{subsec:ATSpaces}
\item[{$\calP^-(\uparrow_n)$}] primitives of $\calA^-(\uparrow_n)$\hfill\ 
  \ref{subsec:ATSpaces}
\item[{$\grad$}] associated graded structure \hfill\ \ref{subsec:Grad}
\item[{$R$}] the invariant of a crossing\hfill\ \ref{subsec:wTFExpansion}
\item[{R4}] a Reidemeister move for
  \newline foams/graphs\hfill\ \ref{subsubsec:wrels}
\item[{$\sder$}] special derivations\hfill\ \ref{subsec:sder}
\item[{$\calS$}] the circuit algebra of skeletons\hfill\ 
  \ref{subsec:CircuitAlgebras}
\item[{$\SAut_n$}] the group $\exp(\sder_n)$\hfill\ \ref{subsec:KTG}
\item[{$S_k$}] complete orientation reversal\hfill\ 
  \ref{subsec:UniquenessForTangles}
\item[{$S_e$}] complete orientation reversal\hfill\ \ref{subsubsec:wops}
\item[{$\sKTG$}] signed long KTGs\hfill\ \ref{subsec:KTG}
\item[{TV}] Twisted Vertex relations\hfill\ \ref{subsubsec:wrels}
\item[{$\tder$}] tangential derivations\hfill\ \ref{subsec:ATSpaces}
\item[{$\attr_n$}] cyclic words\hfill\ \ref{subsec:ATSpaces}
\item[{$\attr^s_n$}] cyclic words mod degree 1\hfill\ \ref{subsec:ATSpaces}
\item[{$\TAut_n$}] the group $\exp(\tder_n)$\hfill\ \ref{subsec:ATSpaces}
\item[{$u$}] a map $\tder_n\to\calP^w(\uparrow_n)$\hfill\ \ref{subsec:ATSpaces}
\item[{$u_e$}] strand unzips\hfill\ \ref{subsubsec:wops}
\item[{$\uT$}] u-tangles\hfill\ \ref{subsec:sder}
\item[{$V$, $V^+$}] the invariant of a (positive) vertex\hfill\ 
  \ref{subsec:wTFExpansion}
\item[{$V^-$}] the invariant of a negative vertex\hfill\  
  \ref{subsec:wTFExpansion}
\item[{VI}] Vertex Invariance\hfill\ \ref{subsec:fgrad}
\item[{$\vT$}] v-tangles\hfill\ \ref{subsec:vw-tangles}
\item[{$W$}] $Z(w)$\hfill\ \ref{subsec:wTFExpansion}
\item[{$W^2$}] Wen squared\hfill\ \ref{subsubsec:wrels}
\item[{$w$}] the wen\hfill\ \ref{subsubsec:wTFgens}
\item[$-$] wenjugation \hfill \ref{subsec:OriFoams}
\item[{$\wT$}] w-tangles\hfill\ \ref{subsec:vw-tangles}
\item[{$\wTF$}] w-tangled foams with wens\hfill\ \ref{subsec:wTF}
\item[{$\wTFo$}] orientable w-tangled foams\hfill\ \ref{subsec:OriFoams}
\item[{$Z$}] expansions \hfill\ throughout
\item[{$Z_\calA$}] an $\calA$-expansion\hfill\ \ref{subsec:Expansions}

\item
\item[{4T}] $4T$ relations\hfill\ \ref{subsec:KTG}
\item[{$\uparrow$}] a ``long'' strand\hfill\ throughout
\item[{$\up$}] the quandle operation\hfill\ 
  \ref{subsec:AlgebraicStructures} 
\item[{$*$}] the adjoint on $\calA^w(\uparrow_n)$\hfill\ \ref{subsec:ATSpaces}

\end{list}
\end{multicols}}

%% file: refs.tex
\draftcut

%% file: WKO2.bbl
\begin{thebibliography}{BND1}

\bibitem[BND]{WKO2} D.~Bar-Natan and Z.~Dancso,
    {\em Finite Type Invariants of W-Knotted Objects II: Tangles and
      the Kashiwara-Vergne Problem,}
  Math.\ Ann.\ {\bf 367} (2017) 1517--1586,
  \arXiv{1405.1955v3}.

\bibitem[BNDv2]{WKO2v2} D.~Bar-Natan and Z.~Dancso,
  \href
    {http://drorbn.net/AcademicPensieve/Projects/WKO2}
    {{\em Finite Type Invariants of W-Knotted Objects II: Tangles and
      the Kashiwara-Vergne Problem,}}
  \url{http://drorbn.net/WKO2} and \arXiv{1405.1955}.
  
\bibitem[BDS]{BDS} D.~Bar-Natan, Z.~Dancso, N.~Scherich,
{\em Ribbon 2-knots, $1+1=2$, and Duflo's Theorem for Arbitrary Lie Algebras}
Alg.\ Geom.\ Topol.\  {\bf 20}(7) 3733--3760. \arXiv{1811.08558}

  
\bibitem[DHR]{DHR} Z.~Dancso, I.~Halacheva, M.~Robertson, 
{\em A Topological Characterisation of the Kashiwara--Vergne Groups}
Trans.\ Amer.\ Math.\ Soc. {\bf 376}(5), 3265--3317. \arXiv{2106.02373}

\end{thebibliography}

\begin{thebibliography}{CCFM}

\bibitem[AM]{AlekseevMeinrenken:KV} A.~Alekseev and
  E.~Meinrenken,
  {\em On the Kashiwara-Vergne conjecture,}
   Inventiones Mathematicae, {\bf 164} (2006) 615--634, \arXiv{0506499}.

\bibitem[AT]{AlekseevTorossian:KashiwaraVergne} A.~Alekseev and
  C.~Torossian,
  {\em The Kashiwara-Vergne conjecture and Drinfel'd's associators,}
   Annals of Mathematics, {\bf 175} (2012) 415--463, \arXiv{0802.4300}.

\bibitem[AET]{AlekseevEnriquezTorossian:ExplicitSolutions} A.~Alekseev, 
B.~Enriquez, and C.~Torossian,
   {\em Drinfel'd's associators, braid groups and an explicit solution of 
the Kashiwara-Vergne equations,}
    Publications Math\'ematiques de L'IH\'ES, {\bf 112-1} (2010) 143--189, 
\arXiv{0903.4067}.







\bibitem[BN1]{Bar-Natan:OnVassiliev} D.~Bar-Natan,
  \href{http://www.math.toronto.edu/~drorbn/LOP.html#OnVassiliev}{{\em
    On the Vassiliev knot invariants,}}
  Topology {\bf 34} (1995) 423--472.


\bibitem[BN2]{Bar-Natan:NAT} D.~Bar-Natan,
  \href{http://www.math.toronto.edu/~drorbn/LOP.html#NAT}{{\em
    Non-associative tangles,}}
  in {\em Geometric topology} (proceedings of the Georgia international
  topology conference), (W.~H.~Kazez, ed.), 139--183, Amer.{} Math.{}
  Soc.{} and International Press, Providence, 1997.



\bibitem[BN3]{Bar-Natan:Associators} D.~Bar-Natan,
  \href{http://www.math.toronto.edu/~drorbn/LOP.html#Associators}{{\em On
  Associators and the Grothendieck-Teichmuller Group I,}}
  Selecta Mathematica, New Series {\bf 4} (1998) 183--212.


\bibitem[BN4]{Bar-Natan:AKT-CFA} D.~Bar-Natan,
  {\em Algebraic Knot Theory --- A Call for Action,}
  web document, 2006,
  \url{http://www.math.toronto.edu/~drorbn/papers/AKT-CFA.html}.

\bibitem[BND]{Bar-NatanDancso:KTG} D.~Bar-Natan and Z.~Dancso,
  {\em Homomorphic expansions for knotted trivalent graphs},
  Journal of Knot Theory and its Ramifications Vol.\ {\bf 22}, No.\ 1 (2013)
  \arXiv{1103.1896}

\bibitem[BGRT]{Bar-NatanGaroufalidisRozanskyThurston:WheelsWheeling}
   D.~Bar-Natan, S.~Garoufalidis, L.~Rozansky and D.~P.~Thurston,
   {\em Wheels, wheeling, and the Kontsevich integral of the unknot,}
   Israel Journal of Mathematics {\bf 119} (2000) 217--237,
   \arXiv{q-alg/9703025}.

\bibitem[BHLR]{Bar-NatanHalachevaLeungRoukema:v-Dims} D.~Bar-Natan,
  I.~Halacheva, L.~Leung, and F.~Roukema,
  \href{http://www.math.toronto.edu/~drorbn/papers/v-Dims}{{\em Some
  Dimensions of Spaces of Finite Type Invariants of Virtual Knots,}}
  submitted.

\bibitem[BLT]{Bar-NatanLeThurston:TwoApplications} D.~Bar-Natan,
  T.~Q.~T.~Le, and D.~P.~Thurston,
  {\em Two applications of elementary knot theory to Lie algebras and
    Vassiliev invariants,}
 Geometry and Topology {\bf 7-1} (2003) 1--31, \arXiv{math.QA/0204311}.


\bibitem[BP]{BerceanuPapadima:BraidPermutation} B.~Berceanu and 
S.~Papadima,
   {\em Universal Representations of Braid and Braid-Permutation Groups,}
   J.~of Knot Theory and its Ramifications {\bf 18-7} (2009) 973--983,
   \arXiv{0708.0634}.


\bibitem[BH]{BrendleHatcher:RingsAndWickets} T.~Brendle and A.~Hatcher,
  {\em Configuration Spaces of Rings and Wickets,}
  \arXiv{0805.4354}.

\bibitem[CL]{ChepteaLe:EvenAssociator}
D.~Cheptea and T.~Q.~T.~Le: {\em A TQFT associated to the LMO invariant of three-dimensional manifolds,}
Commun.{} Math.{} Physics {\bf 272} (2007) 601--634

\bibitem[CS]{CarterSaito:KnottedSurfaces} J.~S.~Carter and M.~Saito,
  {\em Knotted surfaces and their diagrams,}
  Mathematical Surveys and Monographs {\bf 55}, American Mathematical
  Society, Providence 1998.



\bibitem[D]{Dahm:GeneralBraid} D.~M.~Dahm,
{\em A generalization of braid theory}, PhD Thesis, Princeton university, 1962.


\bibitem[Da]{Dancso:KIforKTG} Z.~Dancso,
  {\em On a Kontsevich Integral for Knotted Trivalent Graphs,}
  in {\em Algebraic and Geometric Topology} {\bf 10} (2010) 1317--1365,
  \arXiv{0811.4615}.

\bibitem[Dr1]{Drinfeld:QuantumGroups} V.~G.~Drinfel'd,
  {\em Quantum Groups,}
  in {\em Proceedings of the International Congress of Mathematicians,}
  798--820, Berkeley, 1986.

\bibitem[Dr2]{Drinfeld:QuasiHopf} V.~G.~Drinfel'd,
  {\em Quasi-Hopf Algebras,}
  Leningrad Math.{} J.{} {\bf 1} (1990) 1419--1457.

\bibitem[Dr3]{Drinfeld:GalQQ} V.~G.~Drinfel'd,
  {\em On Quasitriangular Quasi-Hopf Algebras and a Group Closely
   Connected with $\text{Gal}(\bar{\bbQ}/\bbQ)$,}
  Leningrad Math.{} J.{} {\bf 2} (1991) 829--860.



\bibitem[EK]{EtingofKazhdan:BialgebrasI} P.~Etingof and D.~Kazhdan,
  {\em Quantization of Lie Bialgebras, I,}
  Selecta Mathematica, New Series {\bf 2} (1996) 1--41,
  \arXiv{q-alg/9506005}.

\bibitem[FRR]{FennRimanyiRourke:BraidPermutation} R.~Fenn, R.~Rimanyi, and
  C.~Rourke,
  \href{http://msp.warwick.ac.uk/~cpr/}{{\em The Braid-Permutation Group,}}
  Topology {\bf 36} (1997) 123--135.

\bibitem[Gol]{Goldsmith:MotionGroups} D.~L.~Goldsmith,
  \href{http://projecteuclid.org/DPubS?service=UI&version=1.0&verb=Display&handle=euclid.mmj/1029002454}{{\em
  The Theory of Motion Groups,}}
  Mich.{} Math.{} J.{} {\bf 28-1} (1981) 3--17.










\bibitem[Jon]{Jones:PlanarAlgebrasI} V.~Jones,
  {\em Planar algebras, I,}
  New Zealand Journal of Mathematics, to appear, \arXiv{math.QA/9909027}.



\bibitem[KV]{KashiwaraVergne:Conjecture} M.~Kashiwara and M.~Vergne,
  \href{http://www.springerlink.com/content/v73014gx14084624/}{{\em The
    Campbell-Hausdorff Formula and Invariant Hyperfunctions,}}
  Invent.{} Math.{} {\bf 47} (1978) 249--272.


\bibitem[Ka]{Kauffman:VirtualKnotTheory} L.~H.~Kauffman,
  {\em Virtual Knot Theory,}
  European J.{} Comb.{} {\bf 20} (1999) 663--690, \arXiv{math.GT/9811028}.





\bibitem[Kup]{Kuperberg:VirtualLink} G.~Kuperberg,
  {\em What is a Virtual Link?,}
  Algebr.{} Geom.{} Topol.{} {\bf 3} (2003) 587--591,
  \arXiv{math.GT/0208039}.




\bibitem[LM]{LeMurakami:Universal} T.~Q.~T.~Le and J.~Murakami,
  {\em The universal Vassiliev-Kontsevich invariant for framed oriented links,}
  Compositio Math.{} {\bf 102} (1996) 41--64, \arXiv{hep-th/9401016}.




\bibitem[Lei]{Leinster:Higher} Tom Leinster,
  \href{http://www.maths.gla.ac.uk/~tl/book.html}
  {{\em Higher Operads, Higher Categories,}}
  London Mathematical Society Lecture Note Series {\bf 298}, Cambridge
  University Press, ISBN 0-521-53215-9, \arXiv{math.CT/0305049}.




\bibitem[Lev]{Levine:Addendum} J. ~Levine,
{\em Addendum and Correction to: ``Homology Cylinders: an Enlargement of the Mapping Class Group,}
Alg.{} Geom.{} Top.{} {\bf 2} (2002), 1197--1204, \arXiv{math.GT/0207290}



\bibitem[Lod]{Loday:LeibnizAlg} J-L. ~Loday,
{{\em Une version non commutative des algebres de Lie: des algebres
de Leibniz,}}
Enseign.{} math.{} (2) {\bf 39} (3-4): 269--293.


\bibitem[Mc]{McCool:BasisConjugating} J.~McCool,
  {\em On Basis-Conjugating Automorphisms of Free Groups,}
  Can.{} J.{} Math.{} {\bf 38-6} (1986) 1525--1529.

\bibitem[MM]{MilnorMoore:Hopf} J. Milnor and J. Moore,
  {\em On the structure of Hopf algebras,}
  Annals of Math.{} {\bf 81} (1965) 211--264.

\bibitem[MO]{MurakamiOhtsuki:KTGs}
J. Murakami, and T. Ohtsuki, 
{\em Topological quantum field theory for the universal
quantum invariant,} Communications in Mathematical Physics {\bf 188} 3
(1997) 501--520.






\bibitem[Sa]{Satoh:RibbonTorusKnots} S.~Satoh,
  {\em Virtual Knot Presentations of Ribbon Torus Knots,}
  J.~of Knot Theory and its Ramifications {\bf 9-4} (2000) 531--542.





\bibitem[WKO0]{WKO} D.~Bar-Natan and Z.~Dancso,
  {\em Finite Type Invariants of W-Knotted Objects: From Alexander to
    Kashiwara and Vergne,}
  earlier web version of the first two papers of this series
  in one. Paper, videos (wClips) and related files at
  \url{http://www.math.toronto.edu/~drorbn/papers/WKO/}. The
  \arXiv{1309.7155} edition may be older.

\bibitem[WKO1]{Bar-NatanDancso:WKO1} D.~Bar-Natan and Z.~Dancso,
  {\em Finite Type Invariants of W-Knotted Objects I: Braids,
  Knots and the Alexander Polynomial,}
  \url{http://www.math.toronto.edu/drorbn/LOP.html#WKO1}, \arXiv{1405.1956}.
  
\bibitem[WKO2C]{Bar-NatanDancso:WKO2Corr} D.~Bar-Natan and Z.~Dancso,
  {\em Corrigendum to ``Finite Type Invariants of W-Knotted Objects II: Tangles, Foams and the Kashiwara--Vergne problem'',}

\bibitem[WKO3]{Bar-NatanDancso:WKO3} D.~Bar-Natan and Z.~Dancso,
  {\em Finite Type Invariants of W-Knotted Objects III: the Double Tree
    Construction,}
  in preparation.

\end{thebibliography}
